\documentclass[12pt]{amsart}
\usepackage{epsfig}
\usepackage{psfrag}
\usepackage{amssymb}
\usepackage[usenames, dvipsnames]{color}
\usepackage{txfonts}
\usepackage{hyperref}
\usepackage{verbatim}
\usepackage[pagewise]{lineno}
\usepackage{marginnote}
\usepackage{todonotes}
\usepackage[normalem]{ulem}
\usepackage{cancel}
\usepackage{soul}

\theoremstyle{plain}
\newtheorem{theorem}{Theorem}[section]
\newtheorem{lemma}[theorem]{Lemma}

\newtheorem{proposition}[theorem]{Proposition}

\theoremstyle{definition}
\newtheorem{definition}[theorem]{Definition}

\newtheorem{problem}[theorem]{Problem}

\theoremstyle{remark}
\newtheorem{remark}[theorem]{Remark}

\newcommand{\bR}{{\mathbb R}}

\def\lam{\lambda}

\def\th{\theta}

\def\e{\epsilon}

\def\p{\partial}

\def\na{\nabla}
\def\al{\alpha}

\def\O{\Omega}
\def\o{\omega}

\def\be{\begin{equation}}
\def\ee{\end{equation}}
\def\bes{\begin{equation*}}
\def\ees{\end{equation*}}
\def\bali{\begin{aligned}}
\def\eali{\end{aligned}}
\def\al{\begin{aligned}}
\def\eal{\end{aligned}}
\def\erf{\eqref}
\def\lab{\label}

\def\2O{\underline{\O}}

\numberwithin{equation}{section}

\makeatletter
\def\dashint{\operatorname%
{\,\,\text{\bf--}\kern-.98em\DOTSI\intop\ilimits@\!\!}}
\makeatother

\begin{document}


\title[bounded solution]{Bounded solutions to the axially symmetric Navier Stokes equation in a cusp region}

\thanks{}

\author[Q. S. Zhang]{Qi S. Zhang}

\address[Q. S. Zhang]{Department of mathematics, University of California, Riverside, CA 92521, USA}

\email{qizhang@math.ucr.edu}



\begin{abstract}
A domain in $\mathbb{R}^3$ that touches the $x_3$ axis  at one point is found with the following property. For any initial value in a  $C^2$ class, the axially symmetric Navier Stokes equations with Navier slip boundary condition has a finite energy solution that stays bounded for any given time, i.e. no finite time blow up of the fluid velocity occurs. The result seems to be the first case where the Navier-Stokes regularity problem  is solved beyond dimension 2.
\end{abstract}
\maketitle

\tableofcontents

\section{Introduction}

The Navier-Stokes equations (NS) describing the motion of viscous incompressible fluids in a domain $D \subset \bR^3$ is
\be
\lab{nse}
\mu \Delta v -  v \nabla v - \nabla P -\partial_t v =0, \quad div \, v=0, \quad \text{on} \quad D \times (0, \infty)
\ee Here $v$ is the velocity field,  $P$ is the pressure, both of which are the unknowns;
 $\mu>0$ is the viscosity constant, which will be taken as $1$ in this paper. In order to solve the equation, a initial velocity $v_0$ is usually given, together with suitable boundary conditions.
 One can also add a forcing term on the righthand side, then it becomes a nonhomogeneous problem.
Due to Leray \cite{Le2}, if $D=\bR^3$, $v_0 \in L^2(\bR^3)$, the Cauchy problem has a solution in the energy space (c.f. \erf{enorm} below). However, in general it is not known if such solutions stay bounded or regular for all $t>0$. A typical existence theorem for regular solutions always involves a small parameter in the initial condition, as a perturbation of a known regular solution.

In this paper, we will focus on a special case of \eqref{nse}, namely when $v$ and $P$ are independent of the angle in a cylindrical coordinate system $(r,\,\th,\,x_3)$. That is, for $x=(x_1,\,x_2,\,x_3) \in \bR^3$,
$
r=\sqrt{x_1^2+x_2^2}, \quad
\th=\arctan (x_2/x_1),
$ and the basis vectors $e_r,e_\th,e_3$ are:
\[
e_r=(x_1/r, x_2/r, 0),\quad e_\th=(-x_2/r, x_1/r, 0),\quad e_3=(0,0,1).
\]
 In this case, solutions can be written in the form of
\[
v=v_r(r,x_3, t)e_r+v_{\th}(r, x_3, t)e_{\th}+v_3(r, x_3, t)e_3.
\]

Using
 tensor notations and doing vector calculus under the cylindrical system, one finds
\[
\bali
\nabla v &= \p_r v \otimes e_r +\frac{1}{r} \p_\th v \otimes e_\th + \p_{x_3} v \otimes e_3\\
&=(\p_r v_r  e_r + \p_r v_\th e_\th + \p_r v_3 e_3 )\otimes e_r + \frac{1}{r} ( v_r e_\th -
v_\th e_r) \otimes e_\th \\
&\qquad \qquad + (\p_{x_3} v_r  e_r + \p_{x_3} v_\th e_\th + \p_{x_3} v_3 e_3 )\otimes e_3.
\eali
\] It is convenient to denote $e_r \otimes e_r$, $e_r \otimes e_\th$, $e_r \otimes e_3$, ..., $e_3 \otimes e_3$ by the nine single-entry matrices in the standard basis for $3 \times 3$
matrices, going as $J^{11}, J^{12}, J^{13}, J^{21}, ..., J^{33}$. Then $\nabla v$ is given by the 3 by 3 matrix
\be
\lab{nabv}
\na v = \begin{bmatrix} \p_r v_r & - \frac{1}{r} v_\th & \p_{x_3} v_r\\
                          \p_r v_\th &  \frac{1}{r} v_r & \p_{x_3} v_\th \\
                          \p_r v_3 & 0 &  \p_{x_3} v_3
\end{bmatrix},
\ee and $v \na v$ is given by the matrix multiplication $(\na v)  v$ with $v$ being regarded as the column vector $( v_r, v_\th, v_3)^T$.
Therefore, $v_r$, $v_3$ and $v_\theta$ satisfy the axially symmetric Navier-Stokes equations
\be
\begin{aligned}
\lab{eqasns}
\begin{cases}
   \big (\Delta-\frac{1}{r^2} \big )
v_r-(v_r \p_r + v_3 \p_{x_3})v_r+\frac{(v_{\theta})^2}{r}-\partial_r
P-\p_t  v_r=0,\\
   \big   (\Delta-\frac{1}{r^2}  \big
)v_{\theta}-(v_r \p_r + v_3 \p_{x_3} )v_{\theta}-\frac{v_{\theta} v_r}{r}-
\partial_t v_{\theta}=0,\\
 \Delta v_3-(v_r \p_r + v_3 \p_{x_3})v_3-\p_{x_3} P-\p_t v_3=0,\\
 \frac{1}{r} \p_r (rv_r) +\p_{x_3}
v_3=0,
\end{cases}
\end{aligned}
\ee which will be abbreviated as ASNS. Although ASNS looks more complicated than the full 3 dimensional equation,
  a simplification occurs in the 2nd equation where the pressure term drops out.

If the swirl $v_\theta=0$, then it is
known for a while ( O. A. Ladyzhenskaya
\cite{La}, M. R. Uchoviskii and B. I. Yudovich \cite{UY}), that finite energy solutions
to the Cauchy problem of (\ref{eqasns}) in $\mathbb{R}^3$ are smooth for all time $t>0$.
See also the paper by S.
Leonardi, J. Malek, J. Necas,  and  M. Pokorny \cite{LMNP}. By finite energy, we mean the norm
\eqref{enorm} below is finite, i.e., the solution is a Leray-Hopf solution.

In the presence of swirl, it is still not known in general if finite energy solutions
blow up in finite time. By the partial regularity results in \cite{CKN}, possible singularity for suitable weak solutions of ASNS can only occur at the $x_3$ axis. See also \cite{Linf}. The same statement without the word "suitable"  is also true by \cite{BuZh}. Some existence results with a small parameter can be found in \cite{HL} and \cite{HLL}.  Critical or sub-critical regularity conditions can be found in \cite{NP}, \cite{CL} and \cite{JX}.

Despite the difficulty, there has been no lack of research efforts on ASNS. Let us make a brief description of related recent results, starting with the papers by
C.-C. Chen, R. M. Strain, T.-P.Tsai,  and  H.-T. Yau in \cite{CSTY1},
\cite{CSTY2},  G. Koch, N. Nadirashvili, G. Seregin,  and  V.
Sverak in \cite{KNSS}, which appeared around 2008. See also the work by G. Seregin  and
V. Sverak \cite{SS} for a localized version. These authors proved that if
\be
\lab{v<1/r}
|v(x, t)| \le
\frac{C}{r},
\ee then finite energy solutions to the Cauchy problem of ASNS are smooth for all time. Here $C$ is any positive constant.

The proof is based on the fact that the scaling invariant quantity
$\Gamma= r v_\theta$ satisfies the equation
\be
\lab{eqvth}
\Delta \Gamma - b \nabla \Gamma- \frac{2}{r} \p_r
\Gamma-\p_t \Gamma=0,
\ee where $b=v_r e_r + v_3 e_3$. The bound \eqref{v<1/r} says that the equation is essentially scaling invariant and the classical linear regularity theory can be applied after some nontrivial modification.
The above result can be summarized as: type I solutions of ASNS are regular.

Two years later, in the paper \cite{LZ11} by Lei and the author, it was proven that  if
$v_r, v_3$ are in the space of $L^\infty([0, \infty),  BMO^{-1}({\bR}^3))$ and $r v_\th(\cdot, 0) \in L^\infty$, then the solution is regular. Here $BMO$ is the space of functions with bounded mean oscillation, and $BMO^{-1}$ is the space of tempered distributions which can be written as partial derivatives of BMO functions. Well-posedness and other properties of solutions to NS have been studied by Koch-Tataru \cite{KoTa}.

Recently Seregin and Zhou \cite{SeZh} have relaxed the $L^\infty BMO^{-1}$ assumption further to $L^\infty \dot{B}^{-1}_{\infty, \infty}$ assumption.
Let us recall $\dot{B}^{-1}_{\infty, \infty}$ is the Besov space consisted of tempered distributions $f$ such that the norm
\[
\Vert f \Vert_{\dot{B}^{-1}_{\infty, \infty}} = \sup_{t>0} t^{1/2} \sup_{x} \left| \int_{\bR^3} G(x, t, y) f(y) dy \right|
\]is finite. Here $G(x, t, y)=(1/(4\pi t)^{3/2}) \exp(-|x-y|^2/(4t))$ is the standard heat kernel on $\bR^3$.
  In these papers the regularity conditions
are critical and hence are scaling invariant under
standard scaling. Improvements are at most logarithmic so far. See the paper by X.H. Pan
 \cite{Panx}.
In contrast, the energy bound scales as
$-1/2$. So even with axial symmetry, there seems to be a finite scaling gap which makes the ASNS supercritical, just like the full equations.

However in a recent paper \cite{LZ17},  Lei and the author made the following observation.

{\it The vortex stretching term of the ASNS is critical after a  suitable change of dependent variables.}

So the aforementioned scaling gap is $0$. This observation has the effect of making ASNS looks less formidable than the full 3 dimensional one which has a positive scaling gap. Nevertheless all major open problems for the latter are still open for the former.

The observation is based on the study of the equations for the functions $\O=\o_\th/r$ and $J=\o_r/r$ where $\o_\th$ and $\o_r$ are the angular and radial components of the vorticity.
The function $\O$ has been around for longtime \cite{UY}. But $J$ was introduced in the  recent paper by H. Chen-D.Y.Fang-T. Zhang
\cite{CFZ}. By carrying out an energy estimate for these equations, they
proved the following result. Let $v$ be a Leray-Hopf solution to
the Cauchy problem of ASNS with initial data $v_0  \in H^2$
and $\|r  (v_0 \cdot e_\th) \|_{L^\infty} < \infty$. If
\[
  |v_\theta(x, t)| \le C/r^{1-\e},
\] for all $x$ and $t>0$,
then $v$ is regular everywhere.  Here $\e>0$ and $C$ are positive
constants.

The main result in \cite{LZ17} includes the following statement.
Let $\delta_0 \in (0, \frac{1}{2})$ and $C_1 > 1$.
If
\begin{equation}\label{CD}
\sup_{0 \leq t < T}|r v_\theta (r, x_3, t)| \leq C_1|\ln r|^{- 2},\ \ r \leq \delta_0,
\end{equation}
then above $v$ is regular globally in time.
Note that a priori we have $|r v_\theta (r, x_3, t)| \leq C$ by the maximal principle applied on equation
\erf{eqvth}. So there is still a gap of logarithmic nature from regularity.

After \cite{LZ17} was posted on the arxiv, in the paper by Dongyi Wei \cite{Weid}, the power in the log term has
been improved
to $-3/2$. Namely,  if, for some $\delta_0 \in (0, 1/2)$,
\begin{equation}\label{CDwei}
\sup_{0 \leq t < T}|r v_\theta (r, x_3, t)| \leq C_1|\ln r|^{- 3/2},\ \ r \leq \delta_0,
\end{equation} then the above $v$ is regular.

In this paper, for a special class of bounded domains with Navier slip boundary condition, we manage to solve the regularity problem to ASNS. By finite energy, we mean the solutions are in the energy space $\mathbf{E}=L^2_tW^{1, 2}_x \cap L^\infty_t L^2_x$.  Here and throughout, the norm in $\mathbf{E}$ for a function $v$ on $D \times [0, T]$ is taken as
\be
\lab{enorm}
\Vert v \Vert^2_\mathbf{E} = \int^T_0 \int_D |\nabla v|^2 dxdt + \sup_{t \in [0, T]} \int_D |v(x, t)|^2 dx.
\ee  The function $v$ can be  vector or scalar valued, depending on the context,  and $T>0$.

These domains, which touch the $x_3$ axis at one point, are the union of a sequence of rectangles in the $r-x_3$ plane . The ratio between the height and width of the rectangles decreases slowly to $0$ when they approaches the $x_3$ axis.  The side view of the domain is a wedge like region with infinitely many terraces.
More specifically, the domains are given in

\begin{definition} (see also the figure below.)
\be
\lab{dodm}
\al
D_*&= \bigcup^\infty_{m=1}  D_m, \quad \text{with} \qquad  D_m  \equiv \bigcup^m _{j=1}  S_j,  \\
 S_j  &\equiv   \{  (r, x_3)  \, | \,  2^{-j} \le r < 2^{-(j-1)}, \,    0< x_3 < 2^{-\beta (j-1)} \}
\eal
\ee where $\beta \in (1, 1.1]$ is any fixed number.
\end{definition}

\begin{center}
\includegraphics[scale=0.3]{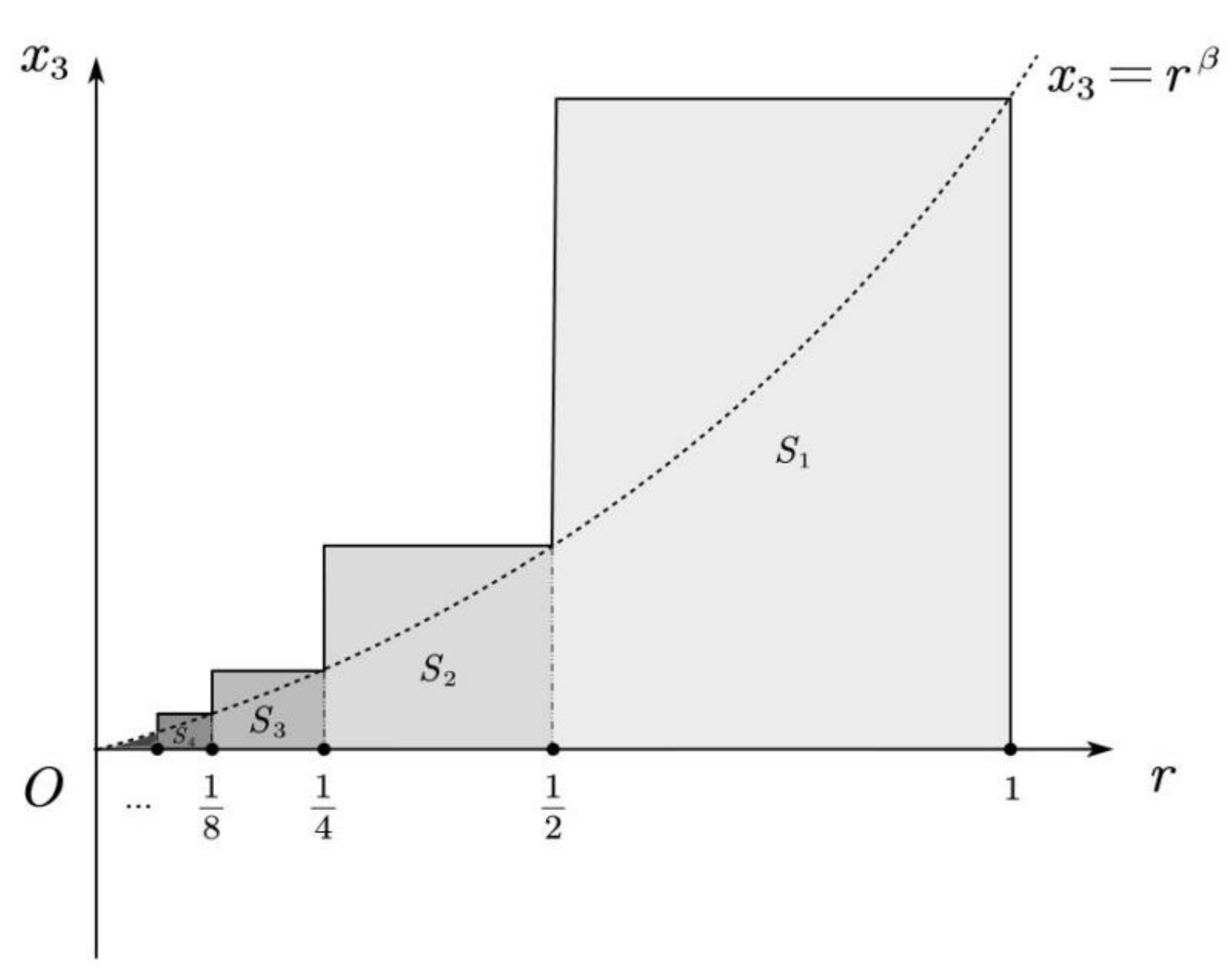}
\end{center}

Recall that the Navier \cite{Nav} slip boundary condition reads
\be
\lab{nvslbc}
(\mathbb{S}(v) n)_{tan} = 0, \qquad v \cdot n = 0, \quad \text{on} \quad \p D_*.
\ee  Here the strain tensor $\mathbb{S}(v)= [\na v + (\na v)^T]/2$ which, by \erf{nabv}, is given by
\be
\lab{Sv}
 \mathbb{S}(v) = \begin{bmatrix} \p_r v_r & \frac{1}{2} (\p_r v_\th - \frac{1}{r} v_\th) & \frac{1}{2}(\p_{x_3} v_r + \p_r v_3) \\
                          \frac{1}{2} (\p_r v_\th - \frac{1}{r} v_\th) &  \frac{1}{r} v_r & \frac{1}{2} \p_{x_3} v_\th \\
                           \frac{1}{2}(\p_{x_3} v_r + \p_r v_3) & \frac{1}{2} \p_{x_3} v_\th &  \p_{x_3} v_3
\end{bmatrix}.
\ee

 In the above $n$ is the unit outward normal on the smooth part of $\p D_*$ and the $(\mathbb{S}(v) n)_{tan}$ stands for the tangential component of the vector $\mathbb{S}(v) n$. Observe that the boundary
$\p D_*$ can be written as the union of horizontal and vertical parts, which are denoted by $\p^H D_*$ and $\p^V D_*$ respectively.  Namely,
\be
\lab{d=dhdv}
\p D_* = \p^H D_* \cup \p^V D_*
\ee  From  \eqref{Sv} and \eqref{nvslbc}, one sees that the Navier slip boundary condition can be expressed explicitly as
\be
\lab{nvbc2}
\bali
\p_{x_3} v_r=\p_{x_3} v_\th=0, \quad v_3=0, \quad \text{on} \quad \p^H D_*;\\
\p_r v_\th=\frac{v_\th}{r}, \, \p_r v_{x_3}=0, v_r=0, \quad \text{on} \quad \p^V D_*.
\eali
\ee

  Due to the cusp shape of $D_*$  near $x_3$ axis and the boundary condition, one needs to define the space for initial values carefully. We pick the space $C^2_{nb}(\overline{D_*}) \subset C^2(D_*)$  for the initial values, which is defined as
\[
\al
C^2_{nb}(\overline{D_*}) &\equiv \{ f \in C^2(D_*) \, | \,  f \, \text{is the limit under the $C^2$ norm of functions
} \\
&\qquad
f_m \in C^2(\overline{D_m}) \quad \text{ satisfying the Navier slip condition on }
D_m, \quad div \, f_m=0 \}.
\eal
\]The above convergence means $\Vert f \, \chi_{D_m} - f_m \Vert_{C^2} \to 0$ as $m \to \infty$.
The main result states that finite energy solutions of ASNS in $D_*$ are bounded in finite time, namely, no finite time blow up occurs.

\begin{theorem}
\label{thbdedsol}
For any initial value in the $C^2_{nb}(\overline{D_*})$ space, the axially symmetric Navier Stokes equations \erf{eqasns} with Navier slip boundary condition on $D_*$ has a finite energy solution such that the velocity $v$ stays bounded for any given time.
\end{theorem}

The result appears to have two implications. First, it seems to be the first case where the regularity problem for finite energy solutions is  solved beyond dimension 2. In the literature ASNS are sometimes considered as  $2\frac{1}{2}$ dimensional equations. Recall that Ladyzhenskaya \cite{La2} (Theorem 6, Chapter 6, Sec.3) already knew that if the domain is bounded away from the symmetric axis,  finite energy solutions to ASNS are regular. But no progress was made afterwards if the domain touches the axis. The second pertains recent interesting non-uniqueness results by Buckmaster and Vicol \cite{BV} on some very weak (infinite energy) solutions of the Navier Stokes equations using the technique of convex integration of
De Lellis and L. Szekelyhidi \cite{DS}. In the paper \cite{BV} the authors also stated a conjecture on non-uniquess and hence non-regularity of some finite energy solutions to the Navier Stokes equations.
See also \cite{ChLu} by Cheskidov and Luo and \cite{HH} by Hou and Huang. Although the above theorem is not exactly a counter example to the conjecture due to the boundary condition and symmetry, it seems to indicate an additional obstacle to the construction of a non-regular solution with finite energy.

Let us briefly describe the proof of the above theorem.
We will first construct solutions $v^{(m)}$ of ASNS in the intermediate regions $D_m$
 with Navier slip boundary condition. Since $D_m$ is bounded away from the $x_3$ axis, the solutions are bounded. The key step is to prove that these bound is independent of $m$, after which one can extract a convergent subsequence to obtain the claimed bounded solution. The usual path to do so is to find a uniform $L^2_tW^{2, 2}_x$ bound for the solutions. However due to the presence of nonconvex corners in $D_m$, such a bound is not expected to be true. See \cite{MS} for an example for 2 dimensional Stokes system. Our observation is that a partial $L^2_tW^{2, 2}$ bound holds, i.e. some components of the $\na^2 v^{(m)}$ are in $L^2_tL^2_x$ space. This will allow us to start an energy estimate on the equations for $\o^{(m)}_\th/r$ and $\o^{(m)}_r/r$. Here $\o^{(m)}_\th$ and $\o^{(m)}_r$ are the angular and radial components of the vorticity $curl \, v^{(m)}$. The cusp shape of $D_m$ and boundary conditions are used together to show that the vortex stretching terms are dominated by the diffusion near the $x_3$ axis. This will lead to a $L^\infty_tL^2_x$ bounds for $\o^{(m)}_\th/r$ and $\o^{(m)}_r/r$. These bounds are used to prove the $L^\infty$ bound for $v^{(m)}$ via the Biot-Savart law. This is easy in the full space case. But the presence of the cusp makes the proof harder. One reason is the standard Sobolev inequality may fail. To proceed, we employ some new observations on the velocity such as the line integral of $v_r$ is $0$ in the $x_3$ direction, and interactions between components of $v$ and $\o$  and boundary values.
 In the theorem, the parameter $\beta$ can be made a slightly larger than $1.1$ but it remains to be seen if it can be chosen as $1$. If $\beta=1$, then the domain $D_*$, after self similar extension to $r=\infty$, resembles the exterior of a cone  globally. One can prove, by the same method, that no finite time blowup occurs under an extra condition $|r v_\theta(\cdot, 0)| \le \epsilon$ for some $\epsilon>0$. Note that no extra assumption is made on $v_r$ or $v_3$.

We end the introduction by listing a number of notations and conventions to be used throughout, which are more or less standard.
The velocity field is usually called $v$ and the vorticity $\nabla \times v$ is called $\o$. We use subscripts to denote their components in coordinates. For example $v_r = v \cdot e_r$, $\o_\th = \o \cdot e_\th$, $\o_3= \o \cdot e_3$; $\O = \o_\th/r$ and $J=\o_r/r$.  We write $b=v_r e_r + v_3 e_3$. $L^p(D)$, $p \ge 1$, denotes the usual Lebesgue space on a domain D which may be a spatial, temporal or space-time domain. Let $X$ be a Banach space defined for functions on $D \subset \bR^3$. $L^p(0, T; X)$ is the Banach space of space-time functions $f$ on the space time domain $D \times [0, T]$ with the norm $\left(\int^T_0 \Vert f(\cdot, t) \Vert^p_X dt\right)^{1/p}$. If no confusion arises, we will also use $L^p X$ to abbreviate $L^p(0, T; X)$. Sometimes we will also use $L^p_xL^q_t$ or $L^q_tL^p_x$ to denote the mixed $p, q$ norm in space time. Let $D \subset \bR^3$ be an open domain, then $H^1(D)=W^{1, 2}(D)=\{f \, | \,  f, |\nabla f| \in L^2(D) \}$ and $H^2(D)=W^{2, 2}(D)=\{f \, |\,  f, |\nabla f|,  |\nabla^2 f|\in L^2(D) \}$, the standard Sobolev spaces on D. Also, interchangeable notations $div \, v = \nabla \cdot v$, $v \nabla v =\sum v_i \p_{x^i} v = v \cdot \nabla v$ will be used. If there is no confusion, the vertical variable $x_3$ may be replaced with $z$. Also $B(x, r)$ denotes the ball of radius $r$ centered at $x$ in a Euclidean space; and $B_X(f, r)$ denotes the open ball in a normed space $X$, centered at $f \in X$ with radius $r$. If $s$ is a number, then $s^-$ means any number which is close but strictly less than $s$.
We use $C$ with or without index to denote a generic constant which may change from line to line and
$\overline{C}$, $C_*$ etc, $\alpha, \beta$ etc to denote important constants which may depend on relevant functions such as initial values, solutions.

\section{Finite energy solutions with partial $L^2_t W^{2, 2}_x$ bounds on approximating domains}

In this section, we study the existence of finite energy solutions to the following initial boundary value problem on $D_m \times (0, T]$, $T>0$. We will also prove some  extra regularity result, which, although not uniform in $m$, will be needed in the next section.

\begin{problem}
\lab{wtdm}

 Find a function $v=v_r e_r + v_\th e_\th + v_3 e_3$ in the energy space $\mathbf{E}= L^2_tW^{1, 2}_x \cap L^\infty_t L^2_x$ on the domain $D_m \times (0, T]$,  such that the following  hold.

(a). $v$  satisfies ASNS \eqref{eqasns} in the weak sense on $D_m \times (0, T]$.

(b). $v$ is subject to the Navier slip boundary condition:
\be
\lab{nvbcdm}
\bali
\p_{x_3} v_r=\p_{x_3} v_\th=0, \quad v_3=0, \quad \text{on} \quad \p^H D_m \times [0, T];\\
\p_r v_\th=\frac{v_\th}{r}, \, \p_r v_{x_3}=0, v_r=0, \quad \text{on} \quad \p^V D_m \times [0, T].
\eali
\ee

(c). The initial value $v(x, 0) = v_0$ is in $W^{2, 2}(D_m)$ and it satisfies the Navier slip boundary condition and the divergence free condition.
\end{problem}

We refer the reader to Section 2 of \cite{Be} for a precise definition of (a) and (b) and the exact function spaces involved. In our case, due to the relative simplicity of the boundary and standard interior and boundary regularity results, the boundary condition (b) can also be understood in the pointwise sense, except at the corners.

Our starting point is the basic energy estimate for velocity $v$.
It is known that
a  solution to Problem \ref{wtdm} exists and it satisfies the energy bound
\be
\lab{en1s}
\int_{D_m} |v(x, T)|^2 dx + 4 \int^T_0 \int_{D_m} |\mathbb{S}v(x, t)|^2 dxdt \le \int_{D_m} |v(x, 0)|^2 dx.
\ee Such a result for piecewise smooth domains can be found in the paper \cite{Be} by Benes and \cite{BK} by  Benes and Kucera e.g., where even more complicated mixed boundary conditions are studied.
 When the angles in the polyhedron is not greater than $\pi$,  it is proven in \cite{Be} that these solutions are strong in the sense that the velocity $v$ is in $L^2_tW^{2, 2}_x$ space.   As mentioned, due to the presence of $3\pi/2$ angle in our case, solution $v$ may not enjoy $L^2_tW^{2, 2}_x$ regularity in general. This lack of regularity is an obstacle in proving the main result, which will be circumvented  by the $L^2_tL^2_x$ bound for $|\na \o_\th|$ and $|\na \o_r|$ in the next proposition.
  More general Navier type boundary conditions on smooth domains was treated  earlier by several authors. The interested reader is referred to the relatively recent papers and their references therein: \cite{XX} by Xiao and Xin, \cite{CQ} by G.Q. Chen and Z.M. Qian and \cite{MR} by Masmoudi and Rousset. See also the paper \cite{AS} by Abe and Seregin on ASNS in the exterior of a cylinder containing the $x_3$ axis,  subject to the Navier slip condition.

From \eqref{en1s} and Korn's inequality, one can quickly derive a $L^2_tL^2_x$ bound for $|\na v|$. However this bound will depend on the index $m$ for the domain $D_m$.

The next proposition is the main result in this section. It states that Problem \ref{wtdm} has a  solution in the energy space such that $|\na \o_\th| + | \na \o_r| \in L^2_tL^2_x$, although $|\na^2 v|$ may not be in $L^2_tL^2_x$. In fact we will prove that all but  the component $\p^2_r v_\th$ in the Hessian of $v$ are  in $L^2_tL^2_x$. Note that the above bound  is not uniform in $m$ and it may also depend on the size of the initial value.

\begin{proposition}
\lab{prsolS}
Given a divergence free initial value $v_0 \in W^{2, 2}(D_m)$ which satisfies the Navier slip boundary condition,  Problem  \ref{wtdm} has a solution in the energy space $\mathbf{E}$ on the domain $D_m \times(0, \infty)$. Moreover,
\be
\lab{naol2}
|\na \o_\th|  + | \na \o_r| + |\na^2 v_r| + |\na^2 v_3| \in L^2_tL^2_x
\ee on $D_m \times(0, T]$ for any $T>0$.  In addition, for $\O \equiv \frac{\o_\th}{r}$ and $J \equiv \frac{\o_r}{r}$,
\be
\lab{Onengl1}
\bali
 \int^T_0 \int_{D_m} |\na \O |^2 dxdt + &\frac{1}{2} \int_{D_m} | \O(x, T) |^2 dx - \frac{1}{2} \int_{D_m} |\O(x, 0)|^2 dx\\
 &=\quad 2 \int^T_0 \int_{D_m} \frac{\p_r \O}{r} \O  dxdt
  -\int^T_0 \int_{D_m} \frac{v^2_\th}{r^2} \p_{x_3} \O  dxdt,
 \eali
\ee
\be
\lab{Jnengl1}
\al
 \int^T_0 \int_{D_m}& |\nabla J|^2 dxdt + \frac{1}{2} \int_{D_m} |J(x, T)|^2 dx
 -\frac{1}{2} \int_{D_m} |J(x, 0)|^2 dx \\
&=  \int^T_0 \int_{D_m} \left\{ J \frac{2}{r}\p_r J + v_\th \, \p_r  \frac{v_r}{r} \p_{x_3} J  -  v_\th \, \p_{x_3} \frac{v_r}{r} \p_r  J\right\} dxdt.
\eal
\ee
\end{proposition}

\proof The proof is divided into a number of steps.

{\it Step 1. Outline of the method.}

As mentioned, the existence of solutions of Problem \ref{wtdm} in the energy space is known. The proof is done by the Galerkin method based on the study of the (linear) Stokes equation in non-smooth domains such as $D_m$. So the main task is to prove \eqref{naol2}. However, since the usual proof of existence does not seem to give \erf{naol2}, we will also need an independent proof of the existence, using a method involving the contraction mapping theorem.

Let $v_0= (v_0)_r e_r + (v_0)_\th e_\th +(v_0)_3 e_3$ be the given initial value.
Let $\mathbf{S}$ be the subspace of the energy space $\mathbf{E}$ for functions on $D_m \times [0, T]$, defined by
\[
\mathbf{S} = \{ v=v_r e_r + v_\th e_\th + v_3 e_3 \, | \,  v \in \mathbf{E}, \, div \, b(\cdot, t) = 0, \, a.e. t,  \, \o_\th e_\th=curl (v_r e_r + v_3 e_3) \in L^2_tW^{1, 2}_x(D_m)  \},
\] equipped with the norm
\[
\Vert v \Vert^2_{\mathbf{S}} = \Vert v \Vert^2_{\mathbf{E}} + \int^T_0\int_{D_m} |\na \o_\th|^2 dxdt.
\]
We will prove that there is a time $T>0$ such that Problem \ref{wtdm} has a solution in the closed unit ball $\overline{B_{\mathbf{S}}(v_0, 2)} \subset \mathbf{E}$, which also satisfies \eqref{naol2}. Afterwards, it will be shown that this solution can be extended to all positive time.

 To this end, we pick scalar functions $v_r, v_3$ such that the vector field
$b \equiv v_r e_r + v_3 e_3 \in \overline{B_{\mathbf{E}}(b(0), 1)}$ with $b(0)=(v_0)_r e_r + (v_0)_3 e_3$.
The strategy is to use $b$ as  given data in the equation for $v_\th$. This linearized equation, with suitable boundary condition and $(v_0)_\th$ as the initial value, determines the vector field $v_\th e_\th
\in \overline{B_{\mathbf{E}}((v_0)_\th e_\th, 1)}$, provided that $T$ is sufficiently small.
Using these $v_r, v_\th, v_3$ as given data in the equation for $\O=\o_\th/r$, with $0$ boundary value  and $\o_\th(0)/r$ as initial value, one finds a vector field $\tilde{\o_\th} e_\th=r \o$. Then the Biot-Savart law with a suitable boundary condition determines a vector in $\tilde{b} = \tilde{v_r} e_r + \tilde{v_3} e_3 \in
\overline{B_{\mathbf{S}}(b(0), 1)}$. The correspondence between $b$ and $\tilde{b}$ gives rise to a contraction mapping from $\overline{B_{\mathbf{S}}(b(0), 1)}$ to itself if $T$ is sufficiently small.
Let us denote by $b$ the fixed point. Then $v \equiv b + v_\th e_\th$ is a solution to Problem \ref{wtdm}. The following diagram illustrates this process:
\[
b \Rightarrow v_\th \Rightarrow \O \Rightarrow \tilde{\o}_\th=r \O \Rightarrow \tilde{b}.
\]In order to prove \eqref{naol2}, we need to use crucially the property that $\tilde{\o}_\th$ has $0$ boundary value and some nonstandard uniqueness results for solutions of elliptic and parabolic equations in rectangular domains with non-convex corners, due to the presence of $3\pi/2$ angles on the boundary of $D_m$.
In the following steps, detail of the proof is carried out.

{\it Step 2. constructing $v_\th$.}

Recall from \eqref{eqasns} and the Navier slip boundary condition that $v_\th$ is determined by the following initial boundary value problem
\be
\lab{ibvth}
\begin{cases}
 (\Delta-\frac{1}{r^2}  \big
)v_{\theta}-(v_r \p_r + v_3 \p_{x_3} )v_{\theta}-\frac{v_{\theta} v_r}{r}-
\partial_t v_{\theta}=0, \quad \text{on} \quad D_m \times (0, T];\\
\p_{x_3} v_\th =0 \quad \text{on} \quad \p^H D_m \times (0, T], \quad \p_r v_\th = v_\th/r \quad \text{on} \quad \p^V D_m \times (0, T];\\
v_\th(x, 0)= (v_0)_\th(x), \quad x \in D_m.
\end{cases}
\ee Notice that the boundary condition for $v_\th$ is of mixed Neumann and Robin type, which also has been around for many years. Nevertheless, it is more complicated than the Neumann condition. However  the function
\[
h \equiv v_\th/r
\] satisfies the standard Neumann boundary condition, since $\p_{x_3} h = \p_{x_3} v_\th /r=0$ on $\p^H D_m$, and $\p_r h = (r \p_r v_\th - v_\th)/r^2=0$ on $\p^V D_m$. Therefore, we will first look for solutions in
$\mathbf{E}$ of the following initial boundary value problem
\be
\lab{ibhth}
\begin{cases}
 (\Delta+\frac{2}{r} \p_r \big
)h- b \na h-\frac{ 2 v_r}{r} h-
\partial_t h=0, \quad \text{on} \quad D_m \times (0, T];\\
\p_n h=0,\, \text{on} \quad \p D_m \times (0, T];\\
h(x, 0)=h_0= (v_0)_\th(x)/r, \quad x \in D_m.
\end{cases}
\ee Then  $v_\th = r h$ will give us a solution to \eqref{ibvth}. Here $\p_n h$ is the derivative of with respect to the exterior normal of $\p D_m$, except at the corners. Notice that $D_m$ is bounded away from the $x_3$ axis. So the appearance of $1/r$ will not create singularity for fixed $m$.

In the absence of the lower order terms, the existence and uniqueness of solutions in the energy space
$\mathbf{E}$ for the initial Neumann problem of the heat equation in bounded piecewise smooth domains is a well known classical result. See \cite{LSU} e.g. In the presence of the lower order terms within suitable
$L^p_xL^q_t$ class, the result is still true with routine modification in the proof. A typical proof usually relies on  a priori energy bound and Galerkin method. Since not all lower order terms in the equation of \eqref{ibhth} are in these suitable $L^p_xL^q_t$ class, we will describe the necessary modification. As mentioned in the introduction, elements in $\mathbf{E}$ can either be scalar or vector valued, depending on the context. For \eqref{ibhth}, elements in $\mathbf{E}$ are scalar valued.

 In \eqref{ibhth}, we treat $-2 v_r/r$ as a potential function for the solution $h$. It lies in the space
 $L^2_x L^\infty_t$ by our choice that $b \in \overline{B_{\mathbf{E}}(b(0), 1)} \subset \mathbf{E}$. It is well known that in spatial dimension 3, the standard class for potential functions which guaranty  the existence and H\"older continuity of solutions in $\mathbf{E}$ is $L^p_x L^q_t$ such that $3/p + 2/q <2$. Therefore $-2 v_r/r$ is in the standard class. In contrast, the standard class for the drift term is $L^p_x L^q_t$ with $3/p + 2/q <1$. Hence the vector field $b$ in the drift term $- b \na h$ is not in the standard class. However, since $ div \, b=0$ in the $L^2_tL^2_x$ sense and $v_r=0$ on $\p^V D_m \times (0, T]$ and $v_3 =0$ on $\p^H D_m \times (0, T]$ in  the sense of boundary trace of functions in $L^2_t W^{1, 2}_x(D_m)$, one can prove that the drift term will be integrated out and standard energy estimate can be proven. With this energy estimate in hand, one can prove existence of solution $h$ in $\mathbf{E}$ by the standard approximation argument. Namely, one can approximate $b$ and $v_r$ by a sequence of smooth vector fields $b^{(n)}$ and functions $v^{(n)}_r$. Then one solves \eqref{ibhth} with $b$ replaced by $b^{(n)}$, and $v_r$ replaced by  $v^{(n)}_r$ respectively to obtain a sequence of functions $\{ h_n \}$ which are uniformly bounded in the norm of $\mathbf{E}$. Then a weak limit of a subsequence of $\{ h_n \}$, say $h$,  is a solution to \eqref{ibhth}. In addition, since $\Vert h_0 \Vert_\infty$ is finite due to $v_0 \in W^{2, 2}(D_m)$, one can prove that $ \Vert h \Vert_\infty $ on $D_m \times [0, T]$ is finite for each $T>0$. The proof of these results are essentially given in \cite{Zcmp04}. The only difference is that one needs to replace the local space time argument by a global one on $D_m \times [0, T]$, which makes the matter simpler. The boundary terms in the integrations all vanish due to the boundary condition of $h$ and $b$.  Now that we have a solution $h$ to \eqref{ibhth}, it is straight forward to verify that
 $v_\th \equiv r h$ is a solution to \eqref{ibvth}. In addition, by the above description of energy estimate,
 \be
 \al
 \lab{vthot}
 \Vert v_\th - (v_0)_\th \Vert_{\mathbf{E}}   &\le \alpha(T, \Vert v_0 \Vert_\infty, \Vert v_0 \Vert_{W^{1, 2}}),\\
 \Vert v_\th - (v_0)_\th \Vert_{\infty} &\le C +  \alpha(T, \Vert v_0 \Vert_\infty, \Vert v_0 \Vert_{W^{1, 2}}).
 \eal
 \ee on $D_m \times [0, T]$. Here $\alpha$ is a function such that $\alpha(T, ., .) \to 0$ as $T \to 0$ and $\alpha(T, ., .) \le C e^{ \alpha_0 T}$ for $T \ge 1$ and  constants $\alpha_0>0$ and $C$ which may depend on $m$.
 Because these inequalities with the exponential bound on $\alpha$ will be used several times, we will give a sketch of the proof in the next paragraph.

 Since $r$ is bounded between two positive constants in $D_m$, we can just prove a similar bound for $h$.
 Using $h$ as a test function in \erf{ibhth} and integrating out the drift term, we see  that
 \[
 \bali
 \int^T_0 \int_{D_m} |\na h |^2 dxdt + &\frac{1}{2} \int_{D_m} | h(x, T) |^2 dx \\
&= \frac{1}{2} \int_{D_m} |h(x, 0)|^2 dx + 2 \int^T_0 \int_{D_m} \frac{\p_r h}{r} h  dxdt
  -\int^T_0 \int_{D_m} \frac{ 2 v_r}{r} h^2  dxdt.
 \eali
 \]As explained, the finiteness of $L^\infty_t L^2_x$ norm of $v_r/r$ implies $v_r/r$ subcritical as a potential function, hence, the following energy bound holds
 \[
 \int^T_0 \int_{D_m} |\na h |^2 dxdt +  \int_{D_m} | h(x, T) |^2 dx \le  \int_{D_m} |h(x, 0)|^2 dx + \overline{C}^2(\Vert v_0 \Vert_{L^2(D_m)}, m) \int^T_0 \int_{D_m}  h ^2 dxdt
 \] Here we have used the energy inequality for $v$, which tells us
  the $L^\infty_t L^2_x$ norm of $v_r/r$ is bounded by $\Vert v_0 \Vert_{L^2(D_m)}$. From here,  Gronwall's inequality infers
 \[
 \Vert h \Vert_{\mathbf{E}} \le C \, \Vert v_0 \Vert_{L^2(D_m)} e^{ C \, \overline{C} \, T}, \qquad \text{on}
 \quad D_m \times [0, T].
\]This and the argument in \cite{Zcmp04}  (mean value inequality with very singular divergence free drift terms) imply
the 2nd inequality in \erf{vthot} and the first one when $T \ge 1$. See also the lines around  \erf{i242p} below. To prove the first one for small $T$, we notice that $\Delta h_0$ can be regarded as an element in $W^{-1, 2}(D_m)$. So the function $h-h_0$ can be regarded as a energy solution to the problem
\be
\lab{ibhth0}
\begin{cases}
\al
 (\Delta+\frac{2}{r} \p_r \big
)(h-h_0)&- b \na (h-h_0)-\frac{ 2 v_r}{r} (h-h_0) -
\partial_t (h-h_0)\\
& = -\Delta h_0 -\frac{2}{r} \p_r h_0 + b \na h_0  + \frac{ 2 v_r}{r} h_0, \quad \text{on} \quad D_m \times (0, T];
\eal\\
\p_n (h-h_0)=0,\, \text{on} \quad \p D_m \times (0, T];\\
h(x, 0)-h_0= 0, \quad x \in D_m.
\end{cases}
\ee Using $h-h_0$ as a test function on \erf{ibhth0}, we deduce, after integration by parts,
\[
 \bali
& \int^T_0 \int_{D_m} |\na (h-h_0) |^2 dxdt + \frac{1}{2} \int_{D_m} | (h-h_0)(x, T) |^2 dx \\
&= \quad 2 \int^T_0 \int_{D_m} \frac{\p_r (h-h_0)}{r} (h-h_0)  dxdt
  -\int^T_0 \int_{D_m} \frac{ 2 v_r}{r} (h-h_0)^2  dxdt -\int^T_0 \int_{D_m} \na h_0 \na (h-h_0) dxdt\\
&\quad
  +\int^T_0 \int_{D_m} [b \na (h-h_0) ] h_0 dxdt + 2 \int^T_0 \int_{D_m} \frac{\p_r h_0}{r} (h-h_0)  dxdt
   -\int^T_0 \int_{D_m} \frac{ 2 v_r}{r} h_0 (h-h_0)  dxdt.
 \eali
\]Applying Cauchy Schwarz to absorb the gradient terms on the right side, we deduce
\[
 \int_{D_m} | (h-h_0)(x, T) |^2 dx \le C(m) T ( \Vert v_0 \Vert^2_{L^2(D_m)} + \Vert (v_0)_\th \Vert^2_{L^\infty(D_m)} + \Vert (v_0)_\th \Vert^2_{W^{1, 2}(D_m)}),
\]which yields the first inequality of \erf{vthot} when $T$ is small.

 \medskip

{\it Step 3.  constructing an intermediate angular vorticity $\tilde{\o_\th}$.}

With the vector field $b$ from Step 1 and $v_\th$ from Step 2,
we will introduce a function
\be
\lab{tiloth}
\tilde{\o_\th} = r \O
\ee where $\O$ is the solution in $\mathbf{E}$ of the initial boundary value problem:
\be
\lab{iboth}
\begin{cases}
 (\Delta+\frac{2}{r} \p_r \big
)\O- b \na \O + \frac{ 2 v_\th}{r^2} \p_{x_3} v_\th -
\partial_t \O=0, \quad \text{on} \quad D_m \times (0, T];\\
\O=0,\,  \text{on} \quad \p D_m \times (0, T];\\
\O(x, 0)= (\o_0)_\th(x)/r, \quad x \in D_m.
\end{cases}
\ee Here $\o_0= curl \, v_0$. This is an initial Dirichlet problem of the heat equation with a nonhomogeneous term and a divergence free drift term. The advantage over the equation of $\o_{\th}$ is the absence of $v_r/r$ term, which will be exploited further in the next section. The existence and uniqueness of solution $\O$ to \eqref{iboth} in $\mathbf{E}$ can be proven in a similar manner to that of \eqref{ibhth}. The key is the a priori energy estimate even with a nonstandard drift term. Since we will need to track dependence of $\O$ on $v_\th$, let us carry out the a priori energy estimate.

If $\O$ is a solution of \eqref{iboth} in $\mathbf{E}$, then it is a legal test function for the equation in \eqref{iboth}. This gives
\be
\lab{Onengl2}
\bali
 \int^T_0 \int_{D_m} |\na \O |^2 dxdt + &\frac{1}{2} \int_{D_m} | \O(x, T) |^2 dx \\
&= \frac{1}{2} \int_{D_m} |\O(x, 0)|^2 dx + 2 \int^T_0 \int_{D_m} \frac{\p_r \O}{r} \O  dxdt
  -\int^T_0 \int_{D_m} \frac{v^2_\th}{r^2} \p_{x_3} \O  dxdt.
 \eali
\ee Here the drift term is integrated out as explained in Step 2. Using the $0$ boundary condition of $\O$ and integration by parts, we see that the 2nd integral on the right hand side is $0$. Then applying Cauchy-Schwarz inequality on the last integral, we deduce
\be
\lab{enOm}
\bali
 \int^T_0 \int_{D_m} |\na \O |^2 dxdt &+ \int_{D_m} | \O(x, T) |^2 dx \\
&\le  \int_{D_m} |\O(x, 0)|^2 dx
  +\int^T_0 \int_{D_m} \frac{v^4_\th}{r^4}  dxdt.
 \eali
\ee From \eqref{vthot}, $v_\th$ is a bounded function and hence the right hand side of \eqref{enOm} is finite, giving us the energy estimate. Note that $r$ is bounded away from $0$ on $D_m$ and $\O(\cdot, 0)=\o_\th(\cdot, 0)/r \in W^{1, 2}(D_m)$ by assumption. As pointed out in Step 2, existence and uniqueness of the solution to \eqref{iboth} in $\mathbf{E}$ follows. Thus $\tilde{\o_\th}$ is well defined, and by direct calculation from \eqref{iboth} and \erf{tiloth},  it satisfies the equation:
\be
\lab{eqotil}
\begin{cases}
 (\Delta-\frac{1}{r^2}  \big
)\tilde{\omega_{\theta}}-(b \nabla)\tilde{\omega_{\theta}}+2\frac{v_{\theta}}
{r}\p_{x_3} v_{\theta}+ \tilde{\omega_{\theta}}\frac{v_r}{r}-\p_t
\tilde{\omega_{\theta}}=0,\\
\tilde{\o_\th}=0,\,  \text{on} \quad \p D_m \times (0, T];\\
\tilde{\o_\th}(x, 0)= (\o_0)_\th(x), \quad x \in D_m.
\end{cases}
\ee We mention that this function $\tilde{\o_\th}$ may not be equal to $curl \, b$ yet.  In the next step, we will use $\tilde{\o_\th}$ to construct a vector field $\tilde{b}$ so that $curl \, \tilde{b}= \tilde{\o_\th}$. Eventually we will prove that the map that assigns $b$ to $\tilde{b}$ has a fixed point.
 For such a fixed point $\tilde{\o_\th} = curl \, b = curl \, \tilde{b}$.

In addition, we can prove higher integrability for $\O$ which will be useful later.
 Since $\O(\cdot, 0)=\o_\th(\cdot, 0)/r \in W^{1, 2}(D_m)$ with $0$ boundary, we can regard it as  $W^{1, 2}$ function in a polygon in the $rx_3$ plane, with $0$ boundary value. By the two dimensional Sobolev inequality, we know $\O(\cdot, 0) \in L^q(D_m)$ for any $q \ge 1$.  Using standard energy estimate for $\O^q$, using the fact that the drift terms are integrated out, we have
 \be
 \lab{olp}
 \Vert \O(\cdot, T) \Vert_{L^q(D_m)} \le \bar{C}(T, q, \Vert v_\th \Vert_{L^\infty}) \Vert \O(\cdot, 0) \Vert_{L^q(D_m)},
 \ee where
 $\bar{C}$ depends on $T$ in an exponential way: $\bar{C}(T, ., .) \le C e^{\alpha_0 T}$ for some constant $\alpha_0>0$. Here $\Vert v_\th \Vert_{L^\infty}$ is on $D_m \times [0, T]$ which has at most exponential growth by \erf{vthot}. Actually, by Moser's iteration, one can also prove that
 $\O(\cdot, t)$ is bounded as soon as $t>0$. See also Lemma \ref{leothinft} below.
 \medskip

{\it Step 4.  Constructing a contraction map from $\overline{B_{\mathbf{S}}(b(0), 1)}$ into itself.}

{\it Step 4.1.}  Definition of the map.

Using the function $\tilde{\o_\th}$ in Step 3, we construct two functions $\tilde{v_r}, \tilde{v_3} \in
\mathbf{E}$ by solving respectively, for each $t \in (0, T]$, the elliptic problems in $W^{1, 2}(D_m)$:
\be
\lab{equrt}
\begin{cases}
\left( \Delta-\frac{1}{r^2} \right) \tilde{v_r} =  \p_{x_3} \tilde{\o_\th}=\p_{x_3} (r \O), \quad \text{in} \quad D_m,\\
\p_{x_3} \tilde{v_r} = 0, \quad \text{on} \quad \p^H D_m; \quad \tilde{v_r} = 0, \quad \text{on} \quad \p^V D_m,
\end{cases}
\ee
\be
\lab{equ3t}
\begin{cases}
\Delta \tilde{v_3} =  - \left(\p_r \tilde{\o_\th} + \frac{\tilde{\o_\th}}{r}\right) = -\frac{1}{r} \p_r (r^2 \O), \quad \text{in} \quad D_m,\\
 \tilde{v_3} = 0, \quad \text{on} \quad \p^H D_m; \quad \p_r \tilde{v_3} = 0, \quad \text{on} \quad \p^V D_m.
\end{cases}
\ee  From Step 3,  $\tilde{\o_\th} \in \mathbf{E}$, hence the right hand side of the above two equations are in $L^2(D_m)$ for a.e. $t \in (0, T]$. For these $t$, problems \eqref{equrt} and \eqref{equ3t} can be regarded as 2 dimensional elliptic problems in a polygon in the $rx_3$ plane, with mixed Dirichlet Neumann boundary values. Note that $\Delta = \p^2_r + \p^2_{x_3} + (1/r) \p_r$. These problems have been well studied. For example, the existence and uniqueness of these solutions can be found in the Grisvard's book \cite{Gris} Chapter 4, together with the following properties.
For $s=3^-/2$, there exists a positive constant $\tilde{C_m}$, depending only on $m$ and $s$,  such that, for a.e. $t$,
\be
\lab{vrv3oth}
\bali
&\Vert \tilde{v_r} \Vert_{W^{1, 2}(D_m)}  + \Vert \tilde{v_3} \Vert_{W^{1, 2}(D_m)} \le  \tilde{C_m}
\Vert \tilde{\o_\th} \Vert_{L^2(D_m)},\\
&\Vert \tilde{v_r} \Vert_{W^{2, s}(D_m)}  + \Vert \tilde{v_3} \Vert_{W^{2, s}(D_m)} \le  \tilde{C_m}
\Vert \tilde{\o_\th} \Vert_{W^{1, 2}(D_m)}.
\eali
\ee In the above the arguments for all functions involved are $(\cdot, t)$. It should be pointed out that the due to the presence of nonconvex  corners of angle $ 3 \pi/2$, solutions $\tilde{v_r}$ and $\tilde{v_3}$ are not in $W^{2, 2}(D_m)$ in general; and the exponent $s=3^-/2$ is determined by this angle. Another remark is that the lower order coefficient $-1/r^2$ in \eqref{equrt} has a good sign that does not affect the solvability of the problem. In addition, by \erf{olp} with say $q=5$ and  Moser's iteration on the whole domain $D_m$ for the equations \erf{equrt} and \erf{equ3t}, we deduce
\be
\lab{vr3linft}
\Vert \tilde{v_r} \Vert_{L^\infty(D_m)}  + \Vert \tilde{v_3} \Vert_{L^\infty(D_m)} \le  \tilde{C_m}
\hat{C} e^{\alpha_0 T}.
\ee Here $\hat{C}$ depends only on the initial value.

Since $\O \in \mathbf{E}$ by \eqref{enOm}, we know $\o_\th \in \mathbf{E}$ and hence, by \eqref{vrv3oth}, the vector field
\be
\tilde{b} \equiv \tilde{v_r} e_r + \tilde{v_3} e_3
\ee also lies in $\mathbf{E} \cap L^\infty_t W^{1, 2}(D_m)$. Now we define a map $\mathbb{L}$ from
$\overline{B_{\mathbf{S}}(b(0), 1)}$ into $\mathbf{E}$ by
\be
\lab{demapl}
\mathbb{L} b = \tilde{b}.
\ee

We must first prove $div \, \tilde{b}=0$ a.e. $t$. Recall that $div \, \tilde{b} = \p_r \tilde{v_r} + \tilde{v_r}/r + \p_{x_3} \tilde{v_3}$. By direct computation from \erf{equrt},
\[
\al
&\left( \Delta-\frac{2}{r^2} \right) \p_r \tilde{v_r}  + \frac{2}{r^3} \tilde{v_r} =  \p_r \p_{x_3} \tilde{\o_\th},\\
& \Delta  (\frac{\tilde{v_r}}{r})  + \frac{2}{r^2}  \p_r \tilde{v_r}- \frac{2}{r^3} \tilde{v_r} =   \p_{x_3} \frac{\tilde{\o_\th}}{r}.
\eal
\]Adding these two equations and  \erf{equ3t}, we find, in pointwise sense,
\be
\lab{eqdiv}
\begin{cases}
\Delta \, div \, \tilde{b} =0, \quad \text{in} \quad D_m,\\
 \p_n div \, \tilde{b} = 0, \quad \text{on} \quad \p D_m, \quad \text{except at corners},
\end{cases}
\ee where $n$ is the unit outward normal of $\p D_m$ except at corners. The Neumann boundary condition is a consequence of the boundary condition of $\tilde{b}$, \erf{equrt} and \erf{equ3t}. Indeed, on $\p^H D_m$ except for corners, we have
\[
\al
\p_{x_3} div \, \tilde{b}& = \p_r \p_{x_3} \tilde{v_r} + \frac{1}{r} \p_{x_3} \tilde{ v_r} + \p^2_{x_3} \tilde{ v_3} = \p^2_{x_3} \tilde{ v_3}\\
&= - \p^2_r \tilde{ v_3} -\frac{1}{r} \p_r \tilde{ v_3} - \left(\p_r \tilde{\o_\th} + \frac{\tilde{\o_\th}}{r}\right), \quad \text{by} \quad \erf{equ3t},\\
&=0, \quad \text{by boundary condition of } \quad \tilde{ v_3}, \, \tilde{\o_\th}=0.
\eal
\]On $\p^V D_m$ except for corners, we have
\[
\al
\p_r div \, \tilde{b}& = \p^2_r  \tilde{v_r} + \frac{1}{r} \p_r \tilde{ v_r} - \frac{1}{r^2}  \tilde{ v_r}+ \p_r \p_{x_3} \tilde{ v_3} \\
& = \p^2_r  \tilde{v_r} + \frac{1}{r} \p_r \tilde{ v_r} - \frac{1}{r^2}  \tilde{ v_r} +\p^2_{x_3} \tilde{ v_r} \quad \text{by} \quad \tilde{ v_r}=0, \, \p_r \tilde{ v_3} =0, \\
&=\p_{x_3} \tilde{\o_\th}, \quad \text{by} \quad \erf{equrt},\\
&=0, \quad \text{by } \quad  \tilde{\o_\th}=0.
\eal
\]

Using $div \, b \in W^{1, s}(D_m)$ with $s=3^-/2$ from \erf{vrv3oth}, we work on \erf{eqdiv} to show that $div \, \tilde{b}$ is identically $0$. Let us mention that
the uniqueness result is not a standard one since $div \, \tilde{b}$ is not known to be in the energy space $W^{1, 2}(D_m)$ at the moment. Actually we are not able to prove the uniqueness of the Neumann problem for $W^{1, s}$ solutions directly.
For $(r, x_3)$ in the interior of $D_m$, consider the function
\[
\al
f&=f(r, x_3) =\int^{x_3}_0  div \, \tilde{b}(r, y_3)  dy_3 = \int^{x_3}_0 [\p_r \tilde{v_r} + \tilde{v_r}/r + \p_{y_3} \tilde{v_3}](r, y_3)  dy_3\\
&=\int^{x_3}_0 [\p_r \tilde{v_r} + \tilde{v_r}/r](r, y_3)  dy_3 + \tilde{v_3}(r, x_3).
\eal
\]which is well defined and smooth since the line segments avoids possible singularity.  Hence
\[
\p_r f = \int^{x_3}_0 [\p^2_r \tilde{v_r} + \p_r \tilde{v_r}/r -  \tilde{v_r}/r^2](r, y_3)  dy_3 + \p_r \tilde{v_3}(r, x_3),
\]which infers, from equation \erf{equrt}, that
\be
\lab{prfdv}
\p_r f = \int^{x_3}_0(-\p^2_{y_3} \tilde{v_r} + \p_{y_3} \tilde{\o_\th})(r, y_3)  dy_3 + \p_r \tilde{v_3}(r, x_3)= \left(-\p_{x_3} \tilde{v_r} + \tilde{\o_\th} + \p_r \tilde{v_3}\right)(r, x_3).
\ee Here we have used the boundary conditions on $\p^H D_m$.
On the other hand, on interior points of $D_m$, we can differentiate equations  \erf{equrt} and \erf{equ3t}
suitably to reach
\[
\begin{cases}
\left( \Delta-\frac{1}{r^2} \right) (\p_{x_3} \tilde{v_r} -\p_r \tilde{v_3} - \tilde{\o_\th})  =  0, \quad \text{in} \quad D^0_m,\\
\quad \p_{x_3} \tilde{v_r} -\p_r \tilde{v_3} - \tilde{\o_\th} =0  \quad \text{on} \quad \p D_m,
\end{cases}
\]except at corners. Notice that $\p_{x_3} \tilde{v_r} -\p_r \tilde{v_3} - \tilde{\o_\th} \in W^{1, s}(D_m)$ with $s=3^-/2$. By the uniqueness result in Proposition \ref{prunid}, we find
\[
\tilde{\o_\th}=\p_{x_3} \tilde{v_r} -\p_r \tilde{v_3}.
\] Substituting this to \erf{prfdv}, we infer in the interior points of $D_m$,
\be
\lab{prf0}
\p_r f(r, x_3)=0.
\ee For these points, we also know from \erf{eqdiv} that
\[
\al
\Delta f&= \p^2_r f + \frac{1}{r} \p_r f + \p^2_{x_3} f\\
&=\int^{x_3}_0 (\p^2_r + \frac{1}{r} \p_r)  div \, \tilde{b} \,  dy_3 + \p_{x_3}  div \, \tilde{b}\\
&=-\int^{x_3}_0  \p^2_{y_3}  div \, \tilde{b} \, dy_3 + \p_{x_3}  div \, \tilde{b}\\
&=0, \quad \text{due to} \quad \p_{x_3} div \, \tilde{b}(r, 0)=0.
\eal
\]Since $f=f(x_3)$  by \erf{prf0}, this shows $\p^2_{x_3} f=0$ so that $f= c x_3$.

Next we observe that for $W^{2, s}$ solutions to \erf{equrt}, we have
\[
\int_{D_m \cap \{r = const. \}} \tilde{v_r}(r, x_3, t) dx_3 =0.
\]This result, which only depends on the $W^{2, s}$ property of $\tilde{v_r}(\cdot, t)$ , is stated as Lemma \ref{leintx3w1s} below. To avoid interrupting the flow of the presentation, we will prove the lemma till the end of the section.   This shows that $f=0$ and hence $div \, \tilde{b}=\p_{x_3} f =0$.

In the next few  sub-steps, we will prove that $\mathbb{L}$ maps $\overline{B_{\mathbf{S}}(b(0), 1)}$ into itself and it is a contraction, provided $T$ is sufficiently small.

{\it Step 4.2.} We prove $\mathbb{L}$ is a contraction.

Let $b^{(i)} \equiv v^{(i)}_r e_r + v^{(i)}_3 e_3 \in \overline{B_{\mathbf{E}}(b(0), 1)}$, $i=1, 2$, be two given vector fields with the same initial value $b_0$.  Then they determine, by Step 2, $v^{(i)}_\th$ and henceforth $\O^{(i)}$ by Step 3. Therefore
$\O^{(2)}- \O^{(1)}$ satisfies
\be
\lab{ibo12}
\bali
\begin{cases}
 &(\Delta+\frac{2}{r} \p_r \big
)(\O^{(2)}- \O^{(1)})- b^{(2)} \na (\O^{(2)}- \O^{(1)}) - ( b^{(2)} - b^{(1)} ) \na \O^{(1)} \\
&\qquad + \frac{1}{r^2} \p_{x_3}  \left( (v^{(2)}_\th)^2 - (v^{(1)}_\th)^2\right)-
\partial_t (\O^{(2)}- \O^{(1)})=0;\\
&(\O^{(2)}- \O^{(1)})=0,\,  \text{on} \quad \p D_m \times (0, T];\\
&(\O^{(2)}- \O^{(1)})(x, 0)= 0, \quad x \in D_m.
\end{cases}
\eali
\ee
Then along the same way as deriving \eqref{enOm}, we arrive at
\be
\bali
 \int^T_0 \int_{D_m} &|\na (\O^{(2)}-\O^{(1)}) |^2 dxdt + \sup_{t \in [0, T]}\int_{D_m} | \O^{(2)}-\O^{(1)} |^2(x, t) dx\\
&\le
  \int^T_0 \int_{D_m} |(v^{(2)}_\th)^2-(v^{(1)}_\th)^2|^2 \, r^{-4}  dxdt
  - 2 \int^T_0 \int_{D_m} ( b^{(2)} - b^{(1)} ) \na \O^{(1)} \,  (\O^{(2)}-\O^{(1)}) dxdt\\
&=\int^T_0 \int_{D_m} |(v^{(2)}_\th)^2-(v^{(1)}_\th)^2|^2 \, r^{-4}  dxdt
  + 2 \int^T_0 \int_{D_m} ( b^{(2)} - b^{(1)} )  \O^{(1)} \,  \na (\O^{(2)}-\O^{(1)}) dxdt.
\eali
\ee Here we have used the divergence free property of $b^{(i)}$ and integration by parts. Using Cauchy Schwarz inequality, we find
\be
\bali
\int^T_0 \int_{D_m} &|\na (\O^{(2)}-\O^{(1)}) |^2 dxdt + 2 \sup_{t \in [0, T]}\int_{D_m} | \O^{(2)}-\O^{(1)} |^2(x, t) dx\\
&\le 2\int^T_0 \int_{D_m} |(v^{(2)}_\th)^2-(v^{(1)}_\th)^2|^2 \, r^{-4}  dxdt
  + 4 \int^T_0 \int_{D_m} ( b^{(2)} - b^{(1)} )^2  (\O^{(1)})^2 dxdt\\
&\le 2 \int^T_0 \int_{D_m} |(v^{(2)}_\th)^2-(v^{(1)}_\th)^2|^2 \, r^{-4}  dxdt
  + 4 \Vert b^{(2)} - b^{(1)} \Vert^2_{L^{10/3}} \Vert \O^{(1)} \Vert^2_{L^5} \quad \text{(by H\"older inequality)}. \\
\eali
\ee From this, Sobolev inequality and interpolation, we deduce, after using \eqref{olp} with $q=5$,
\be
\lab{o2-o1}
\Vert \O^{(2)} - \O^{(1)} \Vert^2_{\mathbf{E}} \le C \int^T_0 \int_{D_m} |(v^{(2)}_\th)^2-(v^{(1)}_\th)^2|^2 \, r^{-4}  dxdt + C  \bar{C}^2(T, \Vert v_\th \Vert_{L^\infty}) \Vert \O(\cdot, 0) \Vert^2_{L^5(D_m)} T^{2/5}  \Vert b^{(2)} - b^{(1)} \Vert^2_{\mathbf{E}}
\ee

 Write $\tilde{\o_\th}^{(i)} = r \O^{(i)}$,  and $\tilde{b}^{(i)}=\tilde{v_r}^{(i)} e_r +
\tilde{v_3}^{(i)} e_3$, $i=1, 2$, where $\tilde{v_r}^{(i)}$ and $\tilde{v_3}^{(i)}$ are determined respectively from
\eqref{equrt} and \eqref{equ3t} with $\tilde{\o_\th}$ replaced by $\tilde{\o_\th}^{(i)}$.
Then, just like \eqref{vrv3oth}, we have
\be
\lab{lb2b10}
\bali
&\Vert \mathbb{L} b^{(2)} - \mathbb{L} b^{(1)} \Vert_{\mathbf{E}} = \Vert \tilde{b}^{(2)} - \tilde{b}^{(1)} \Vert_{\mathbf{E}} \le \tilde{C_m} (1+T)  \sup_{t \in [0, T]} \Vert \tilde{\o_\th}^{(2)} - \tilde{\o_\th}^{(1)}
\Vert_{L^2(D_m)}\\
&\le \tilde{C_m} (1+T)  \sup_{t \in [0, T]} \Vert \O^{(2)} - \O^{(1)}
\Vert_{L^2(D_m)}\\
&\le C  \bar{C}(T, \Vert v_\th \Vert_{L^\infty}) \Vert \O(\cdot, 0) \Vert_{L^5(D_m)} T^{1/5}  \Vert b^{(2)} - b^{(1)} \Vert_{\mathbf{E}} \\
&\qquad  \qquad + \tilde{C_m} (1+T)  \left(\int^T_0 \int_{D_m} |(v^{(2)}_\th)^2-(v^{(2)}_\th)^2|^2 \,  r^{-4}  dxdt \right)^{1/2}.
\eali
\ee Therefore,
\be
\lab{lb2b11}
\bali
\Vert \mathbb{L} b^{(2)} - \mathbb{L} b^{(1)} \Vert_{\mathbf{E}} &\le C  \bar{C}(T, \Vert v_\th \Vert_{L^\infty}) \Vert \O(\cdot, 0) \Vert_{L^5(D_m)} T^{1/5}  \Vert b^{(2)} - b^{(1)} \Vert_{\mathbf{E}} \\
&\qquad+ C \tilde{C_m} (1+T)  \left[\Vert v_0 \Vert_{\mathbf{E}} + \alpha(T,  \Vert v_0 \Vert_{\mathbf{E}}, \Vert v_0 \Vert_{W^{2, 2}})\right]  \left(\int^T_0 \int_{D_m} |v^{(2)}_\th-v^{(1)}_\th|^2 \,  r^{-4}  dxdt \right)^{1/2}.
\eali
\ee Here we just used \eqref{o2-o1} and \eqref{vthot} in the last lines. To bound the last integral, we
write $h_i = v^{(i)}_\th/r$, $i=1, 2$. From \eqref{ibhth}, we see that the function $h_2-h_1$ satisfies the initial boundary value problem
\be
\lab{ibh21}
\bali
\begin{cases}
& (\Delta+\frac{2}{r} \p_r \big
)(h_2-h_1)- b^{(2)} \na (h_2-h_1) - (b^{(2)}-b^{(1)}) \na h_1\\
&\qquad -\frac{ 2 v^{(2)}_r}{r} (h_2-h_1)
- \frac{2 ( v^{(2)}_r -v^{(1)}_r) }{r} h_1-
\partial_t (h_2-h_1)=0, \quad \text{on} \quad D_m \times (0, T];\\
&\p_n (h_2-h_1)=0,\, \text{on} \quad \p D_m \times (0, T];\\
&(h_2-h_1)(x, 0)=0, \quad x \in D_m.
\end{cases}
\eali
\ee Using $h_2-h_1$ as a test function in the preceding equation, we see that
\be
\lab{enh2h1}
\bali
\int^T_0 &\int_{D_m} |\na (h_2-h_1) |^2 dxdt + \frac{1}{2} \int_{D_m} (h_2-h_1)^2(x, T) dx \\
&=\int^T_0 \int_{D_m} \frac{2}{r} \p_r (h_2-h_1)  \, (h_2-h_1) dxdt
-\int^T_0 \int_{D_m} (b^{(2)}-b^{(1)}) \na h_1 \, (h_2-h_1) dxdt\\
&\qquad - \int^T_0 \int_{D_m} \frac{ 2 v^{(2)}_r}{r} (h_2-h_1)^2 dxdt
- \int^T_0 \int_{D_m} \frac{2 ( v^{(2)}_r -v^{(1)}_r) }{r} h_1 (h_2-h_1) dxdt\\
&\equiv I_1 + I_2+ I_3 + I_4.
\eali
\ee Next we find suitable bounds for $I_1, ..., I_4$. For simplicity we denote $D_m \times [0, T]$ by $Q_T$ here.

By Cauchy Schwarz inequality
\be
\lab{hi1}
|I_1| \le (1/4) \Vert \na (h_2-h_1) \Vert^2_{L^2(Q_T)} + C(m) \Vert  h_2-h_1 \Vert^2_{L^2(Q_T)};
\ee Using integration by parts and divergence free property,
\be
\lab{hi2}
\bali
|I_2| &= \left| \int^T_0 \int_{D_m} (b^{(2)}-b^{(1)})  h_1 \, \na (h_2-h_1) dxdt \right| \\
&\le (1/4) \Vert \na (h_2-h_1) \Vert^2_{L^2(Q_T)} + C \Vert h_1 \Vert^2_{L^\infty(Q_T)} T  \Vert  b^{(2)}- b^{(1)} \Vert^2_{\mathbf{E}};
\eali
\ee By H\"older inequality, Sobolev inequality and interpolation,
\be
\lab{hi3}
\bali
|I_3| & \le C(m) \Vert v^{(2)}_r \Vert_{L^{10/3}(Q_T)}  \left(\int^T_0 \int_{D_m} |h_2-h_1|^{10/3} dxdt\right)^{3/5}  (T |D_m|)^{1/10}\\
&\le C C(m) T^{1/10} \,  \Vert v^{(2)}_r \Vert_{\mathbf{E}} \, \Vert h_2-h_1 \Vert^2_{\mathbf{E}} ;
\eali
\ee
\be
\lab{hi4}
|I_4|  \le C(m)  \, \Vert h_2-h_1 \Vert^2_{L^2(Q_T)} + C T \Vert  h_1 \Vert^2_{L^\infty(Q_T)}  \Vert  b^{(2)}- b^{(1)} \Vert^2_{\mathbf{E}}.
\ee Substituting \eqref{hi1}, ..., \eqref{hi4} into \eqref{enh2h1}, we find, after observing that $I_3$ can be absorbed by the right hand side of \eqref{enh2h1} when $T$ is sufficiently small, that
\[
\bali
\int_{D_m} &(h_2-h_1)^2(x, T) dx \\
&\le C C(m) \int^T_0 \int_{D_m} |h_2-h_1|^2 dx dt +
 C \Vert h_1 \Vert^2_{L^\infty(Q_T)} T  \Vert  b^{(2)}- b^{(1)} \Vert^2_{\mathbf{E}}.
\eali
\] By the usual trick, this inequality actually holds when $T$ is replaced by any time $T' \in [0, T]$. So Gronwall's inequality then infers, for these $T'$,
\be
\lab{vth12-l2}
\int_{D_m} |v^{(2)}_\th-v^{(1)}_\th|^2(x, T') r^{-2} dx=\int_{D_m} (h_2-h_1)^2(x, T') dx \le \bar{C}(T, m, \Vert v_0 \Vert_{W^{2, 2}(D_m)}) T \Vert  b^{(2)}- b^{(1)} \Vert^2_{\mathbf{E}}.
\ee Here the constant $\bar{C}$ is bounded when $T, m$ are bounded.  Substituting this to the last term on the right hand side of \eqref{lb2b11}, we find that
\[
\Vert \mathbb{L} b^{(2)} - \mathbb{L} b^{(1)} \Vert_{\mathbf{E}} \le 0.5 \Vert  b^{(2)}- b^{(1)} \Vert_{\mathbf{E}}
\]if $T$ is sufficiently small.
Using \erf{o2-o1} and \erf{vth12-l2}, we also deduce
\[
\Vert \tilde{\o}^{(2)}_\th - \tilde{\o}^{(1)}_\th \Vert^2_{\mathbf{E}} \le  C \bar{C}_m \bar{C}^2(T,  \Vert v_\th \Vert_{L^\infty}) \Vert \O(\cdot, 0) \Vert^2_{L^5(D_m)} T^{2/5}  \Vert b^{(2)} - b^{(1)} \Vert^2_{\mathbf{E}}.
\]Combining the last two inequalities, we conclude
\[
\Vert \mathbb{L} b^{(2)} - \mathbb{L} b^{(1)} \Vert_{\mathbf{S}} \le 0.75 \Vert  b^{(2)}- b^{(1)} \Vert_{\mathbf{E}} \le 0.75 \Vert  b^{(2)}- b^{(1)} \Vert_{\mathbf{S}}.
\]
Therefore the map $\mathbb{L}$ is a contraction under the norm of $\mathbf{S}$ if $T$ is sufficiently small.
\medskip

{\it Step 4.3.} We prove
$\mathbb{L}$ maps $\overline{B_{\mathbf{S}}(b(0), 1)}$ into itself if $T$ is sufficiently small.

 Pick $b \in \overline{B_{\mathbf{E}}(b(0), 1)}$, then  $\mathbb{L} b = \tilde{b} = \tilde{v_r} e_r + \tilde{v_3} e_3$, where $\tilde{v_r}$ and $\tilde{v_3}$ are given by \erf{equrt} and \erf{equ3t} respectively.

Let $\O_0 \equiv (\o_0)_\th/r$ where $\o_0 = curl \, v_0$. By \eqref{iboth}, $\O-\O_0$ satisfies, in the energy space,

\be
\lab{iboo0}
\begin{cases}
 \left(\Delta+\frac{2}{r} \p_r \right)(\O-\O_0)- b \na (\O-\O_0) + \frac{ 2 v_\th}{r^2} \p_{x_3} v_\th -
\partial_t (\O-\O_0) \\
\qquad = -\left(\Delta+\frac{2}{r} \p_r \right)\O_0 + b \na \O_0, \quad \text{on} \quad D_m \times (0, T];\\
(\O-\O_0)=0,\,  \text{on} \quad \p D_m \times (0, T];\\
(\O-\O_0)(x, 0)=0, \quad x \in D_m.
\end{cases}
\ee Here $\Delta \O_0$ is regarded as an element in $H^{-1}_0(D_m)$. Using $\O-\O_0$ as a test function in the above equation,  we obtain, in a similar manner as \eqref{enOm},
\[
\al
 \int^T_0 \int_{D_m}& |\na (\O -\O_0) |^2 dxdt + \int_{D_m} | (\O-\O_0)(x, T) |^2 dx \\
&\le
  \int^T_0 \int_{D_m} \frac{v^4_\th}{r^4}  dxdt +
   2  \int^T_0 \int_{D_m} \left[ \left(\Delta+\frac{2}{r} \p_r \right)\O_0 - b \na \O_0 \right]
   ( \O -\O_0) dxdt.
\eal
\]This implies, after integration by parts, that
\[
\al
 \int^T_0 \int_{D_m}& |\na (\O -\O_0) |^2 dxdt + \int_{D_m} | (\O-\O_0)(x, T) |^2 dx \\
 &\le
  \int^T_0 \int_{D_m} \frac{v^4_\th}{r^4}  dxdt
   -
   2  \int^T_0 \int_{D_m}  \na \O_0 \na (\O-\O_0) dxdt \\
 & \qquad
   -4 \int^T_0 \int_{D_m}  \frac{\O_0}{r} \p_r (\O-\O_0) dxdt  + 2  \int^T_0 \int_{D_m} b  \O_0 \na
   ( \O -\O_0) dxdt,
\eal
\]which infers, via Cauchy Schwarz inequality, that
\[
\al
 \int^T_0 \int_{D_m}& |\na (\O -\O_0) |^2 dxdt + \int_{D_m} | (\O-\O_0)(x, T) |^2 dx \\
&\le
 64 \int^T_0 \int_{D_m} \frac{v^4_\th}{r^4}  dxdt
  + 64 \int^T_0 \int_{D_m} \left(|\na \O_0|^2 +
 \O^2_0 r^{-2} \right) dxdt + 64 \int^T_0 \int_{D_m}
 |b|^2 \O^2_0  dxdt.
\eal
\]By H\"older, Sobolev inequality and interpolation
\[
\al
 \int^T_0 \int_{D_m}
 |b|^2 \O^2_0  dxdt &\le \Vert b \Vert^2_{L^{10/3}(Q_T)}  \left(\int^T_0 \int_{D_m} |\O_0|^5 dxdt \right)^{2/5}\\
 &\le C_m T^{2/5} \Vert b \Vert^2_{\mathbf{E}} \, \Vert v_0 \Vert^2_{W^{2, 2}}.
\eal
\]Combining the previous two inequalities and using \erf{vthot} and the initial condition, we deduce
\be
\lab{o0ol2}
\int^T_0 \int_{D_m} |\na (\O -\O_0) |^2 dxdt + \sup_{t \in [0, T]} \int_{D_m} | (\O-\O_0)(x, t) |^2 dx \le C_m T^{2/5} \bar{C}(T, \Vert v_0 \Vert^2_{W^{2, 2}}).
\ee

According to the definition in \erf{demapl},  $\mathbb{L} b = \tilde{b} = \tilde{v_r} e_r + \tilde{v_3} e_3$, where $\tilde{v_r}$ and $\tilde{v_3}$ are given by \erf{equrt} and \erf{equ3t} respectively.
Consequently $\tilde{v_r}-(v_0)_r$ and $\tilde{v_3}-(v_0)_3$ satisfy
\be
\lab{equrt0}
\begin{cases}
\left( \Delta-\frac{1}{r^2} \right) (\tilde{v_r}-(v_0)_r) =  \p_{x_3} (\tilde{\o_\th}- (\o_0)_\th) =\p_{x_3} (r (\O-\O_0)), \quad \text{on} \quad D_m,\\
\p_{x_3} (\tilde{v_r}-(v_0)_r) = 0, \quad \text{on} \quad \p^H D_m; \quad (\tilde{v_r}-(v_0)_r) = 0, \quad \text{on} \quad \p^V D_m,
\end{cases}
\ee and
\be
\lab{equ3t0}
\begin{cases}
\Delta (\tilde{v_3}-(v_0)_3) =  - \left(\p_r [\tilde{\o_\th}-(\o_0)_\th] + \frac{\tilde{\o_\th}-(\o_0)_\th}{r}\right) = -\frac{1}{r} \p_r (r^2 [\O-\O_0]), \quad \text{on} \quad D_m,\\
 \tilde{v_3}-(v_0)_3 = 0, \quad \text{on} \quad \p^H D_m; \quad \p_r(\tilde{v_3}-(v_0)_3) = 0, \quad \text{on} \quad \p^V D_m.
\end{cases}
\ee
Using
\erf{vrv3oth}, we obtain, for all $t \in [0, T]$,
\[
\Vert \tilde{v_r}-(v_0)_r  \Vert_{W^{1, 2}(D_m)}  + \Vert \tilde{v_3} -(v_0)_3  \Vert_{W^{1, 2}(D_m)} \le  \tilde{C_m}
\Vert \tilde{\o_\th} -(\o_0)_\th \Vert_{L^2(D_m)}.
\]This and \erf{o0ol2} yield
\[
\Vert \tilde{b} - b(0) \Vert_{\mathbf{S}} \le C_m T^{2/5} \bar{C}(T, \Vert v_0 \Vert^2_{W^{2, 2}}).
\]This proves
$\mathbb{L}$ maps $\overline{B_{\mathbf{S}}(b(0), 1)}$ into itself if $T$ is sufficiently small.

Now the contraction mapping theorem tells us that the map $\mathbb{L}$ has a unique fixed point
$b =v_r e_r + v_3 e_3 \in \overline{B_{\mathbf{S}}(b(0), 1)}$. Namely $b=\mathbb{L} b$, i.e.
\be
\lab{vr3=vr3t}
v_r=\tilde{v_r}, \qquad v_3=\tilde{v_3}.
\ee

 Let $v_\th$ be given by \erf{ibvth}. Now we can claim that
\[
v \equiv v_r e_r + v_\th e_\th + v_3 e_3
\] is a solution to Problem \ref{wtdm} such that $| \na \o_\th | \in L^2(D_m \times [0, T])$, where
$\o_\th = (curl \, v)_\th$.
The reason is the following. Let $\tilde{\o_\th}$ be given by \erf{eqotil}. Then, due to \erf{vr3=vr3t}, $v_r$ and $v_3$ are given by \erf{equrt} and \erf{equ3t} respectively. In vector form, these two equations can be written as
\[
\Delta b = - curl \, (\tilde{\o_\th} e_\th).
\]
 On the other hand, by the above definition of $\o_\th$ and cylindrical curl formula, we actually have $\o_\th e_\th = curl \, b$. Therefore
\[
\Delta b = - curl \, (\o_\th e_\th).
\]Subtraction of the last two equations gives us $curl \,[ (\tilde{\o_\th}- \o_\th) e_\th ]=0.$ Therefore
\[
\Delta \,[ (\tilde{\o_\th}- \o_\th) e_\th] =0.
\]By our boundary condition $\o_\th= \p_3 v_r - \p_r v_3 =0$ on $\p D_m$ and by construction $\tilde{\o_\th}=0$ on $\p D_m$. Since $b \in \mathbf{S}$,  we know that $\tilde{\o_\th}(\cdot, t)$ and $\o_\th(\cdot, t)$ are in $W^{1, 2}(D_m)$ for a.e. $t$.
Hence, by standard uniqueness result of solutions of Laplace equation in $W^{1, 2}$ space, we have $\tilde{\o_\th}(\cdot, t)=\o_\th(\cdot, t)$, a. e. $t$. Now interior regularity of solutions tells us:
\be
\lab{o=otil}
\o_\th=\tilde{\o_\th} \quad \text{in}  \quad  D_m \times [0, T].
\ee By \erf{enOm} and $\tilde{\o_\th} = r \O$, we know that $| \na \o_\th | \in L^2(D_m \times [0, T])$.
Finally, by \erf{eqotil} and \erf{o=otil}, we see that $\o_\th$ satisfies
\be
\lab{eqotil2}
\begin{cases}
 (\Delta-\frac{1}{r^2}  \big
)\omega_{\theta}-(b \nabla)\omega_{\theta}+2\frac{v_{\theta}}
{r}\p_{x_3} v_{\theta}+ \omega_{\theta} \frac{v_r}{r}-\p_t
\omega_{\theta}=0,\\
\o_\th=0,\,  \text{on} \quad \p D_m \times (0, T];\\
\o_\th(x, 0)= (\o_0)_\th(x), \quad x \in D_m.
\end{cases}
\ee Converting the above equation into the vector equation for $\o_\th e_\th= curl \, b$, we check, by vector calculus identity in the cylindrical system and \erf{eqotil2}
\be
\lab{curl[]}
curl \, [\Delta b - b \na b + \frac{v^2_\th}{r} e_r - \p_t b] = 0, \quad \text{pointwise in} \quad D_m \times (0, T],
\ee where $b \na b = (v_r \p_r v_r + v_3 \p_{x_3} v_r) e_r + (v_r \p_r v_3 + v_3 \p_{x_3} v_3) e_3$.
Since $D_m$ can be regarded as a simply connected domain in $rx_3$ plane, the term inside the brackets of
\erf{curl[]} is a gradient field. Hence we know that
$v_r$ and $v_3$ satisfy their respective equations in \erf{eqasns}. Note also that $v_\th$ satisfies its corresponding equation in \erf{eqasns} by construction, c.f. \erf{ibvth}. This proves the claim, namely
 $v$ is a solution to Problem \ref{wtdm} in $\mathbf{S}$ on the time interval $(0, T]$. Due to the bounds
 \erf{vthot}, \erf{olp} and \erf{vr3linft}, this fixed point argument can be continued indefinitely.
Indeed, choose any moment $t_0>0$ such that a solution already exists. We can repeat the above fixed point argument starting from $t_0$. From the a priori bounds\erf{vthot}, \erf{olp} with $q=5$, and \erf{lb2b10} and \erf{o2-o1} with initial time $0$ replaced by $t_0$, we see that the fixed point argument works at least in the interval $[t_0, t_0 + \epsilon e^{-2 \alpha_0 t_0}]$. Here $\alpha_0 >0$ and $\e>0$ are constants which are independent of $t_0$.
 Thus
 we have proven the existence of solutions to Problem \ref{wtdm} in $\mathbf{S}$ for all time. Due to the
 definition of $\mathbf{S}$, we have also proven $|\na \o_\th| \in L^2_tL^2_x$. It remains to prove
 $|\na \p_{x_3} v_r | + | \na \o_r| \in L^2_tL^2_x$.
\medskip

{\it Step 5.  $|\na \p_{x_3} v_r | \in L^2_tL^2_x$.}

According to \erf{vrv3oth}, for a.e. $t$, we have $v_r(\cdot, t) \in W^{2, s}(D_m)$, $\p_{x_3} \o_\th(\cdot, t) \in  L^2(D_m) $ with $s=3^-/2$. From here, in this step, we will always work on these time $t$, and for simplicity,  will suppress the time variable $t$ unless there is a confusion.
 By \erf{equrt}, $\p_{x_3} v_r$ can be regarded as a $W^{1, s}$ solution to the problem:
\be
\lab{eqp3ur}
\begin{cases}
\left( \Delta-\frac{1}{r^2} \right) \p_{x_3} v_r =  \p^2_{x_3} \o_\th, \quad \text{on} \quad D_m,\\
\p_{x_3} v_r = 0, \quad \text{on} \quad \p^H D_m; \quad \p_{x_3} v_r = 0, \quad \text{on} \quad \p^V D_m.
\end{cases}
\ee Notice that the boundary value for $\p_{x_3} v_r$ is actually $0$.
On the other hand, since  $\p_{x_3} \o_\th \in L^2(D_m)$, we can find a $W^{1, 2}$ solution to the following problem.
\be
\lab{eq3f}
\begin{cases}
\left( \Delta-\frac{1}{r^2} \right) f =  \p^2_{x_3} \o_\th, \quad \text{on} \quad D_m,\\
f= 0, \quad \text{on} \quad \p D_m.
\end{cases}
\ee The existence and uniqueness of $W^{1, 2}$ solutions to \erf{eq3f} is standard if the right hand side of the equation is in $L^2(D_m)$. See Chapter 4 of \cite{Gris} e.g.  In our case, although the right hand side is only in $H^{-1}$, the result still holds and the proof is more or less the same. But for completeness, we will give a quick proof by an approximating procedure. First we construct a sequence of functions $g_j \in W^{1, 2}_0(D_m)$ such that $g_j \to \p_{x_3} \o_\th$ in $L^2(D_m)$. Now that $\p_{x_3} g_j$ is in $L^2(D_m)$, the standard theory mentioned above states that the problem
\be
\lab{eq3f_j}
\begin{cases}
\left( \Delta-\frac{1}{r^2} \right) f_j =  \p_{x_3} g_j, \quad \text{on} \quad D_m,\\
f_j= 0, \quad \text{on} \quad \p D_m.
\end{cases}
\ee has an unique solution in $W^{1, 2}_0(D_m)$. Using $f_j$ as a test function in the above equation, we obtain
\[
\int_{D_m} | \na f_j|^2 dx + \int_{D_m} r^{-2} f_j^2 dx =  \int_{D_m} g_j \p_{x_3} f_j dx,
\]which yields the uniform energy bound
\[
\int_{D_m} | \na f_j|^2 dx + \int_{D_m} r^{-2} f_j^2 dx \le 2 \int_{D_m} g^2_j \le 2 \int_{D_m} |\p_{x_3}
\o_\th|^2 dx + C.
\]Therefore, a subsequence of $\{f_j\}$ converges weakly in $W^{1, 2}$ to a $W^{1, 2}_0$ function $f$.
Using a $C^1_0(D_m)$ test function, it is easy to see that this $f$ is a solution to \erf{eq3f}.

So $\p_{x_3} v -f \in W^{1, s}_0(D_m)$ is a solution to the homogeneous problem
\be
\lab{eq3v-f}
\begin{cases}
\left( \Delta-\frac{1}{r^2} \right)(\p_{x_3} v - f) =  0, \quad \text{on} \quad D_m,\\
\p_{x_3} v - f= 0, \quad \text{on} \quad \p D_m.
\end{cases}
\ee According to Proposition \ref{prunid} below, uniqueness holds for the above problem. Thus we have proven $\p_{x_3} v_r =f \in W^{1, 2}(D_m)$, as stated. We comment that, since the solution space $W^{1, s}_0(D_m)$ with $s=3^-/2$ is not the standard energy space and $\p D_m$ has  corners with $3 \pi/2$ angle, more efforts are needed for the  proof of the uniqueness. To avoid interrupting the flow of the main argument, we have moved this to Section 4.

\medskip

{\it Step 6.  $|\na \p_{x_3} v_\th | \in L^2_tL^2_x$ or $|\na J|=|-\na \p_{x_3} (v_\th/r) | \in L^2_tL^2_x$.}
\medskip

The proof is based on a weak-strong uniqueness argument.
It is convenient to consider the function $J = \frac{\o_r}{r}= - \frac{\p_{x_3} v_\th}{r}$, which satisfies the following equation pointwise:
\be
\lab{eqJbc}
\begin{cases}
\Delta J  -(b\cdot\nabla) J +\frac{2}{r}\p_r J +
 (\o_r \p_r + \o_3 \p_{x_3}) \frac{v_r}{r} - \p_t J
=0, \quad \text{on} \quad  D_m \times [0, T]\\
J=0, \quad \text{on} \quad \p^H D_m \times [0, T]; \qquad  \p_n J=0, \quad \text{on} \quad \p^V D_m \times [0, T]\\
J(x, 0) = - \p_{x_3} v_\th(x, 0)/r.
\end{cases}
\ee This equation can be derived from the vorticity equation for $\o_r$ and will also play an important role in the next section. Note the boundary condition for $J$ is a homogeneous mixed Dirichlet-Neumann condition.  The rest of the step is divided into 3 sub-steps.

{\it Step 6.1.}  Here we prove the following assertion:  the Laplace transform of $J(x, \cdot)$ is well defined if the real part of the phase variable $\lam$ is sufficiently large, and it is in $ W^{1, s}(D_m)$ with $s=3^-/2$ if $\lam$ is real and sufficiently large.

There may be different ways to do it. But we will use the Laplace transform to convert the problem to an elliptic one which was already studied before. In order to do this, we will extend the solution $v$ to all time interval $(0, \infty)$.
Due to the bounds \erf{vthot} and \erf{vr3linft}, this is always possible.
We need the following component of the energy inequality
\[
\al
 & 2 \int^\infty_0 \int_{D_m} |\p_{x_3} v_\th|^2 dxdt + 2 \int^T_0 \int_{D_m} |\p_r v_\th - \frac{1}{r} v_\th |^2 dxdt + 2 \int^\infty_0 \int_{D_m} |\p_{x_3} v_r + \p_r v_3 |^2 dxdt \\
 &\qquad + 4 \int^\infty_0 \int_{D_m} \left( |\p_r v_r |^2+\frac{v^2_r}{r^2} \right) dxdt
 + 4 \int^\infty_0 \int_{D_m} |\p_{x_3} v_3 |^2 dxdt
 \le \int_{D_m} |v(x, 0)|^2 dx.
 \eal
\]
which is a consequences of \eqref{en1s} and the formula for the strain tensor
\eqref{Sv}. It is helpful to convert the above inequality into one involving $h=v_\th/r$. Since $ r \p_r h = \p_r v_\th - v_\th/r$, the preceding inequality implies
\[
 2 \int^\infty_0 \int_{D_m} r^2 |\p_{x_3} h|^2 dxdt + 2 \int^\infty_0 \int_{D_m} r^2 |\p_r h |^2 dxdt
 \le \int_{D_m} |v(x, 0)|^2 dx,
\]and consequently
\be
\lab{dhl22}
  \int^T_0 \int_{D_m}  |\na h|^2 dxdt
 \le \bar{C_m}\int_{D_m} |v(x, 0)|^2 dx,
\ee for some $\bar{C_m}$ depending on $m$.

Next we write the equation \erf{ibhth}  for $h$ as
\[
\Delta h -
\partial_t h= - \frac{2}{r} \p_r h +  b \na h+\frac{ 2 v_r}{r} h \equiv R,
\]By \erf{dhl22}, the standard energy bound for $L^\infty_tL^2_x$ norm of $v_\th$ and the $L^\infty$ bound for $v_\th$ \erf{vthot} and $b$ \erf{vr3linft}, we know that $R \in L^2_tL^2_x$ and there is a constant $\beta>0$ such that
\be
\lab{Rh12}
\int^\infty_0 e^{- \beta t} \Vert R(\cdot, t) \Vert^2_{L^2(D_m)} dt + \int^\infty_0 e^{- \beta t} \Vert h(\cdot, t) \Vert^2_{W^{1, 2}(D_m)} dt < \infty.
\ee For $\lam$ such that $Re \, \lam \ge \beta$, denote by $w=w(x, \lam)$ the Laplace transform of $h=h(x, t)$ with respect to the $t$ variable i.e.
\[
w=w(x, \lam)=\int^\infty_0 e^{-\lam t} h(x, t) dt \equiv \mathcal{L}  h(x, \lam).
\]This is well defined since the integral is absolutely continuous and it is a $W^{1, 2}(D_m)$ valued, analytic function for $Re \, \lam \ge \beta$.
By \erf{Rh12}, $w=w(\cdot, \lam)$ is a $W^{1, 2}$ solution to the elliptic problem
\be
\lab{eqwlam}
\Delta w(x, \lam) - \lam w(x, \lam) = - h(x, 0) + \mathcal{L} R (x, \lam), \quad \text{in} \quad D_m; \quad \p_n w= 0 \quad \text{on} \quad \p D_m.
\ee By \erf{Rh12} again, the right hand side of this equation is in $L^2(D_m)$. If $\lam \ge \beta$ is real, from \cite{Gris} Chapter 4 again, we know that, for $s=3^-/2$,
\[
\Vert w(\cdot, \lam) \Vert_{W^{2, s}(D_m)} \le C \left(\Vert h(\cdot, 0) \Vert_{L^2(D_m)} +
\Vert  \mathcal{L} R (\cdot, \lam) \Vert_{L^2(D_m)} + (\lam +1) \Vert  w(\cdot, \lam) \Vert_{L^2(D_m)} \right)
\]Since $\lam \ge \beta>0$, we can use $w=w(x, \lam)$ as a test function in \erf{eqwlam} to deduce
\[
\Vert  w(\cdot, \lam) \Vert_{L^2(D_m)} \le 4 \lam^{-1/2} ( \Vert h(\cdot, 0) \Vert_{L^2(D_m)} +
\Vert  \mathcal{L} R (\cdot, \lam) \Vert_{L^2(D_m)}).
\]A combination of the previous two inequalities yields, for $\lam \ge \beta$,
\be
\Vert w(\cdot, \lam) \Vert_{W^{2, s}(D_m)} \le C (1+ \lam^{1/2} + \lam^{-1/2} ) \left(\Vert h(\cdot, 0) \Vert_{L^2(D_m)} +
\Vert  \mathcal{L} R (\cdot, \lam) \Vert_{L^2(D_m)} \right)<\infty.
\ee Since $J=-\p_{x_3} h$, this proves the assertion at the beginning of this sub-step.

{\it Step 6.2.} Our next task is to find a $L^2_t W^{1, 2}(D_m)$ solution to the following problem
\be
\lab{eqgbc}
\begin{cases}
\Delta g   - \p_t g + F
=0, \quad \text{on} \quad  D_m \times [0, T]\\
g=0, \quad \text{on} \quad \p^H D_m \times [0, T]; \qquad  \p_n g=0, \quad \text{on} \quad \p^V D_m \times [0, T]\\
g(x, 0) = - \p_{x_3} v_\th(x, 0)/r,
\end{cases}
\ee
where
\be
\lab{F=F123}
\al
F &\equiv -(b\cdot\nabla) J +\frac{2}{r}\p_r J + (\o_r \p_r + \o_3 \p_{x_3}) \frac{v_r}{r}\\
&= -(b\cdot\nabla) J +\frac{2}{r}\p_r J + [-\p_{x_3} v_\th \, \p_r + \frac{1}{r} \p_r(r v_\th) \, \p_{x_3}] \frac{v_r}{r}
\eal
\ee is the lower order terms in equation \erf{eqJbc}. Note also that $g$ and $J$ share the same initial value and boundary condition. Since $J=-\p_{x_3} v_\th/r$ and $v_\th$ is only in the energy space $\mathbf{E}$, we can only understand $\na J$ in pointwise sense except at the corners at this moment.
Let $\{ v^{(k)}_\th \}$ be a sequence of $L^2_t C^2_x(\bar{D_m})$  functions such that $ v^{(k)}_\th \to v_\th$ in $L^2_tW^{1, 2}_x(D_m)$. Such approximation sequence always exists since $D_m$ is a bounded Lipschitz domain and hence a $W^{1, 2}$ Sobolev extension domain. c.f. Calderon \cite{Cal}.
Since $D_m$ can be regarded as a rectangular polygon in the $rx_3$ plane, we can also require
\be
\lab{vkvwqi}
\Vert v^{(k)}_\th(\cdot, t) \Vert_{L^\infty(D_m)} \le \Vert v_\th(\cdot, t) \Vert_{L^\infty(D_m)}.
\ee The reason is the following. Fix a small positive number $\delta$. For any positive integer $k$, let
$D_{mk}$ be the polygon contained in $D_m$ such that the distance between $\p D_{mk}$ and $\p D_m$ is $\delta/k$.
Choose a one to one linear map $L_k$ that maps $D_{m}$ onto $D_{mk}$. For example, one can move the origin to a suitable point inside $D_{mk}$ and just shrink the horizontal and vertical variables in the order of $\delta/k$.  Then $ v^{(k)}_\th(r, x_3, t) =
v_\th(L_k(r, x_3), t)$ will satisfy all the above requirements.

Denote by $J^{(k)} = -\p_{x_3} v^{(k)}_\th/r$ and
\[
F_k= -(b\cdot\nabla) J^{(k)} +\frac{2}{r}\p_r J^{(k)} + [-\p_{x_3} v^{(k)}_\th \, \p_r + \frac{1}{r} \p_r(r v^{(k)}_\th) \, \p_{x_3}] \frac{v_r}{r}.
\]Since $F_k \in L^2_tL^2_x$, by standard Galerkin method, the problem
\be
\lab{eqgkbc}
\begin{cases}
\Delta g_k   - \p_t g_k + F_k
=0, \quad \text{on} \quad  D_m \times [0, T], \quad T>0,\\
g_k=0, \quad \text{on} \quad \p^H D_m \times [0, T]; \qquad  \p_n g_k=0, \quad \text{on} \quad \p^V D_m \times [0, T]\\
g_k(x, 0) = - \p_{x_3} v_\th(x, 0)/r,
\end{cases}
\ee has a unique solution in $L^2_t W^{1, 2}_{dn}(D_m) \cap \mathbf{E}$. Here $W^{1, 2}_{dn}(D_m)$ is the closure, under
$W^{1, 2}$ norm, of the set
\[
\{ u \in C^1(\bar{D_m}) \, | \, u(x)=0, \, x \in \p^H D_m;  \,
\p_n u(x)=0, \, x \in \p^V D_m \}.
\] Using $g_k$ as a test function in \erf{eqgkbc}, we obtain
\be
\lab{gknengl}
\al
 \int^T_0 \int_{D_m}& |\nabla g_k|^2 dxdt + \frac{1}{2} \int_{D_m} |g_k(x, T)|^2 dx
 -\frac{1}{2} \int_{D_m} |g(x, 0)|^2 dx =   \int^T_0 \int_{D_m} F_k g_k dxdt\\
&=  \int^T_0 \int_{D_m} \left\{-(b\cdot\nabla) J^{(k)} +\frac{2}{r}\p_r J^{(k)} + [-\p_{x_3} v^{(k)}_\th \, \p_r + \frac{1}{r} \p_r(r v^{(k)}_\th) \, \p_{x_3}] \frac{v_r}{r}\right\} g_k dxdt\\
&=  \int^T_0 \int_{D_m} \left\{(b\cdot\nabla g_k)  J^{(k)} + \underbrace{g_k \frac{1}{r}\p_r J^{(k)}}_{Y_2} + v^{(k)}_\th \, \p_r  \frac{v_r}{r} \p_{x_3} g_k  -  v^{(k)}_\th \, \p_{x_3} \frac{v_r}{r} \p_r  g_k \right\} dxdt.
\eal
\ee In the last step, we have used integration by parts except for the 2nd term $Y_2$, which are justified due to the following facts:
$div  \, b=0$, $\na \p_{x_3} v_r \in L^2_t L^2_x(D_m)$ from Step 5, $J^{(k)} \in L^2_t C^1_x(\bar{D_m})$ by assumption and by the boundary conditions for $v_r, v_3$ and $g_k$. There is also a cancellation of the second order derivatives of $v_r/r$:  $\p_r \p_{x_3} v_r/r$ which is in $L^2_tL^2_x$ by Step 5.
Using Cauchy-Schwarz inequality on the last line of \erf{gknengl} except for the 2nd term in the bracelet $Y_2$, we
arrive at
\[
\al
 \int^T_0 \int_{D_m}& |\nabla g_k|^2 dxdt +  \int_{D_m} |g_k(x, T)|^2 dx
 - \int_{D_m} |g(x, 0)|^2 dx \\
&\le C \int^T_0 \int_{D_m} \left( |v^{(k)}_\th|^2 \, | \na \frac{v_r}{r} |^2 + |b|^2 |J^{(k)}|^2 \right)
dxdt  + 4 \int^T_0 \int_{D_m}  g_k \frac{1}{r}\p_r J^{(k)} dxdt.
\eal
\]The last term requires a slightly different integration by parts since $g_k$ may not be $0$ on the vertical boundary $\p^V D_m$. We can treat it as follows
\[
\al
 4 \int^T_0 \int_{D_m}  g_k \frac{2}{r}\p_r J^{(k)} dxdt=-4 \int^T_0 \int_{D_m}  g_k \frac{1}{r}\p_r \frac{\p_{x_3} v^{(k)}_\th}{r} dxdt=4 \int^T_0 \int_{D_m}  \p_{x_3} g_k \frac{1}{r}\p_r \frac{ v^{(k)}_\th}{r} dxdt
\eal
\] since $g_k=0$ on the horizontal boundary $\p^H D_m$. Substituting this identity into the previous inequality and using Cauchy Schwarz again, we deduce
\[
\al
& \int^T_0 \int_{D_m} |\nabla g_k|^2 dxdt +  2 \int_{D_m} |g_k(x, T)|^2 dx
  \\
&\le C \int^T_0 \int_{D_m} \left( |v^{(k)}_\th|^2 \, | \na \frac{v_r}{r} |^2 + |b|^2 |J^{(k)}|^2
+ \frac{1}{r^2}  \left|\p_r \frac{ v^{(k)}_\th}{r}\right|^2 \right)
dxdt + 2 \int_{D_m} |g(x, 0)|^2 dx\\
&\le C \bar{C_m} \int^T_0 [\Vert v^2_\th(\cdot, t) + |b(\cdot, t)|^2 \Vert_{L^\infty(D_m)}+1]
\int_{D_m} \left(  | \na v_r |^2 + |\na \frac{v^{(k)}_\th}{r}|^2 + v^2_r
 \right)dxdt + 2 \Vert g(\cdot, 0) \Vert^2_{L^2(D_m)},
\eal
\] where we have used \erf{vkvwqi}. From here, using \erf{vthot},  \erf{vr3linft} and  the convergence
 of $v^{(k)}_\th$ to $v_\th$ in $L^2_tW^{1, 2}(D_m)$, we deduce
\be
\lab{energk}
\al
& \int^T_0 \int_{D_m} |\nabla g_k|^2 dxdt +  2 \int_{D_m} |g_k(x, T)|^2 dx\\
&\le C \bar{C_m} \int^T_0 e^{\alpha_0 t}
\int_{D_m} \left(  | \na v_r |^2 + |\na \frac{v_\th}{r}|^2 + v^2_r
 \right)dxdt + 2 \Vert g(\cdot, 0) \Vert^2_{L^2(D_m)}.
\eal
\ee
 This says that the energy norm of $g_k$ on $D_m \times [0, T]$ can be bounded from above by an uniform constant. Hence we can extract a subsequence converging to a function $g$ in weak $L^2_tW^{1, 2}(D_m)$ sense. This function $g$ is a solution to \erf{eqgbc} in $L^2_tW^{1, 2}(D_m)$. By standard interior and boundary regularity theory of the heat equation, we also know that $g$ solves \erf{eqgbc} except possibly at the corners of $D_m$. Moreover the term $b \na J$ is understood in the weak sense: for test functions
 $\phi \in L^2_t W^{1, 2}(D_m)$, the following holds
 \[
 \int^T_0\int_{D_m} (b \na J) \phi dxdt = - \int^T_0\int_{D_m} (b \na \phi) J dxdt.
 \]

{\it Step 6.3.} From \erf{energk}, the solution $g$ also satisfies the energy estimate
\be
\lab{energ1}
\al
& \int^T_0 \int_{D_m} |\nabla g|^2 dxdt +  2 \int_{D_m} |g(x, T)|^2 dx\\
&\le C \bar{C_m} \int^T_0 e^{\alpha_0 t}
\int_{D_m} \left(  | \na v_r |^2 + |\na \frac{v_\th}{r}|^2 + v^2_r
 \right)dxdt + 2 \Vert g(\cdot, 0) \Vert^2_{L^2(D_m)}.
\eal
\ee Now we prove $J=g$. Once this is done, since $J=-\p_{x_3} v_\th/r$, the energy bound \erf{energ1} will imply $|\na \p_{x_3} v_\th| \in L^2_tL^2_x$, completing the proof of Step 6.
Since $J$ satisfies \erf{eqJbc} and $g$ satisfies \erf{eqgbc}, the difference $J-g$ satisfies the following equation pointwise
\be
\lab{eqj-g}
\begin{cases}
\Delta (J-g)   - \p_t (J-g)
=0, \quad \text{on} \quad  D_m \times [0, \infty)\\
g=0, \quad \text{on} \quad \p^H D_m \times [0, \infty); \qquad  \p_n g=0, \quad \text{on} \quad \p^V D_m \times [0,\infty)\\
(J-g)(x, 0) =0.
\end{cases}
\ee

From \erf{energ1}, the energy bound for $h=v_\th/r$ \erf{dhl22} and standard energy bound for $v_r$,
we find that, for all $T>0$, there is a constant $C_1$ such that
\be
\lab{energ2}
\al
 \int^T_0 \int_{D_m} e^{-2 \alpha_0 t} (|\nabla g|^2 + g^2) dxdt
\le C_1 \bar{C_m}  \Vert v_0 \Vert^2_{W^{1,2}(D_m)}.
\eal
\ee This can be seen by taking $T$ to be positive integers and by splitting the time integral into the sum of integrals on $[j-1, j]$ with $j=1, 2, ... , T$.  Inequality \erf{energ2} and the assertion at Step 6.1 infer that the $W^{1, s}(D_m)$ valued Laplace transform for
$J-g$ is well defined if the phase variable $\lam$ satisfies $Re \, \lam \ge \beta$ for some suitable $\beta>0$. Here $s=3^-/2$ as before. Therefore the Laplace transform
$\mathcal{L}(J-g)(\cdot, \lam)$ is a $W^{1, s}(D_m)$ solution to the problem:
\[
\begin{cases}
\Delta [\mathcal{L }(J-g)](x, \lam)   - \lam \, \mathcal{L}(J-g)(x, \lam)
=0, \quad \text{on} \quad  D_m\\
\mathcal{L }(J-g)(x, \lam)=0, \quad \text{on} \quad \p^H D_m; \qquad  \p_n \mathcal{L }(J-g)(x, \lam)=0, \quad \text{on} \quad \p^V D_m.
\end{cases}
\] According to Proposition \ref{prunidn} below, if $\lam$ is real and sufficiently large, then $\mathcal{L}(J-g)(\cdot, \lam)=0$. But the Laplace transform is a $W^{1, s}(D_m)$ valued analytic function in $\lam$ with $Re \, \lam \ge \beta$.
This fact is again due to the assertion in Step 6.1 and \erf{energ2}. Thus $\mathcal{L}(J-g)(\cdot, \lam)=0$ for all these $\lam$, implying that
$J=g$. As mentioned, this completes the proof of Step 6.
\medskip

{\it Step 7.  $|\na^2 v_r| + |\na^2 v_3| \in L^2_tL^2_x$ and \erf{Onengl1} and \erf{Jnengl1}. }

Since it is already proven that $|\na \o_\th|, \, |\na \p_{x_3} v_r| \in L^2_tL^2_x$, and
$\o_\th=\p_{x_3} v_r - \p_r v_3$, it is clear that $|\na \p_r v_3| \in L^2_tL^2_x$.
Using the divergence free condition $\p_r v_r + \frac{1}{r} v_r + \p_{x_3} v_3=0$, one sees
\[
\al
\p^2_r v_r &= - \p_r( v_r/r) - \p_r \p_{x_3} v_3 \in L^2_tL^2_x,\\
\p^2_{x_3} v_3 &= - \p_{x_3}( v_r/r) -  \p_{x_3} \p_r v_r \in L^2_tL^2_x.
\eal
\] Hence $|\na^2 v_r| + |\na^2 v_3| \in L^2_tL^2_x$.

Finally \erf{Onengl2} is just \erf{Onengl1}. Also, since $J=g$, by taking $k \to \infty$ in \erf{gknengl}, we find
\be
\lab{Jnengl2}
\al
 \int^T_0 \int_{D_m}& |\nabla J|^2 dxdt + \frac{1}{2} \int_{D_m} |J(x, T)|^2 dx
 -\frac{1}{2} \int_{D_m} |J(x, 0)|^2 dx \\
&=  \int^T_0 \int_{D_m} \left\{(b\cdot\nabla J)  J + J \frac{2}{r}\p_r J + v_\th \, \p_r  \frac{v_r}{r} \p_{x_3} J  -  v_\th \, \p_{x_3} \frac{v_r}{r} \p_r  J\right\} dxdt.
\eal
\ee But
$
\int^T_0 \int_{D_m} (b\cdot\nabla J)  J dxdt= 0
$ since $div \, b=0$ and $v_r=0$ on $\p^V D_m$ and $v_3 =0$ on $\p^H D_m$. This integration by parts is justified because we already proved $b \in L^\infty \cap L^2_tW^{2, 2}(D_m)$ and $J \in L^2_tW^{1, 2}(D_m)$.
Therefore \erf{Jnengl1} follows from \erf{Jnengl2}.
This completes the proof of Step 7 and the proposition.
\qed

Now let us state and prove the lemma that was used in Step 4.1.

\begin{lemma}
\lab{leintx3w1s}
  Fixing $t$ and $s=3^-/2$, let $\tilde{v}_r=\tilde{v}_r(\cdot, t)$ be a $W^{2, s}(D_m)$ solution  to the problem
 \be
\lab{equrt2}
\begin{cases}
\left( \Delta-\frac{1}{r^2} \right) \tilde{v_r} =  \p_{x_3} \tilde{\o_\th}, \quad \text{in} \quad D_m,\\
\p_{x_3} \tilde{v_r} = 0, \quad \text{on} \quad \p^H D_m; \quad \tilde{v_r} = 0, \quad \text{on} \quad \p^V D_m,
\end{cases}
\ee where $\tilde{\o_\th}=\tilde{\o_\th}(\cdot, t) \in W^{1, 2}_0(D_m)$ and is smooth except at corners.

  Then the following holds
\[
\int_{D_m \cap \{r = const. \}} \tilde{v}_r(r, x_3, t) dx_3 =0.
\]
\end{lemma}

\proof For simplicity, we drop the $\, \tilde{} \,$  in the proof. For fixed $r>0$, denote by $L_r = \{ (r, x_3) \, | \, (r, x_3) \in D_m \}$ a vertical cross section of $D_m$, which is regarded as a polygon in the $rx_3$ plane. By standard interior and boundary regularity result, $v_r$ is smooth on $L_r$ if no corners of $D_m$ are contained in $L_r$. Since $\o_\th =0, \, \p_{x_3} v_r =0$ on $\p^H D_m$, we know that
\[
\int_{L_r} \p_{x_3} \o_\th \, dx_3 = 0, \quad \int_{L_r} \p^2_{x_3} v_r \, dx_3 = 0.
\]Using these and integrating along $L_r$,  the equation (c.f. \erf{equrt2}):
\[
\p^2_r v_r + \frac{1}{r} \p_r v_r + \p^2_{x_3} v_r - \frac{1}{r^2}  v_r = \left( \Delta-\frac{1}{r^2} \right) v_r =  \p_{x_3} \o_\th,
\]we find
\be
\lab{eqgr}
\p^2_r \int_{L_r} v_r dx_3 + \frac{1}{r} \p_r \int_{L_r} v_r dx_3   - \frac{1}{r^2} \int_{L_r} v_r dx_3 = 0.
\ee

Let us write $g_t(r) \equiv \int_{L_r} v_r(r, x_3, t) dx_3$.  Since $v_r(\cdot, t) \in  W^{2, s}(D_m)$ with $s=3^-/2$, by the trace theorem, $g_t(r)$ and $g'_t(r)$ are continuous function of $r$, and $g_t(r)$ is smooth unless $L_r$ contains the corners of $D_m$. For these $t$, we can solve the following ode from \erf{eqgr}
\[
g''_t(r) + \frac{1}{r} g'_t(r) -  \frac{1}{r^2} g_t(r)=0,
\]piecewise. Let $L_{r_1}, L_{r_2}, ..., L_{r_{m+1}}$ be those vertical segments containing corners of $D_m$. Then
\[
g_t(r)= c_i(t) r + \bar{c}_i(t) r^{-1}, \quad r_i < r \le r_{i+1}, \,  i=1, 2, ..., m,
\]where $c_i(t), \bar{c}_i(t)$ are independent of $r$. By continuity of $g_t(r)$ and $g'_t(r)$, for each $i=1, 2, ..., m-1$, the following hold:
\[
\al
 c_i(t) r_i + \bar{c}_i(t) r^{-1}_i &=  c_{i+1} (t) r_i + \bar{c}_{i+1} (t) r^{-1}_i, \\
 c_i(t)  - \bar{c}_i(t) r^{-2}_i &=  c_{i+1} (t) - \bar{c}_{i+1} (t) r^{-2}_i.
 \eal
\]Therefore $c_i(t)=c_{i+1}(t)$ and  $\bar{c}_i(t)=\bar{c}_{i+1}(t)$ so that
\[
g_t(r)= c_1(t) r + \bar{c}_1(t) r^{-1}.
\]According to the Navier boundary condition, $g_t(r_1)=g_t(r_{m+1})=0$ since $L_{r_1}$ is the left most vertical boundary of $D_m$ and $L_{r_{m+1}}$ is the right most. Hence $c_1(t) = \bar{c}_1(t) =0$ and
$g_t(r) \equiv 0$, proving the lemma.
\qed

\section{Proof of the main result}

\subsection{A priori bounds for $\Vert r v_\th \Vert_{L^\infty}$, $\Vert \o_\th/r \Vert_{L^\infty_t L^2_x}$, $\Vert \o_r/r \Vert_{L^\infty_t L^2_x}$, $\Vert v_r/r \Vert_{L^\infty_t L^6_x}$}

In this section, we prove a number a priori bounds for finite energy solutions to the following initial boundary value problem of  ASNS \eqref{eqasns} on
$D_m \times (0, T]$, $T>0$. These a priori bounds include the usual energy bound for $v$, and $L^2_tL^2_x$ bounds for $|\nabla \o_\th|$ and $|\nabla \o_r|$. The main point is to show that these bounds are independent of $m$, which will allow us to prove the main result in the next subsection after letting $m \to \infty$.

\begin{lemma}
\lab{lemuneng}  Let $v$ be a solution to Problem \ref{wtdm}. Then the following bounds are true.
\be
\lab{endm}
\bali
 (a). \quad & \int_{D_m} |v(x, T)|^2 dx \le \int_{D_m} |v(x, 0)|^2 dx;\\
 (b). \quad & 2 \int^T_0 \int_{D_m} |\p_{x_3} v_\th|^2 dxdt + 2 \int^T_0 \int_{D_m} |\p_r v_\th - \frac{1}{r} v_\th |^2 dxdt + \underbrace{2 \int^T_0 \int_{D_m} |\p_{x_3} v_r + \p_r v_3 |^2 dxdt}_{L_3} \\
 &\qquad + 4 \int^T_0 \int_{D_m} \left( |\p_r v_r |^2+\frac{v^2_r}{r^2} \right) dxdt
 + 4 \int^T_0 \int_{D_m} |\p_{x_3} v_3 |^2 dxdt
 \le \int_{D_m} |v(x, 0)|^2 dx;\\
& \text{if, in addition, $v$ comes from Proposition \ref{prsolS}, then}\\
 (c). \quad &   2 \int^T_0 \int_{D_m} \left( |\na v_r |^2+\frac{v^2_r}{r^2} + |\na v_3|^2
 + |r \, \p_r \frac{v_\th}{r} |^2 + |\p_{x_3} v_\th|^2 \right) dxdt
 \le \int_{D_m} |v(x, 0)|^2 dx.
 \eali
\ee
\end{lemma}

\proof Bounds (a) and (b) are immediate consequences of \eqref{en1s} and the formula for the strain tensor
\eqref{Sv}. So we just need to prove (c). Starting with part (b), we will just need to work out the third term $L_3$ on the left hand side of (b), which will be denoted by $L_3$.

Let us compute, using integration by parts,
\[
\al
\int^T_0 \int_{D_m} &\p_{x_3} v_r  \p_r v_3  dxdt = - \int^T_0 \int_{D_m}  v_r  \p_r  \p_{x_3} v_3  dxdt\\
&=\int^T_0 \int_{D_m}  v_r  \p_r (\p_r v_r + \frac{1}{r} v_r) dxdt = \int^T_0 \int \int \chi_{D_m}  v_r  \p_r (\p_r v_r + \frac{1}{r} v_r) r drdx_3dt\\
&= - \int^T_0 \int \int \chi_{D_m} \p_r v_r  \, (\p_r v_r + \frac{1}{r} v_r) r drdx_3dt -
\int^T_0 \int \int \chi_{D_m}  v_r  \, (\p_r v_r + \frac{1}{r} v_r)  drdx_3dt\\
&= - \int^T_0 \int \int \chi_{D_m} |\p_r v_r|^2   r drdx_3dt - 2 \int^T_0 \int \int \chi_{D_m} v_r \, \p_r v_r    drdx_3dt \\
&\qquad  - \int^T_0 \int \int \chi_{D_m} \frac{1}{r} |v_r|^2   drdx_3dt.
\eal
\] We comment that integration by parts is legal due to the property that $|\na^2 v_r|, \, |\na^2 v_3|
\in L^2_tL^2_x$ from Section 2 and the boundary condition $v_r =0$ on $\p^V D_m$ and $v_3 =0$ on $\p^H D_m$ in the trace sense for a.e.  $t$. Doing integration by parts and using $v_r =0$ on $\p^V D_m$ again, we see that  the 2nd from last term in the previous identity is $0$. Hence
\be
\lab{vrv31}
\int^T_0 \int_{D_m} \p_{x_3} v_r  \p_r v_3  dxdt =
- \int^T_0 \int_{D_m} |\p_r v_r|^2   dxdt - \int^T_0 \int_{D_m} \frac{1}{r^2} |v_r|^2   dxdt.
\ee Similarly, since $\p_{x_3} v_r=0$ on $\p^V D_m$ and $v_3=0$ on $\p^H_{D_m}$ in the trace sense, we deduce, after using the divergence free condition, that
\be
\lab{vrv32}
\al
\int^T_0 \int_{D_m} &\p_{x_3} v_r  \p_r v_3  dxdt = - \int^T_0 \int_{D_m}  v_3  \p_{x_3} \p_r (r v_r) \frac{1}{r}  dxdt\\
&=\int^T_0 \int_{D_m}  v_3  \p_{x_3} \p_{x_3} v_3  dxdt = -\int^T_0 \int_{D_m}     |\p_{x_3} v_3|^2 dxdt.
\eal
\ee Inequality (c) is derived after substituting \erf{vrv31} and \erf{vrv32} into $L_3$ in (b) separately (after expanding the square) and adding the resulting two inequalities.
\qed

The next lemma states that the line integral of $v_r$ in $x_3$ variable is $0$ for a.e. $t$.

\begin{lemma}
\lab{leintx3}
  Let $v$ be a solution to Problem \ref{wtdm} coming from Proposition \ref{prsolS}. Then
for a.e. $t$, the following holds
\[
\int_{D_m \cap \{r = const. \}} v_r(r, x_3, t) dx_3 =0.
\]
\end{lemma}

\proof This is done in Lemma \ref{leintx3w1s} in the previous section if $v_r(\cdot, t) \in W^{2, s}(D_m)$. Now that $v_r(\cdot, t) \in W^{2, 2}(D_m)$ for a. e. $t$, the conclusion follows.

\qed

The following lemma states that solutions to Problem \ref{wtdm} has the property that $v_\th$ is uniformly bounded near the right boundary of $D_m$ for all $m \ge 2$. Eventually, by the end of the section, we will have proven $v_\th$ is uniformly bounded in all $D_m$.

\begin{lemma}
\lab{levth1}
Let $v$ be a solution to Problem \ref{wtdm} with $m \ge 2$. There exists a positive constant $\bar{C}_1$, depending only on $\Vert v_\th(\cdot, 0) \Vert_{L^\infty(D_1)}$ and $\Vert v(\cdot, 0) \Vert_{L^2(D_m)}$ such that
\be
\lab{vthwqd1}
\Vert v_\th(\cdot, t) \Vert_{L^\infty(D_1 \cap \{ 3/4<r<1 \})} \le \bar{C}_1, \, \forall t >0.
\ee
\end{lemma}
\proof We will use Moser's iteration to prove the $L^\infty$ bound for $\Gamma = r v_\th$ in $D_1 \times [0, \infty)$. The Navier boundary condition for $v_\th$ on the right vertical boundary of $D_m$ and the lack of good bounds on $b=v_r e_r + v_3 e_3$ will  present a small obstacle. In order to proceed,  we will use a trace inequality and adopt the method in \cite{Zcmp04}, which can prove local boundedness of solutions under weaker than usual condition on the drift term $b$. Alternatively one can also work with
\erf{ibhth}, the equation for $h=v_\th/r$  , which satisfies the more friendly Neumann boundary condition. But one still needs to deal with the drift term and an additional potential term $-v_r/r$, using dimension reduction and energy bound \erf{endm} (c).

Let us start with the equation for $\Gamma$:
\be
\lab{eqvth2}
\begin{cases}
\Delta \Gamma - b \nabla \Gamma- \frac{2}{r} \p_r
\Gamma-\p_t \Gamma=0, \quad \text{in} \quad D_m \times (0, \infty)\\
\Gamma(x, 0) = r v_\th(x, 0), x \in D_m.
\end{cases}
\ee For any rational number $p>1$ in the form of $2k/l$ where $k, l$ are positive integers, we know that $\Gamma^p$ is a sub-solution, namely
\be
\lab{vthpsub}
\begin{cases}
\Delta \Gamma^p - b \nabla \Gamma^p- \frac{2}{r} \p_r
\Gamma^p-\p_t \Gamma^p \ge 0, \quad \text{in} \quad D_m \times (0, \infty)\\
\Gamma^p(x, 0) = (r v_\th(x, 0))^p, x \in D_m.
\end{cases}
\ee

Let $\phi=\phi(r)$ be a cut off function defined  on the interval $[0, 1]$ on the $r$ axis such that $0 \le \phi \le 1$,  $\phi(r)=1, \, r \in [3/4, 1]$; $\phi(r)=0, \, r \in [0, 1/2]$ and that $\Vert \phi'/\phi^{0.99} \Vert_{L^\infty}<\infty$. For $T \ge 1$, let $\eta=\eta(t)$ be a cut off function in time, supported in $[T-1, T]$ which will be specified later. Since  a nonuniform $L^\infty$ bound for $v_\th$ in Section 2 (Step 2 in the proof of Proposition \ref{prsolS}) is already proved, we can use $\Gamma^p (\phi \eta)^2$ as a test function on \erf{vthpsub} to deduce
\be
\lab{gamp3}
\al
&LS \equiv \int^T_{T-1}\int_{D_m} |\na (\Gamma^p \phi \eta)|^2 dxdt + \frac{1}{2} \int_{D_m} \Gamma^{2p} (\phi \eta)^2 (x, T) dx\\
& \le \underbrace{\int^T_0\int_{\p D_m} \p_n \Gamma^p \,  \Gamma^{p} (\phi \eta)^2 dSdt}_{R_1}
- \underbrace{\int^T_{T-1}\int_{D_m} b \na \Gamma^{p} \,  \Gamma^{p} (\phi \eta)^2 dxdt}_{R_2}
- \underbrace{\int^T_{T-1}\int_{D_m} \frac{2}{r} \p_r \Gamma^{p} \,  \Gamma^{p} (\phi \eta)^2 dxdt}_{R_3}\\
&\qquad + \int^T_{T-1}\int_{D_m} \Gamma^{2p} (|\na  (\phi \eta)|^2 +\eta' \eta) dxdt + \frac{1}{2} \int_{D_m} \Gamma^{2p} (\phi \eta)^2 (x, T-1) dx,
\eal
\ee where $n$ is the exterior normal of $\p D_m$. Next we will find bounds for $R_1$, $R_2$ and $R_3$.

To bound $R_1$, we use the boundary conditions $\p_{x_3} \Gamma = r \p_{x_3} v_\th =0$ on $\p^H D_m$.
Therefore
\[
\al
R_1 &= \int^T_{T-1} \int_{\p^V D_m} \p_n \Gamma^p \,  \Gamma^{p} (\phi \eta)^2 dS dt \\
&= - \sum^m_{j=1} \int^T_{T-1} \int_{L_j}  \p_r \Gamma^p \Gamma^{p} (\phi \eta)^2 r dx_3dt + \int^T_{T-1} \int_{L_{m+1}}  \p_r \Gamma^p \Gamma^{p} (\phi \eta)^2 r dx_3dt\\
&=\int^T_{T-1} \int_{L_{m+1}}  \p_r \Gamma^p \Gamma^{p} (\phi \eta)^2 r dx_3dt. \quad \text{(because $L_1, ..., L_m$ are cut off by $\phi$)}.
\eal
\] Here  $\p^V D_m = \cup^m_{j=1} L_j \cup L_{m+1}$, with $L_1,..., L_m$ being the vertical boundary segments to the left of $D_m$ and $L_{m+1}$ being the only vertical boundary segment to the right of $D_m$. Notice that $\p_n = -\p_r$ on $L_1, ..., L_m$ and $\p_n=\p_r$ on $L_{m+1}$. On $L_{m+1}$ and $L_j, \, j=1, 2, ..., m$, the Navier boundary condition reads
\be
\lab{prgambian}
\p_r \Gamma^p= \p_r ( r^p v^p_\th) = p r^{p-1} v^p_\th + r^p p v^{p-1}_\th \p_r v_\th
= 2 p r^{p-1} v^p_\th = \frac{2 p}{r} \Gamma^p.
\ee Consequently
\be
\lab{r1bian=}
R_1= 2p \int^T_{T-1} \int_{L_{m+1}}   \Gamma^{2p}  \eta^2  dx_3dt.
\ee

Next, via integration by parts, we see that
\be
\lab{r3bian=}
R_3 = \int^T_{T-1} \int_{L_{m+1}}   \Gamma^{2p}  \eta^2  dx_3dt
 - \int^T_{T-1} \int_{D_m } \Gamma^{2p} \p_r (\phi \eta)^2 \frac{1}{r} dxdt.
\ee

Using integration by parts again together with the boundary condition $v_r = 0$ on $\p^V D_m$, $v_3=0$ on $\p^H D_m$, we noticed that
\be
\lab{r2bian=}
R_2 = -\int^T_{T-1}\int_{D_m} v_r   \,  \Gamma^{2p} \phi \p_r \phi \, \eta^2 \, dxdt.
\ee The combination of \erf{r1bian=}, \erf{r2bian=}, \erf{r3bian=} and \erf{gamp3} yields
\be
\lab{ls2}
\al
LS & \le \underbrace{(2p-1) \int^T_{T-1} \int_{L_{m+1}}   \Gamma^{2p}  \eta^2  dx_3dt}_{I_1} + \underbrace{\int^T_{T-1}\int_{D_m} v_r   \,  \Gamma^{2p} \phi \p_r \phi \, \eta^2 \, dxdt}_{I_2} \\
&\qquad + \int^T_{T-1}\int_{D_m} \Gamma^{2p} \left[ |\na  (\phi \eta)|^2 + \eta' \eta + \p_r (\phi \eta)^2 \frac{1}{r}
\right] dxdt + \frac{1}{2} \int_{D_m} \Gamma^{2p} (\phi \eta)^2 (x, T-1) dx
\eal
\ee We need to absorb $I_1$ and $I_2$ by $LS$.

Since $\phi=1$ on $L_{m+1}$ and $\phi=0$ when $r=1/2$, we can compute
\[
\al
\int_{L_{m+1}}   \Gamma^{2p} dx_3& = \int_{L_{m+1}}   \Gamma^{2p} \phi^2 dx_3 - \int_{D_m \cap \{ r =1/2 \}}   \Gamma^{2p} \phi^2 dx_3\\
&=\int^1_0  \int^1_{1/2} \p_r \left( \Gamma^{2p} \phi^2 \right) dr dx_3 =
2 \int^1_0  \int^1_{1/2} (\Gamma^p \phi) \,  \p_r \left( \Gamma^p \phi \right) dr dx_3.
\eal
\]Therefore
\be
\lab{i142p}
\al
I_1 &\le 4 (2p-1) \int^T_{T-1}\int^1_0  \int^1_{1/2} (\Gamma^p \phi \eta) \, \sqrt{r}  \sqrt{r} \left|\p_r \left( \Gamma^p \phi \eta \right) \right| dr dx_3 dt\\
&\le 0.5 \int^T_{T-1}\int_{D_m} |\na (\Gamma^p \phi \eta)|^2 dxdt
+ 8(2p-1)^2 \int^T_{T-1}\int_{D_m}  (\Gamma^p \phi \eta)^2 dxdt
\eal
\ee

 Next we use the argument on p254-p255 in \cite{Zcmp04} with $m=4/3$, $b=v_r e_r+v_3 e_3$ and $\o=\Gamma^p$ and Corollary 2 there to deduce
\be
\lab{i242p}
\al
I_2 &\le 0.25 \int^T_{T-1}\int_{D_m} |\na (\Gamma^p \phi \eta)|^2 dxdt
+  C \Vert v_r \Vert^4_{L^\infty_tL^2_x(D_1)} \Vert \p_r \phi/\phi^{0.99} \Vert^4_\infty
\int^T_{T-1}\int_{D_1} |\Gamma^p \phi|^2 dxdt\\
&\le 0.25 \int^T_{T-1}\int_{D_m} |\na (\Gamma^p \phi \eta)|^2 dxdt
+  C \Vert v(\cdot,0) \Vert^4_{L^2(D_m)} \Vert \p_r \phi/\phi^{0.99} \Vert^4_\infty
\int^T_{T-1}\int_{D_1} |\Gamma^p \phi|^2 dxdt.
\eal
\ee In the last step, part (a) of the energy inequality \erf{endm} has been used. Also note that $v_3$ drops out since $\phi$ depends only on $r$.
Now we can plug \erf{i142p} and \erf{i242p} into \erf{ls2} to reach, since $\phi$ is supported in $D_1$,
\be
\lab{mositgam}
\al
&\int^T_{T-1}\int_{D_1} |\na (\Gamma^p \phi \eta)|^2 dxdt + 2 \int_{D_1} \Gamma^{2p} (\phi \eta)^2 (x, T) dx\\
&\le 32(2p-1)^2 \int^T_{T-1}\int_{D_1}  (\Gamma^p \phi \eta)^2 dxdt  + 4C \Vert v(\cdot,0) \Vert^4_{L^2(D_m)} \Vert \p_r \phi/\phi^{0.99} \Vert^4_\infty
\int^T_{T-1}\int_{D_1} |\Gamma^p \phi|^2 dxdt       \\
& \qquad + 4 \int^T_{T-1}\int_{D_1} \Gamma^{2p} \left[ |\na  (\phi \eta)|^2 + \eta' \eta + \p_r (\phi \eta)^2 \frac{1}{r}
\right] dxdt + 2 \int_{D_1} \Gamma^{2p} (\phi \eta)^2 (x, T-1) dx
\eal
\ee

Take $T=1$ and $\eta=1$ first. Then the last term in \erf{mositgam}  is bounded by 2 $\Vert v_\th(\cdot, 0) \Vert^{2p}_{L^\infty(D_1)}$.  By standard Moser's iteration from \erf{mositgam}, it is clear that \erf{vthwqd1} holds
for $t \in [0, 1]$. Note that only the standard $L^2$ Sobolev inequality is needed because  the spatial domain is $D_1$. For $T>1$, we take suitable sequences of $\phi$ and $\eta$ such that $\eta(T-1)=0$. Then the last term in  \erf{mositgam} drops out. Moser's iteration again tells us, for a positive number $q>1$,
\[
\Vert \Gamma(\cdot, T) \Vert^2_{L^\infty(D_1) \cap \{3/4<r<1\}} \le C [1+ \Vert v(\cdot,0) \Vert^4_{L^2(D_m)}]^q
\int^T_{T-1}\int_{D_1} \Gamma^2 dxdt.
\]Since $\Gamma=r v_\th$, this and the energy inequality imply \erf{vthwqd1} for all $T>1$, completing the proof of the lemma.
\qed

\begin{lemma}
\lab{levthj}
Let $v$ be a solution to Problem \ref{wtdm} and $\Gamma = r v_\th$. Then there is a positive constant $\bar{C}_1$, depending only on $\Vert v_\th(\cdot, 0) \Vert_{L^\infty(D_1)}$ and $\Vert v(\cdot, 0) \Vert_{L^2(D_m)}$ such that
\[
\Vert \Gamma(\cdot, t) \Vert_{L^\infty(D_m)} \le
\Vert \Gamma(\cdot, 0) \Vert_{L^\infty(D_m)} + \bar{C}_1, \quad \forall t>0.
\]
\end{lemma}
\proof
Recall from \erf{vthpsub} that for any positive even integer $p$, $\Gamma^p$ is a subsolution, namely
\be
\lab{eqvthp2}
\begin{cases}
\Delta \Gamma^p - b \nabla \Gamma^p- \frac{2}{r} \p_r
\Gamma^p-\p_t \Gamma^p \ge 0, \quad \text{in} \quad D_m \times (0, \infty)\\
\Gamma^p(x, 0) = (r v_\th(x, 0))^p, x \in D_m.
\end{cases}
\ee This implies
\be
\lab{dtgp}
\al
\p_t \int_{D_m} \Gamma^p(x, t) dx &\le \int_{D_m} \Delta \Gamma^p(x, t) dx - \int_{D_m} b \na \Gamma^p(x, t) dx - \int_{D_m} \frac{2}{r} \p_r \Gamma^p(x, t) dx\\
&\equiv T_1 + T_2 + T_3.
\eal
\ee Since we already proved a nonuniform $L^\infty$ bound for $v_\th$ in Section 2 (Step 2 in the proof of Proposition \ref{prsolS}), the above integrals are well defined.
We will find uniform bounds for $T_2$, $T_1$ and $T_3$ respectively.

Using the boundary condition $v_r = 0$ on $\p^V D_m$, $v_3=0$ on $\p^H D_m$, we see that
\be
\lab{t2dtgp}
T_2= \int_{D_m} div \, b  \, \Gamma^p(x, t) dx =0.
\ee To bound $T_1$, we use the boundary conditions $\p_{x_3} \Gamma = r \p_{x_3} v_\th =0$ on $\p^H D_m$.
Therefore
\be
\lab{t1lj}
T_1 = \int_{\p^V D_m} \p_n \Gamma^p dS = - \sum^m_{j=1} \int_{L_j}  \p_r \Gamma^p r dx_3 + \int_{L_{m+1}}  \p_r \Gamma^p r dx_3.
\ee Here $n$ is the outward normal of $\p D_m$ and $\p^V D_m = \cup^m_{j=1} L_j \cup L_{m+1}$, with $L_1,..., L_m$ being the vertical boundary segments to the left of $D_m$ and $L_{m+1}$ being the only vertical boundary segment to the right of $D_m$. Notice that $\p_n = -\p_r$ on $L_1, ..., L_m$ and $\p_n=\p_r$ on $L_{m+1}$, which accounts to the sign in \erf{t1lj}. By \erf{prgambian}, on $L_j$, $j=1, 2, ..., m+1$, the Navier boundary condition reads
$
\p_r \Gamma^p = \frac{2 p}{r} \Gamma^p.
$ Substituting this to \erf{t1lj} gives
\be
\lab{t1lj2}
T_1 = -2p  \sum^m_{j=1} \int_{L_j}   \Gamma^p  dx_3 + 2p \int_{L_{m+1}}  \Gamma^p  dx_3.
\ee Similarly,
\be
\lab{t3lj2}
T_3 = - \int_{D_m} \frac{2}{r} \p_r \Gamma^p(x, t) dx = 2 \sum^m_{j=1} \int_{L_j} \Gamma^p dx_3
- 2  \int_{L_{m+1}} \Gamma^p dx_3.
\ee Substituting   \erf{t2dtgp}, \erf{t1lj2}  and \erf{t3lj2} into \erf{dtgp}, we find, since $p$ is an even integer, that
\[
\al
\p_t \int_{D_m} \Gamma^p(x, t) dx &\le -2(p-1)  \sum^m_{j=1} \int_{L_j}   \Gamma^p  dx_3 + 2(p-1) \int_{L_{m+1}}  \Gamma^p  dx_3\\
&\le 2(p-1) \int_{L_{m+1}}  \Gamma^p  dx_3.
\eal
\]Therefore
\[
\al
\Vert \Gamma(\cdot, t) \Vert_{L^p(D_m)} &\le \left[ 2(p-1) \int^t_0 \int_{L_{m+1}} \Gamma^p dx_3ds +
\int_{D_m} \Gamma^p(x, 0) dx \right]^{1/p}\\
&\le  (2(p-1))^{1/p} \left[ \int^t_0 \int_{L_{m+1}} \Gamma^p dx_3ds \right]^{1/p}  +
\Vert \Gamma(\cdot, 0) \Vert_{L^p(D_m)}.
\eal
\]This shows, after letting $p \to \infty$, that
\[
\Vert \Gamma(\cdot, t) \Vert_{L^\infty(D_m)} \le \Vert v_\th \Vert_{L^\infty_tL^\infty(D_1 \cap \{ 3/4<r<1\})}
+\Vert \Gamma(\cdot, 0) \Vert_{L^\infty(D_m)},
\]which, together with Lemma \ref{levth1}, completes the proof of the lemma.
\qed

\begin{lemma}
\lab{levr/rOm}
  Let $v$ be a solution to Problem \ref{wtdm} coming from Proposition \ref{prsolS}. Then
for a.e. $t$, the following hold
\be
\lab{vr/rtoOm}
\Vert \nabla \frac{v_r}{r} \Vert_{L^2(D_m)} \le \Vert \frac{\o_\th}{r} \Vert_{L^2(D_m)}, \qquad
\Vert \nabla \p_{x_3} \frac{v_r}{r} \Vert_{L^2(D_m)} \le \Vert \partial_{x_3} \frac{\o_\th}{r} \Vert_{L^2(D_m)}.
\ee
\end{lemma}
\proof Again we choose those $t$ so that $v_r(\cdot, t) \in W^{2, 2}(D_m)$. According to Proposition \ref{prsolS} those $t$ form a set of full measure on the time axis. For simplicity, we write $\O=\o_\th/r$ and suppress the $t$ variable. It is known that the following relation is true:
\be
\lab{eqvr/r}
(\Delta + \frac{2}{r} \p_r ) \frac{v_r}{r}
=\p_{x_3} \O.
\ee  For example it can be derived from \erf{equrt} easily.

Using $\frac{v_r}{r}$ as a test function on \erf{eqvr/r} and using the boundary condition $v_r=0$ on $\p^V D_m$, $\p_{x_3} v_r=0$ on $\p^H D_m$ and $\O=0$ on $\p D_m$, we deduce
\[
\int_{D_m} \left| \na \frac{v_r}{r} \right|^2 dx = \int_{D_m} \O \, \p_{x_3} \frac{v_r}{r} dx
\]which yields the first inequality in \erf{vr/rtoOm} immediately.

To prove the 2nd inequality, we use $\p^2_{x_3} \frac{v_r}{r}$ as a test function on equation $\erf{eqvr/r}$. Using the notation $f=\frac{v_r}{r}$ for simplicity, we see that
\[
\int_{D_m}\p_{x_3} \O \,  \p^2_{x_3} f  dx =\int_{D_m}  \Delta f \, \p^2_{x_3} f dx + \int_{D_m} \frac{2}{r} \p_r f \, \p^2_{x_3} f dx
\equiv T_1 + T_2.
\]Using integration by parts and the boundary condition $\p_{x_3} f=0$ on $\p^V D_m$, we find
\[
T_2 =- \int_{D_m} \frac{2}{r} \p_r \p_{x_3} f \, \p_{x_3} f dx = - \int_{D_m} \frac{1}{r} \p_r (\p_{x_3} f)^2 \,  dx =0.
\]Next we work on $T_1$. Note that $\p_{x_3} f \in W^{1, 2}_0(D_m)$ due to the Navier boundary condition. Hence there exists a sequence of $C^2_0(\bar{D_m})$ functions,say $\{h_j\}$, such that $h_j$ converges to
$\p_{x_3} f$  in $W^{1, 2}_0(D_m)$ norm as $j \to \infty$. Thus, since $\na f_j$ vanishes on $\p D_m$,
\[
\al
T_1 = \lim_{j \to \infty} \int_{D_m}  \Delta f \, \p_{x_3} h_j dx = - \lim_{j \to \infty} \int_{D_m}  \na f \, \p_{x_3} \na h_j dx = \lim_{j \to \infty} \int_{D_m}  \na \p_{x_3} f \,  \na h_j dx
=\int_{D_m}  |\na \p_{x_3} f |^2 dx
\eal
\]The last three identities tell us
\[
\int_{D_m}  |\na \p_{x_3} f |^2 dx = \int_{D_m}\p_{x_3} \O \,  \p^2_{x_3} f  dx
\]which  proves the 2nd inequality in \eqref{vr/rtoOm}.

\qed

We need a lemma which states that a uniform Sobolev inequality holds for a subclass of $W^{1, 2}(D_m)$ functions on the domain $D_m$. Due to the thinness of $D_m$ near the $x_3$ axis, some extra condition is needed.
\begin{lemma}
\lab{sobdm}
Let $D_m$ be the domain in \eqref{dodm} and $f \in W^{1, 2}(D_m)$. Suppose either $f=0$ on $\p^H D_m$,  the horizontal boundary of $D_M$, or $\int_{D_m \cap \{r=const\}} f(r, x_3) dx_3=0$. Then there is a uniform constant $\mathbf{s}_0$, independent of $m$ or $f$, such that
\[
\Vert f \Vert_{L^6(D_m)} \le \mathbf{s}_0 \Vert \nabla f \Vert_{L^2(D_m)}.
\]
\end{lemma}

\proof  By definition $D_m = \cup^m_{i=1} S_i$ and we can write
\[
S_i = \cup^{m_i}_{k=1} S_{ik}
\]where $m_i$ is some positive integer and $S_{ik}$ are rectangles in the $rx_3$ plane with the following properties: the above union is non-overlapping and the ratio of the height and width of $S_{ik}$ is between $1$ and $2$. Note that $S_{ik}$ is not exactly an open rectangle, which may contain a vertical edge.
Denote the height of $S_{ik}$ by $h_i$. Since $h_i$ is between $1$ and $2$ times the width of $S_{ik}$, the standard Sobolev inequality holds on $S_{ik}$. Namely, there exists a uniform positive constant $C_0$ such that
\[
\left( \int_{S_{ik}} f^6 dx \right)^{1/3} \le C_0 \int_{S_{ik}} |\nabla f|^2 dx + \frac{C_0}{h^2_i}
\int_{S_{ik}} f^2 dx.
\]Due to the extra assumptions on $f$, we can apply the Poincar\'e inequality in the $x_3$ direction, which says, for a positive constant $C_1$,
\[
\int_{S_{ik}} f^2 dx \le C_1 h^2_i \int_{S_{ik}} |\p_{x_3} f|^2 dx.
\]A combination of the preceding two inequalities yields:
\[
\left( \int_{S_{ik}} f^6 dx \right)^{1/3} \le C_0 (1+C_1) \int_{S_{ik}} |\nabla f|^2 dx.
\]
Taking $\mathbf{s}_0= \sqrt{C_0(1+C_1)}$ and using the elementary inequality
\[
\left( \sum^j_{i=1} a_i \right)^{1/3} \le \sum^j_{i=1} a^{1/3} _i, \qquad a_i \ge 0, \quad j=1, 2, 3, ...,
\]we find
\[
\bali
\left( \int_{D_m} f^6 dx \right)^{1/3} &= \left( \sum^m_{i=1} \sum^{m_i}_{k=1} \int_{S_{ik}} f^6 dx \right)^{1/3} \le \sum^m_{i=1} \sum^{m_i}_{k=1} \left(  \int_{S_{ik}} f^6 dx \right)^{1/3}\\
& \le
\mathbf{s}^2_0 \sum^m_{i=1} \sum^{m_i}_{k=1}  \int_{S_{ik}} |\nabla f|^2 dx = \mathbf{s}^2_0 \int_{D_m} |\nabla f|^2 dx,
\eali
\]which proves the lemma. \qed

Now we turn to the key estimate on the vorticity.
Let $\o=\nabla \times v=\o_r e_r + \o_\theta e_\theta + \o_3 e_3$ be the vorticity.
Define
\be
\lab{defJO}
J= \frac{\omega_r}{r},  \quad \Omega=\frac{\omega_{\theta}}{r}.
\ee Then the triple $J, \O, \o_3$
 satisfy the system: for $b=v_r e_r + v_3 e_3$,
\begin{equation}
\label{eqjoo}
\begin{cases}
\Delta J  -(b\cdot\nabla) J +\frac{2}{r}\p_r J +
 (\o_r \p_r + \o_3 \p_{x_3}) \frac{v_r}{r} - \p_t J
=0,\\
\Delta \Omega -(b\cdot\nabla)\Omega+\frac{2}{r}\p_r
\Omega - \frac{2v_{\theta}}{r} J -\p_t \Omega=0,\\
\Delta \o_3-(b\cdot\nabla) \o_3+ \o_{r}\p_r v_3 + \o_3 \p_{x_3}
v_3 -\p_t \o_3=0.
\end{cases}
\end{equation}
These follow from direct computation based on the vorticity equation
\begin{align}
\lab{eqvort}
\begin{cases}
  \big (\Delta-\frac{1}{r^2} \big
)\omega_r-(b\cdot\nabla)\omega_r+\omega_r
\p_r v_r +\omega_3\p_{x_3} v_r -\p_t
\omega_r
=0,\\
   \big  (\Delta-\frac{1}{r^2}  \big
)\omega_{\theta}-(b\cdot\nabla)\omega_{\theta}+2\frac{v_{\theta}}
{r}\p_{x_3} v_{\theta}+\omega_{\theta}\frac{v_r}{r}-\p_t
\omega_{\theta}=0,\\
 \Delta\omega_3-(b\cdot\nabla)\omega_3+\omega_3\p_{x_3}
v_3+\omega_{r}\p_r v_3 -\p_t \omega_3=0,
\end{cases}
\end{align} and the relations
\be
\lab{w-v}
\o_r= -\p_{x_3} v_\theta, \quad \o_\theta= \p_{x_3} v_r-\p_r v_3, \quad \o_3 =
\p_r v_\theta +\frac{v_\theta}{r}.
\ee

\begin{lemma}
\lab{leotr/r}
Let $v$ be a solution to Problem \ref{wtdm} coming from Proposition \ref{prsolS}.
There exists a constant $\lam_0=\lam_0(\Vert r v_\th (\cdot, 0) \Vert_{L^\infty(D_m)}, \Vert v_\th(\cdot, 0) \Vert_{L^\infty(D_1)}, \, \Vert v(\cdot, 0) \Vert_{L^2(D_m)}, \beta) $, depending only on the initial value in terms of the stated quantities and independent of $m$ such that
\[
\Vert \O(\cdot, t) \Vert_{L^2(D_m)} + \Vert J(\cdot, t) \Vert_{L^2(D_m)} \le e^{\lam_0 t}
\left[\Vert \O(\cdot, 0) \Vert_{L^2(D_m)} + \Vert J(\cdot, 0) \Vert_{L^2(D_m)} \right].
\]Moreover
\[
\Vert \frac{v_r(\cdot, t)}{r} \Vert_{L^6(D_m)} \le \mathbf{s}_0 e^{\lam_0 t}
\left[\Vert \O(\cdot, 0) \Vert_{L^2(D_m)} + \Vert J(\cdot, 0) \Vert_{L^2(D_m)} \right],
\]where $\mathbf{s}_0$ is the Sobolev constant in Lemma \ref{sobdm}.
\end{lemma}
\proof
 The proof is based on inequalities \erf{Onengl1} for $\O$ and \erf{Jnengl1} for $J$ in Proposition \ref{prsolS}.
Using integration by parts on \erf{Onengl1} and the boundary condition that $\O=0$ on $\p D_m$, we find
\be
\lab{oo1}
\al
 \int^T_0 \int_{D_m} &|\na \O |^2 dxdt + \frac{1}{2} \int_{D_m} | \O(x, T) |^2 dx - \frac{1}{2} \int_{D_m} |\O(x, 0)|^2 dx\\
 &= -\int^T_0 \int_{D_m} \frac{v^2_\th}{r^2} \p_{x_3} \O  dxdt = - 2 \int^T_0 \int_{D_m} \frac{v_\th}{r}  \O \, J dxdt,
 \eal
\ee where we have used $J=- \p_{x_3} v_\th/r$.

To control the function $J$, we need to do integration by parts on the first integral on the right hand side of \erf{Jnengl1}.
\[
\al
\int^T_0 \int_{D_m} &J \frac{2}{r}\p_r J dxdt = \int^T_0 \int \int \chi_{D_m} \p_r J^2 drdx_3dt\\
&=-  \sum^m_{j=1} \int^T_0\int_{L_j}   J^2  dx_3dt + \int^T_0\int_{L_{m+1}}  J^2  dx_3dt.
\eal
\] Here as before  the vertical part of the boundary $\p^V D_m = \cup^m_{j=1} L_j \cup L_{m+1}$, with $L_1,..., L_m$ being the vertical boundary segments to the left of $D_m$ and $L_{m+1}$ being the only vertical boundary segment to the right of $D_m$. Since $r=1$ on $L_{m+1}$, the above identity implies
\[
\al
\int^T_0 \int_{D_m} &J \frac{2}{r}\p_r J dxdt \le 4 \int^T_0\int_{L_{m+1}}  J^2 \left(r-\frac{1}{2}\right)^2 dx_3dt = 4 \int^T_0\int^1_0 \int^1_{1/2} \p_r \left[  J^2 \left(r-\frac{1}{2}\right)^2 \right] dr dx_3dt\\
&= 4 \int^T_0\int^1_0 \int^1_{1/2}  \left[ 2 J \p_r J \left(r-\frac{1}{2}\right)^2 + 2 J^2 \left(r-\frac{1}{2}\right) \right] dr dx_3dt\\
&\le \frac{1}{2} \int^T_0\int_{D_m} |\na J|^2 dxdt + 10  \int^T_0\int_{D_m} J^2 dxdt.
\eal
\]Substituting this to the right hand side of \erf{Jnengl1}, we find
\[
\al
 \int^T_0 \int_{D_m}& |\nabla J|^2 dxdt +  \int_{D_m} |J(x, T)|^2 dx
 - \int_{D_m} |J(x, 0)|^2 dx \\
&\le  2 \int^T_0 \int_{D_m} \left\{ v_\th \, \p_r  \frac{v_r}{r} \p_{x_3} J  -  v_\th \, \p_{x_3} \frac{v_r}{r} \p_r  J\right\} dxdt + 1000  \int^T_0\int_{D_m} J^2 dxdt.
\eal
\]This infers, after applying Cauchy-Schwarz, that
\be
\lab{jj1}
\al
 \int^T_0 \int_{D_m}& |\nabla J|^2 dxdt +  2 \int_{D_m} |J(x, T)|^2 dx
 - 2\int_{D_m} |J(x, 0)|^2 dx \\
&\le  4 \int^T_0 \int_{D_m}  v^2_\th \, \left( \left|\p_r  \frac{v_r}{r} \right|^2 +
 \left|\p_{x_3}  \frac{v_r}{r} \right|^2 \right) dxdt + 2000  \int^T_0\int_{D_m} J^2 dxdt.
\eal
\ee A combination of \erf{oo1} and \erf{jj1} yields
\be
\lab{jowz}
\aligned
& \int_{D_m} \left(J^2+ \O^2 \right)  dx \bigg |^T_0
+  \frac{1}{2} \int^T_0 \int_{D_m} \left( | \nabla J|^2 +|\nabla \O|^2 \right)  dxdt\\
&\le  - \int^T_0 \int_{D_m} \frac{ 4 v_\theta}{r} \O J dxdt  + 2 \int^T_0 \int_{D_m}  v^2_\th \, \left( \left|\p_r  \frac{v_r}{r} \right|^2 +
 \left|\p_{x_3}  \frac{v_r}{r} \right|^2 \right) dxdt + 1000  \int^T_0\int_{D_m} J^2 dxdt.
\endaligned
\ee We will show that the first two terms on the right hand side of \erf{jowz} can be absorbed by the left hand side, modulo lower order terms.
According to Lemma \ref{levthj},
 there is a positive constant $\bar{C}_1=\bar{C}_1(\Vert v_\th(\cdot, 0) \Vert_{L^\infty(D_1)}, \, \Vert v(\cdot, 0) \Vert_{L^2(D_m)})$,
\[
\Vert \Gamma(\cdot, t) \Vert_{L^\infty(D_m)} \le
\Vert \Gamma(\cdot, 0) \Vert_{L^\infty(D_m)} + \bar{C}_1, \quad \forall t>0.
\] Since $\Gamma = r v_\th$, it follows that
\be
\lab{vth1/r}
| v_\th(r, x_3, t) | \le
 \frac{\Vert r v_\th (\cdot, 0) \Vert_{L^\infty(D_m)}+\bar{C}_1}{r} \equiv \frac{C_*(v_0)}{r}, \quad \forall t>0, \, (r, x_3) \in D_m.
\ee The key point is that $C_*(v_0)=\Vert r v_\th (\cdot, 0) \Vert_{L^\infty(D_m)} +\bar{C}_1(\Vert v_\th(\cdot, 0) \Vert_{L^\infty(D_1)}, \, \Vert v(\cdot, 0) \Vert_{L^2(D_m)})$ depends only on the initial velocity $v_0$.

Plugging \erf{vth1/r} into \erf{jowz}, we deduce
\be
\lab{jowz2}
\aligned
& \int_{D_m} \left(J^2+ \O^2 \right)  dx \bigg |^T_0
+  \frac{1}{2} \int^T_0 \int_{D_m} \left( | \nabla J|^2 +|\nabla \O|^2 \right)  dx dt\\
&\le 2 C_*(v_0) \int^T_0 \underbrace{\int_{D_m} \frac{1}{r^2} \O^2 dx}_{T1}  dt + 2 C_*(v_0) \int^T_0
\underbrace{\int_{D_m} \frac{1}{r^2} J^2 dx}_{T2} dt\\
&\qquad + 2 C^2_*(v_0)\int^T_0 \underbrace{\int_{D_m}  \frac{1}{r^2} \, \left( \left|\p_r  \frac{v_r}{r} \right|^2 +
 \left|\p_{x_3}  \frac{v_r}{r} \right|^2 \right) dx}_{T_3} dt + 1000  \int^T_0\int_{D_m} J^2 dxdt.
\endaligned
\ee Let us estimate $T_1, T_2$ and $T_3$ respectively.

Recall from \erf{dodm} that $D_m = \cup^m_{j=1} S_j$, the union of rectangles
\[
S_j = \{  (r, x_3)  \, | \,  2^{-j} \le r < 2^{-(j-1)}, \,    0< x_3 < 2^{-\beta (j-1)} \}
\]in the $rx_3$ plane. Denote by $h_j$ the height of $S_j$ (maximum value of $x_3$) and $r_j$ the width of $S_j$.  According to the boundary condition, the functions $\O=\o_\th/r$, $J=-\p_{x_3} v_\th/r$ and $\p_{x_3} v_r/r$ all vanishes on $\p^H D_m$ in point-wise or trace sense a.e. $t$. Recall from  Lemma \ref{leintx3}, $\int_{D_m \cap \{r = const. \}} v_r(r, x_3, t) dx_3 =0$, which infers
\[
\int_{D_m \cap \{r = const. \}} \p_r [v_r(r, x_3, t)/r] dx_3 =0.
\] Hence the one dimensional Poincar\'e inequality
\be
\lab{idpoin}
\int_{S_j \cap \{r = const. \}} f^2(x_3) dx_3 \le h^2_j \int_{S_j \cap \{r = const. \}} |\p_{x_3} f|^2(x_3) dx_3
\ee
 holds for all four of these functions
\[
f=f(x_3)= \O(r, x_3, t), \, \text{or} \, \,  J(r, x_3, t), \, \text{or} \, \,  \p_{x_3} v_r(r, x_3, t)/r, \, \text{or} \, \,  \p_r [v_r(r, x_3, t)/r] .
\]Let $j_0$ be a positive integer to be determined later.
We can estimate, using \erf{idpoin},
\be
\lab{t1j}
\al
T_1 &= \sum^m_{j=1}  \int_{S_j} \frac{1}{r^2} \O^2 dx = \sum^m_{j=j_0}  \int_{S_j} \frac{1}{r^2} \O^2  dx + \sum^{j_0-1}_{j=1}  \int_{S_j} \frac{1}{r^2} \O^2  dx\\
&\le  \sup_{j_0 \le j \le m} \frac{h^2_j}{r^2_j}  \int_{D_m}  |\p_{x_3} \O|^2  dx + 4^{j_0}  \int_{D_m} \O^2  dx
\eal
\ee Analogously
\be
\lab{t2j}
T_2 \le  \sup_{j_0 \le j \le m} \frac{h^2_j}{r^2_j}  \int_{D_m}  |\p_{x_3} J|^2  dx + 4^{j_0}  \int_{D_m} J^2  dx;
\ee
\be
\lab{t3j}
\al
T_3
&\le  \sup_{j_0 \le j \le m} \frac{h^2_j}{r^2_j}  \int_{D_m}
 \left|\na \p_{x_3}  \frac{v_r}{r} \right|^2  dx   + 4^{j_0}  \int_{D_m} \left|\na  \frac{v_r}{r} \right|^2  dx\\
&\le  \sup_{j_0 \le j \le m} \frac{h^2_j}{r^2_j}  \int_{D_m}
 \left| \p_{x_3}  \O \right|^2  dx   + 4^{j_0}  \int_{D_m} \O^2  dx
\eal
\ee where we have used Lemma \ref{levr/rOm}. Substituting \erf{t1j}, \erf{t2j} and \erf{t3j} into \ref{jowz2}, we find
\[
\aligned
& \int_{D_m} \left(J^2+ \O^2 \right)  dx \bigg |^T_0
+  \frac{1}{2} \int^T_0 \int_{D_m} \left( | \nabla J|^2 +|\nabla \O|^2 \right)  dx dt\\
&\le [2 C_*(v_0) + 2 C^2_*(v_0)] \sup_{j_0 \le j \le m} \frac{h^2_j}{r^2_j} \int^T_0 \int_{D_m} |\p_{x_3} \O|^2 dx  dt + 2 C_*(v_0) \sup_{j_0 \le j \le m} \frac{h^2_j}{r^2_j} \int^T_0
\int_{D_m} |\p_{x_3} J|^2 dxdt\\
&\qquad + [2 C_*(v_0) + 2 C^2_*(v_0)] 4^{j_0}\int^T_0 \int_{D_m}  \O^2 +
  [2 C_*(v_0) 4^{j_0} + 1000]  \int^T_0\int_{D_m} J^2 dxdt.
\endaligned
\]By our choice of $\beta \in (1, 1.1]$ in the construction of the domain $D_m$, we can pick the smallest integer $j_0$ such that
\[
j_0 \ge (\beta-1)^{-1} \ln(4^{2.1} [2 C_*(v_0) + 2 C^2_*(v_0)])/\ln 4.
\]Then
\[
[2 C_*(v_0) + 2 C^2_*(v_0)] \sup_{j_0 \le j \le m} \frac{h^2_j}{r^2_j} \le 1/4.
\]Consequently
\[
\aligned
& \int_{D_m} \left(J^2+ \O^2 \right)  dx \bigg |^T_0
+  \frac{1}{4} \int^T_0 \int_{D_m} \left( | \nabla J|^2 +|\nabla \O|^2 \right)  dx dt\\
&\le
 [2 C_*(v_0) + 2 C^2_*(v_0)] 4^{j_0}\int^T_0 \int_{D_m}  \O^2 +
  [2 C_*(v_0) 4^{j_0} + 1000]  \int^T_0\int_{D_m} J^2 dxdt.
\endaligned
\]Taking
\[
\lam_0=[2 C_*(v_0) + 2 C^2_*(v_0)] 4^{j_0} + 1000.
\]and applying Gronwall's inequality, we have verified the first inequality in the statement of the lemma.
From the definition of $C_*(v_0)$ in \erf{vth1/r}, we see that $\lam_0$ depends only on the 3 stated quantities about the initial value and the parameter $\beta$ for $D_m$.
 The second one follows from the first one, Lemma \ref{levr/rOm} and  Lemma \ref{sobdm}, which is applicable since the line integral in the $x_3$ direction for $v_r/r$ is $0$ by Lemma \ref{leintx3}.  \qed

\subsection{A priori bounds for $\Vert  v_\th \Vert_{L^\infty_t L^\infty_x}$, $\Vert  v_r \Vert_{L^\infty_t L^\infty_x}$,
$\Vert  v_3 \Vert_{L^\infty_t L^\infty_x}$}

Based on the a priori bounds from the previous subsection, we will prove an a priori bound for $\Vert v \Vert_{L^\infty_t L^\infty_x}$. Let us start with $v_\th$. We will apply Moser's iteration on the equation for $v_\th$ in
\erf{eqasns}:
\be
\lab{eqvth2}
\begin{cases}
   \big (\Delta-\frac{1}{r^2} \big )
v_\th-(v_r \p_r + v_3 \p_{x_3})v_\th - \frac{v_r v_\th}{r}-\p_t  v_\th=0,\\
v_\th(x, 0)=(v_0)_\th(x).
\end{cases}
\ee Intuitively, we regard $v_r/r$ as a potential function for $v_\th$. According to Lemma \ref{leotr/r},
$v_r/r \in L^\infty_t L^6_x$, which makes it a subcritical potential. Hence there is a chance to prove $L^\infty$ bound for $v_\th$. There is an extra hurdle to face in though, due to the lack of uniform Sobolev inequality for functions like $v_\th$ on $D_m$.  Instead, we will use a Sobolev inequality with a smaller power and extra weight $1/r^2$, which is weaker than the standard Sobolev inequality but is sufficient to carry out Moser's iteration. The constant in this weaker inequality is uniform for all $D_m$. The exponent $13/5$ may be improved, which will result in an improvement in the range of $\beta$ in the definition of the domain in the main theorem.

\begin{lemma}
\lab{lewsob}
Let $f \in W^{1, 2}(D_m)$, $m=1, 2, 3, ...$. There exists a uniform constant $C_s$ such that
\[
\left( \int_{D_m} |f|^{13/5} dx \right)^{10/13} \le C_s \int_{D_m} \left( |\na f |^2 + \frac{1}{r^2} f^2 \right) dx.
\]
\end{lemma}
\proof By definition again $D_m= \bigcup^m_{j=1} S_j$ where $S_j$ can be regarded as a rectangle in the $rx_3$ plane with width $r_j = 2^{-j}$ and height $h_j= 2^{-\beta (j-1)}$ with $\beta \in (1, 1.1]$. By the standard Sobolev inequality on $S_j$, there exists a uniform positive constant $C_0$ such that
\[
\left( \int_{S_j} f^6 dx \right)^{1/3} \le C_0 \int_{S_j} |\nabla f|^2 dx + \frac{C_0}{h^2_j}
\int_{S_j} f^2 dx.
\]Observe that $h_j$ is much smaller than $r_j$ for large $j$ since $\beta \in (1, 1.1]$.
This and H\"older inequality imply
\be
\lab{13/5}
\al
\left( \int_{S_j} |f|^{13/5} dx \right)^{10/13} &\le \left( \int_{S_j} f^{(13/5) (30/13)} dx \right)^{1/3} \left( \int_{S_j}  dx \right)^{17/39} \le C \left( \int_{S_j} f^6 dx \right)^{1/3} (r^2_j h_j)^{17/39}\\
&\le C \left[ C_0 \int_{S_j} |\nabla f|^2 dx + \frac{C_0}{h^2_j}
\int_{S_j} f^2 dx \right] (r^2_j h_j)^{17/39}.
\eal
\ee By the choice that $\beta \in (1, 1.1]$, we have $r^{103/61}_j \le h_j = o(r_j)$ for $j \ge 100$. Therefore
\[
(r^2_j h_j)^{17/39}/ h^{2}_j = r^{34/39}_j/h^{61/39}_j \le r^{-69/39}_j \le r^{-2}_j.
\]Substituting this to \erf{13/5}, we find
\[
\left( \int_{S_j} |f|^{13/5} dx \right)^{10/13} \le C \left[ C_0 \int_{S_j} |\nabla f|^2 dx + \frac{C_0}{r^2_j}
\int_{S_j} f^2 dx \right]
\]Summing up, we deduce
\[
\left( \int_{D_m} |f|^{13/5} dx \right)^{10/13} \le \sum^m_{j=1} \left(  \int_{S_j} |f|^{13/5} dx \right)^{10/13} \le C C_0 \left[  \int_{D_m} |\nabla f|^2 dx +
\int_{D_m} \frac{1}{r^2} f^2 dx \right]
\] The lemma is proven by taking $C_s = C C_0$. \qed

The next lemma provides a $L^\infty$ bound for $v_\th$.

\begin{lemma}
\lab{levthinft}
Let $v$ be a solution to Problem \ref{wtdm} coming from Proposition \ref{prsolS}. Then, there is an absolute constant $C$ such that
\be
\lab{vtht>0}
\al
\Vert v_\th(\cdot, T) \Vert_{L^\infty(D_m)} \le
\begin{cases}
C \left[ 1 +C_s^{13/5}  \Vert \frac{v_r}{r} \Vert_{L^\infty_{0<t<T}L^6(D_m)}^{18/5}
 \right]^{8/3} |D_m|^{1/2} \, \Vert (v_0)_\th \Vert_{L^\infty(D_m)}, \quad T \le 1,\\
  C \left[ 1 +C_s^{13/5}  \Vert \frac{v_r}{r} \Vert_{L^\infty_{T-0.75<t<T}L^6(D_m)}^{18/5}
 \right]^{8/3} \left(\int_{D_m} |(v_0)_\th|^2 dx\right)^{1/2}, \quad T \ge 1.
\end{cases}
\eal
\ee Here $C_s$ is the Sobolev constant in Lemma \ref{lewsob}.
\end{lemma}
\proof
We will work on the equation \erf{eqvth2}.
Except for the term $-v_r v_\th/r$, the treatment of other terms are similar to those in Lemma \ref{levth1}.

For any rational number $p>1$ in the form of $2k/l$ where $k, l$ are positive integers, we know that $ v_\th^p$ is a sub-solution, namely
\be
\lab{vthpsub2}
\begin{cases}
\Delta v^p_\th - \frac{1}{r^2} v_\th^p - b \nabla v^p_\th -p \frac{v_r}{r} v^p_\th -\p_t v^p_\th \ge 0, \quad \text{in} \quad D_m \times (0, \infty)\\
v^p_\th(x, 0) = [(v_0)_\th(x)]^p, x \in D_m.
\end{cases}
\ee

 For $T \ge 1$, let $\eta=\eta(t)$ be a cut off function in time, supported in $[T-1, T]$ which will be specified later. Since  a nonuniform $L^\infty$ bound for $v_\th$ is already proven (Step 2 in the proof of Proposition \ref{prsolS}), we can use $v_\th^p  \eta^2$ as a test function on \erf{vthpsub2} to deduce
\be
\lab{gamp4}
\al
&LS \equiv \int^T_{T-1}\int_{D_m} |\na (v_\th^p  \eta)|^2 dxdt + \int^T_{T-1}\int_{D_m} \frac{1}{r^2} v_\th^{2p} \eta^2 (x, t) dxdt + \frac{1}{2} \int_{D_m} v_\th^{2p} \eta^2 (x, T) dx\\
& \le \underbrace{\int^T_{T-1}\int_{\p D_m} \p_n v_\th^p \,  v_\th^{p} \eta^2 dSdt}_{R_1}
- \underbrace{\int^T_{T-1}\int_{D_m} b \na v_\th^{p} \,  v_\th^{p}  \eta^2 dxdt}_{R_2}
\\
&\qquad - \underbrace{p \int^T_{T-1}\int_{D_m} \frac{v_r}{r} \,  v_\th^{2p} \eta^2 dxdt}_{R_3} +
\int^T_{T-1}\int_{D_m}  v_\th^{2p} \eta' \eta dxdt + \frac{1}{2} \int_{D_m} v_\th^{2p} \eta^2 (x, T-1) dx,
\eal
\ee where $n$ is the exterior normal of $\p D_m$. Next we will find bounds for $R_1$, $R_2$ and $R_3$.

First, using integration by parts  together with the boundary condition $v_r = 0$ on $\p^V D_m$, $v_3=0$ on $\p^H D_m$, we see that
\be
\lab{r2bian2=}
R_2 = 0.
\ee

To bound $R_1$, we use the boundary conditions $\p_{x_3} v_\th  =0$ on $\p^H D_m$.
Therefore
\[
\al
R_1 &= \int^T_{T-1} \int_{\p^V D_m} \p_n v_\th^p \,  v_\th^{p}  \eta^2 dS dt \\
&= - \sum^m_{j=1} \int^T_{T-1} \int_{L_j}  \p_r v_\th^p v_\th^{p}  \eta^2 r dx_3dt + \int^T_{T-1} \int_{L_{m+1}}  \p_r v_\th^p v_\th^{p} \eta^2 r dx_3dt.
\eal
\] Here again $\p^V D_m = \cup^m_{j=1} L_j \cup L_{m+1}$, with $L_1,..., L_m$ being the vertical boundary segments to the left of $D_m$ and $L_{m+1}$ being the only vertical boundary segment to the right of $D_m$. Notice that $\p_n = -\p_r$ on $L_1, ..., L_m$ and $\p_n=\p_r$ on $L_{m+1}$. On $L_{m+1}$ and $L_j, \, i=1, 2, ..., m$, the Navier boundary condition reads
\be
\lab{vthpbian}
\p_r v_\th^p=    p v^{p-1}_\th \p_r v_\th
=  p r^{-1} v^p_\th.
\ee Consequently
\be
\lab{r1bian2=}
\al
R_1&= p \int^T_{T-1} \int_{L_{m+1}}   v_\th^{2p}  \eta^2  dx_3dt - p \sum^m_{j=1} \int^T_{T-1} \int_{L_j}  v_\th^{2p}    \eta^2  dx_3dt\\
&\le p \int^T_{T-1} \int_{L_{m+1}}   v_\th^{2p}  \eta^2  dx_3dt.
\eal
\ee

 For simplicity we write $f\equiv v_\th^p$. Since $r=1$ on $L_{m+1}$, we see that
\[
\al
&R_1 \le 4 p \int^T_{T-1}\int_{L_{m+1}}  f^2 \left(r-\frac{1}{2}\right)^2 \eta^2 dx_3dt = 4 p \int^T_{T-1}\int^1_0 \int^1_{1/2} \p_r \left[  f^2 \left(r-\frac{1}{2}\right)^2 \right] dr dx_3 \eta^2 dt\\
&= 4 p \int^T_{T-1}\int^1_0 \int^1_{1/2}  \left[ 2 f \p_r f \left(r-\frac{1}{2}\right)^2 + 2 f^2 \left(r-\frac{1}{2}\right) \right] \eta^2 dr dx_3dt.
\eal
\]This implies
\be
\lab{r131}
R_1 \le \frac{1}{2} \int^T_{T-1}\int_{D_m} |\na f|^2 \eta^2 dxdt + 10 p^2 \int^T_{T-1}\int_{D_m} f^2 \eta^2 dxdt.
\ee

To estimate $R_3$, let us notice that
\[
\al
\int_{D_m} &\frac{|v_r|}{r} \,  v_\th^{2p}  dx = \int_{D_m} \frac{|v_r|}{r} \,  f^2  dx =
\int_{D_m} \frac{|v_r|}{r} \,  f^{13/9}  f^{5/9}  dx \\
&\le \left(\int_{D_m} \left| \frac{v_r}{r} \right|^6 dx \right)^{1/6} \,
\left( \int_{D_m} \,  f^{13/5} dx \right)^{5/9} \, \left( \int_{D_m} \,  f^2 dx \right)^{5/18}  \quad \text{(by H\"older inequality)}\\
&\le \e^{18/13} \left( \int_{D_m} \,  f^{13/5} dx \right)^{10/13} + \e^{-18/5} \Vert \frac{v_r}{r}(\cdot, t) \Vert_{L^6(D_m)}^{18/5} \, \int_{D_m} \,  f^2 dx\\
&\le  C_s  \e^{18/13} \int_{D_m} \left( |\na f |^2 + \frac{1}{r^2} f^2 \right) dx
 + \e^{-18/5} \Vert \frac{v_r}{r}(\cdot, t) \Vert_{L^6(D_m)}^{18/5} \, \int_{D_m} \,  f^2 dx
 \quad \text{(by Lemma \ref{lewsob}).}
\eal
\]Thus
\[
|R_3| \le  p C_s \e^{18/13} \int^T_{T-1} \int_{D_m} \left( |\na f |^2 + \frac{1}{r^2} f^2 \right) dx + p \e^{-18/5} \Vert \frac{v_r}{r} \Vert_{L^\infty_tL^6(D_m)}^{18/5} \, \int^T_{T-1} \int_{D_m} \,  f^2 dx.
\]From this inequality, we choose $\e$ so that $ p C_s \e^{18/13} = 1/8$ or $p \e^{-18/5} = ( 8 C_s)^{13/5} p^{18/5}$ to obtain
\be
\lab{vthr3}
|R_3| \le  \frac{1}{8} \int^T_{T-1} \int_{D_m} \left( |\na f |^2 + \frac{1}{r^2} f^2 \right) \eta^2 dxdt +
  (8 C_s)^{13/5} p^{18/5} \Vert \frac{v_r}{r} \Vert_{L^\infty_tL^6(D_m)}^{18/5} \, \int^T_{T-1} \int_{D_m} \,  f^2 \eta^2dxdt.
\ee

The combination of \erf{vthr3}, \erf{r131}, \erf{r2bian2=},  and \erf{gamp4} yields
\be
\lab{vth2p}
\al
&\frac{1}{4} \int^T_{T-1}\int_{D_m} |\na (f  \eta)|^2 dxdt + \frac{1}{4} \int^T_{T-1}\int_{D_m} \frac{1}{r^2} f^2 \eta^2 (x, t) dxdt + \frac{1}{2} \int_{D_m} f^2 \eta^2 (x, T) dx\\
&\le \left[ 10 p^2+(8 C_s)^{13/5} p^{18/5} \Vert \frac{v_r}{r} \Vert_{L^\infty_tL^6(D_m)}^{18/5}\right] \, \int^T_{T-1} \int_{D_m} \,  f^2 \eta^2 dxdt + \\
&\qquad + \int^T_{T-1}\int_{D_m}  f^2 \eta' \eta dxdt + \frac{1}{2} \int_{D_m} f^2 \eta^2 (x, T-1) dx.
\eal
\ee Next
\[
\al
\int^T_{T-1}\int_{D_m}& (f \eta)^{32/13} dxdt = \int^T_{T-1}\int_{D_m} (f \eta)^2 (f \eta)^{6/13} dxdt\\
&\le \int^T_{T-1} \left( \int_{D_m} |f \eta|^{13/5} dx \right)^{10/13} \left( \int_{D_m} |f \eta|^{(6/13)(13/3)} dx \right)^{3/13} dt \\
&\le \sup_{t \in [T-1, T]} \left( \int_{D_m} |f \eta|^2(x, t) dx \right)^{3/13} \, \int^T_{T-1} \left( \int_{D_m} |f \eta|^{13/5} dx \right)^{10/13} dt \\
&\le C_s \sup_{t \in [T-1, T]} \left( \int_{D_m} |f \eta|^2(x, t) dx \right)^{3/13} \,
\int^T_{T-1} \int_{D_m} \left[|\na (f  \eta)|^2 + \frac{1}{r^2} f^2 \eta^2\right]  dxdt
\eal
\]
where we have used Lemma \ref{lewsob} again. Plugging \erf{vth2p} to the right hand side in the preceding inequality gives, since $f=v^p_\th$,
\be
\lab{vthmoseri}
\al
&\left(\int^T_{T-1}\int_{D_m} (v^p_\th \eta)^{32/13} dxdt\right)^{13/16} \le 4\left[ 10 p^2+(8 C_s)^{13/5} p^{18/5} \Vert \frac{v_r}{r} \Vert_{L^\infty_tL^6(D_m)}^{18/5} +
\sup |\eta'| \right] \\
&\qquad \times \, \left[ \int^T_{T-1} \int_{D_m} \,  v^{2p}_\th \eta dxdt + 2 \int_{D_m} v^{2p}_\th \eta^2 (x, T-1) dx \right].
\eal
\ee
Take $T=1$ and $\eta=1$ first. Then the integrand in the last term in the above  is bounded by 2 $\Vert (v_0)_\th(\cdot) \Vert^{2p}_{L^\infty(D_m)}$.  Writing $\alpha_0=\Vert (v_0)_\th(\cdot) \Vert_{L^\infty(D_m)}$, then \erf{vthmoseri} implies
\be
\lab{vthmoseri1}
\al
&\left(\int^1_0\int_{D_m} [(|v_\th| \vee \alpha_0)^{2p}]^{16/13} dxdt\right)^{13/16} \le C \left[ 10 p^2+(8 C_s)^{13/5} p^{18/5} \Vert \frac{v_r}{r} \Vert_{L^\infty_tL^6(D_m)}^{18/5}
 \right] \\
&\qquad \times \,  \int^1_0 \int_{D_m} \,  (|v_\th| \vee \alpha_0)^{2p} dxdt.
\eal
\ee Here $|v_\th| \vee \alpha_0 =\max \{|v_\th|(x, t), \alpha_0 \}$ and $C$ is an absolute constant.
 By  Moser's iteration from \erf{vthmoseri1} with $p=(16/13)^j$, $j=0, 1, 2, ...$, we deduce
\be
\lab{vtht01}
\Vert v_\th \Vert_{L^\infty(D_m \times [0, 1])} \le  C \left[ 1 +C_s^{13/5}  \Vert \frac{v_r}{r} \Vert_{L^\infty_{0<t<1}L^6(D_m)}^{18/5}
 \right]^{8/3} \left(\int^1_0\int_{D_m} (|v_\th|+\alpha_0)^2 dxdt\right)^{1/2},
\ee where $C$ is another absolute constant.
 For $T>1$, we take suitable sequences of $\eta$ such that $\eta(T-1)=0$. Then the last term in  \erf{vthmoseri} drops out. Moser's iteration again tells us,
\be
\lab{vtht>1}
\al
\Vert v_\th \Vert_{L^\infty(D_m \times [T-0.5, T])} &\le  C \left[ 1 +C_s^{13/5}  \Vert \frac{v_r}{r} \Vert_{L^\infty_{0<t<T}L^6(D_m)}^{18/5}
 \right]^{8/3} \left(\int^T_{T-0.75}\int_{D_m} |v_\th|^2 dxdt\right)^{1/2},\\
& \le C \left[ 1 +C_s^{13/5}  \Vert \frac{v_r}{r} \Vert_{L^\infty_{T-0.75<t<T}L^6(D_m)}^{18/5}
 \right]^{8/3} \left(\int_{D_m} |(v_0)_\th|^2 dx\right)^{1/2}
\eal
\ee where $C$ is an absolute constant and the energy inequality is used.  The lemma follows from
\erf{vtht01} and \erf{vtht>1}.
\qed

The next lemma provides a $L^\infty$ bound for $\o_\th$, which is needed to bound $v_r$ and $v_3$.

\begin{lemma}
\lab{leothinft}
Let $v$ be a solution to Problem \ref{wtdm} coming from Proposition \ref{prsolS}. Then, there is an absolute constant $C$ such that
\be
\lab{otht>0}
\al
\Vert \o_\th(\cdot, T) \Vert_{L^\infty(D_m)} \le
\begin{cases}
 C \Lambda^{5/4}_0 \left(\int_{D_m} (|v_0|+\Vert (\o_0)_\th \Vert_{L^\infty(D_m)}+1)^2 dx\right)^{1/2}, \quad T \le 1,\\
   C (\Lambda+1)^{5/4}_0 \left(\int_{D_m} (|v_0|+1)^2 dx\right)^{1/2}, \quad T \ge 1,
\end{cases}
\eal
\ee where
\[
\Lambda_0 \equiv C_s^{13/5}  \Vert \frac{v_r}{r} \Vert_{L^\infty_{T-1<t<T}L^6(D_m)}^{18/5} +   \Vert v_\th \Vert^4_{L^\infty(D_m \times [T-1, T])} (\Vert \frac{\p_{x_3} v_\th}{r} \Vert^4_{L^\infty_{T-1<t<T} L^2_x}
+1)
\]Here $C_s$ is the Sobolev constant in Lemma \ref{lewsob}.
\end{lemma}
\proof The proof is based on the equation for $\o_\th$ which was already stated as in \erf{eqotil2} e.g.

\be
\lab{eqoth2}
\begin{cases}
 (\Delta-\frac{1}{r^2}  \big
)\omega_{\theta}-(b \nabla)\omega_{\theta}+2\frac{v_{\theta}}
{r}\p_{x_3} v_{\theta}+ \omega_{\theta} \frac{v_r}{r}-\p_t
\omega_{\theta}=0,\quad  D_m \times (0, T]\\
\o_\th=0,\,  \text{on} \quad \p D_m \times (0, T];\\
\o_\th(x, 0)= (\o_0)_\th(x), \quad x \in D_m.
\end{cases}
\ee The procedure  is similar to that of Lemma \ref{levthinft} with two differences. One is that the boundary value of $\o_\th$ is $0$ which makes the proof shorter in this aspect. The other, however,  is the presence of the extra nonhomogeneous term $2\frac{v_{\theta}}
{r}\p_{x_3} v_{\theta}$ which needs to be dealt with.

For any rational number $p>1$ in the form of $2k/l$ where $k, l$ are positive integers, we know that $ f=\o_\th^p$ is a sub-solution to the above equation, namely
\be
\lab{subsothp}
\begin{cases}
 (\Delta-\frac{1}{r^2}  \big
)f-(b \nabla)f + p  \frac{v_r}{r} f + 2p\frac{v_{\theta}}
{r}\p_{x_3} v_\th   \, \o^{p-1}_\th  -\p_t
f \ge 0,\quad  D_m \times (0, T]\\
f=0,\,  \text{on} \quad \p D_m \times (0, T];\\
f(x, 0)= (\o_0)^p_\th(x), \quad x \in D_m.
\end{cases}
\ee
As in the previous lemma, for $T \ge 1$, let $\eta=\eta(t)$ be a cut off function in time, supported in $[T-1, T]$ which will be specified later. Since  a nonuniform $L^q$ bound with any $q>1$ for $\o_\th$  (c.f. \erf{olp}) and actually $L^\infty$ bound is already proved in Section 2, we can use $f \eta^2$ as a test function on \erf{subsothp} to deduce
\be
\lab{egothp1}
\al
&LS \equiv \int^T_{T-1}\int_{D_m} |\na (f  \eta)|^2 dxdt + \int^T_{T-1}\int_{D_m} \frac{1}{r^2} f^2 \eta^2 (x, t) dxdt + \frac{1}{2} \int_{D_m} f^2 \eta^2 (x, T) dx\\
& \le
- \underbrace{\int^T_{T-1}\int_{D_m} b \na f \,  f  \eta^2 dxdt}_{T_1}
+ \underbrace{p \int^T_{T-1}\int_{D_m} \frac{v_r}{r} f^2 \eta^2 dxdt}_{T_2}\\
&\qquad +  \underbrace{2 p \int^T_{T-1}\int_{D_m}  \frac{v_{\theta}}
{r}\p_{x_3} v_\th   \, \o^{p-1}_\th f \eta^2 dxdt}_{T_3} +
\int^T_{T-1}\int_{D_m}  f^2 \eta' \eta dxdt + \frac{1}{2} \int_{D_m} f^2 \eta^2 (x, T-1) dx,
\eal
\ee  As in the previous lemma, $T_1=0$ by integration by parts. The term $T_2$ can be estimated exactly
as the term $R_3$ in the previous lemma, although the function $f$ is different. So, like \erf{vthr3}, we have
 \be
\lab{otht2}
|T_2| \le  \frac{1}{8} \int^T_{T-1} \int_{D_m} \left( |\na f |^2 + \frac{1}{r^2} f^2 \right) \eta^2 dxdt +
  (8 C_s)^{13/5} p^{18/5} \Vert \frac{v_r}{r} \Vert_{L^\infty_tL^6(D_m)}^{18/5} \, \int^T_{T-1} \int_{D_m} \,  f^2 \eta^2dxdt.
\ee We remark that the power $18/5$ can be improved since the standard Sobolev inequality holds for $f$ due to the $0$ boundary condition. Since a lower power is not needed we will not pursue it here.

So we are left to bound $T_3$ which requires the $L^\infty_tL^2_x$ bound for $J=-\p_{x_3} v_\th/r$ in
Lemma \ref{leotr/r} and the $L^\infty$ bound for $v_\th$ in Lemma \ref{levthinft}.
Since $| \o^{p-1}_\th f| \le (|\o^p_\th| +1) f =(f+1) f $, we see that
\be
\lab{ot3}
|T_3 | \le 2 p \Vert v_\th \Vert_{L^\infty(Q_1)} \int^T_{T-1} \int_{D_m} |J| f (f+1) dx \eta^2 dt,
\ee where $Q_1=D_m \times [T-1, T]$. For the spatial integral on the right hand side, we can proceed as follows. By H\"older inequality, for $f=f(\cdot, t)$, $J=J(\cdot, t)$ and $\e \in (0, 1)$,
\[
\al
\int_{D_m} |J| (f+1)^2 dx &= \int_{D_m} |J| (f+1)^{3/2} (f+1)^{1/2}  dx\\
&\le \left( \int_{D_m} |J|^2 dx \right)^{1/2} \left( \e \int_{D_m} (f+1)^6 dx \right)^{1/4}
\left( \e^{-1} \int_{D_m} (f+1)^2 dx \right)^{1/4}\\
&\le \e^{1/3} \left( \int_{D_m} (f+1)^6 dx \right)^{1/3} + \e^{-1} \Vert J(\cdot, t) \Vert^4_{L^2_x}
\int_{D_m} (f+1)^2 dx\\
&\le 2 \e^{1/3} \left( \int_{D_m} f^6 dx \right)^{1/3} + [\e^{-1} \Vert J(\cdot, t) \Vert^4_{L^2_x}
+2] \int_{D_m} (f+1)^2 dx.
\eal
\]Here we have used $|D_m|>1/2$. Substituting this to \erf{ot3} and using the standard 3 dimensional Sobolev inequality for functions with $0$ boundary, we deduce
\[
\al
|T_3|  &\le 4  \e^{1/3} p \Vert v_\th \Vert_{L^\infty(Q_1)} \mathbb{S}^2_0 \int^T_{T-1}  \int_{D_m} |\na f|^2 \eta^2 dx dt \\
&\qquad + 2p \Vert v_\th \Vert_{L^\infty(Q_1)} [\e^{-1} \Vert J \Vert^4_{L^\infty_t L^2_x}
+2]  \int^T_{T-1}\int_{D_m} (f+1)^2 \eta^2 dxdt.
\eal
\]Here $\mathbb{S}_0$ is the $L^2$ Sobolev constant in $\mathbb{R}^3$. This infers, after taking $\e$ such that $4  \e^{1/3} p \Vert v_\th \Vert_{L^\infty(Q_1)} \mathbb{S}_0 = 1/4$, that
\be
\lab{otht32}
\al
|T_3|  &\le \frac{1}{4} \int^T_{T-1} \int_{D_m} |\na f|^2 \eta^2 dx dt \\
&\qquad +  4 (16\mathbb{S}^2_0)^3   \Vert v_\th \Vert^4_{L^\infty(Q_1)} [\Vert J \Vert^4_{L^\infty_t L^2_x}
+1] p^4 \int^T_{T-1}\int_{D_m} (f+1)^2 \eta^2  dxdt.
\eal
\ee After substituting \erf{otht2} and \erf{otht32} into \erf{egothp1}, we arrive at

\be
\lab{egothp2}
\al
&\frac{1}{2} \int^T_{T-1}\int_{D_m} \left[ |\na (f  \eta)|^2 dxdt + \frac{1}{r^2} f^2 \eta^2 \right] dxdt + \frac{1}{2} \int_{D_m} f^2 \eta^2 (x, T) dx\\
& \le \Lambda \, p^4 \int^T_{T-1} \int_{D_m} \,  (f+1)^2 \eta^2dxdt +
\int^T_{T-1}\int_{D_m}  f^2 \eta' \eta dxdt + \frac{1}{2} \int_{D_m} f^2 \eta^2 (x, T-1) dx,
\eal
\ee where
\be
\lab{deLamb}
\Lambda \equiv \left[
(8 C_s)^{13/5}  \Vert \frac{v_r}{r} \Vert_{L^\infty_{T-1<t<T}L^6(D_m)}^{18/5} + 4 (16\mathbb{S}^2_0)^3   \Vert v_\th \Vert^4_{L^\infty(D_m \times [T-1, T])} (\Vert \frac{\p_{x_3} v_\th}{r} \Vert^4_{L^\infty_{T-1<t<T} L^2_x}
+1)  \right].
\ee

We treat the case $T=1$ first. Taking $\eta=1$ in \erf{egothp2}, we obtain.
\[
\al
\int^1_0\int_{D_m}& \left[ |\na f |^2 dxdt + \frac{1}{r^2} f^2 \right] dxdt +  \int_{D_m} f^2 (x, T) dx\\
& \le 2 \Lambda \, p^4 \int^1_0 \int_{D_m} \,  (f+1)^2 dxdt +
 \int_{D_m} f^2(x, 0) dx,
\eal
\]which implies, via the standard Sobolev inequality and interpolation that
\[
\left(\int^1_0\int_{D_m} \o_\th^{p 10/3} dxdt \right)^{3/5} \le \mathbb{S}^{6/5}_0  2 \Lambda \, p^4 \int^1_0 \int_{D_m} \,  (|\o_\th|^p+1+\beta_0)^{2} dxdt,
\] where $\beta_0 \equiv \sup_{D_m} |(\o_0)_\th(x)|$. Since $|D_m|>1/2$, this implies, for some absolute constant $C$,
\[
\left(\int^1_0\int_{D_m} [|\o_\th| \vee (1 + \beta_0) ]^{2p (5/3)} dxdt \right)^{3/5} \le C \Lambda \, p^4 \int^1_0 \int_{D_m} \,  [|\o_\th| \vee (1 + \beta_0) ]^{2p} dxdt.
\]Here $|\o_\th| \vee \beta_0 =\max \{|\o_\th|(x, t), \beta_0 \}$. From this inequality, by  Moser's iteration with $p=(5/3)^j$, $j=0, 1, 2, ...$, we deduce
\be
\lab{otht01}
\Vert \o_\th \Vert_{L^\infty(D_m \times [0, 1])} \le  C \Lambda^{5/4} \left(\int^1_0\int_{D_m} (|\o_\th|+\beta_0+1)^2 dxdt\right)^{1/2},
\ee where $C$ is another absolute constant.
 For $T>1$, we take suitable sequences of $\eta$ such that $\eta(T-1)=0$. Then in \erf{egothp2} implies, since the last term  drops out,
 \[
 \left(\int^T_{T-1}\int_{D_m} (|\o_\th| \vee 1)^{2p (5/3)} \eta^2 dxdt \right)^{3/5} \le C  (\Lambda+\sup |\eta'|) \, p^4 \int^T_{T-1} \int_{D_m} \,  (|\o_\th| \vee 1)^{2p} \eta dxdt,
 \]where is an absolute constant depending only on $\mathbb{S}_0$. Moser's iteration again tells us,
\be
\lab{otht>1}
\al
\Vert \o_\th \Vert_{L^\infty(D_m \times [T-0.5, T])} &\le  C (\Lambda+1)^{5/4} \left(\int^T_{T-0.75}\int_{D_m} (|\o_\th|+1)^2 dxdt\right)^{1/2},
\eal
\ee where $C$ may have changed in value but is still an absolute constant. Since $\Lambda$ is given by
\erf{deLamb}, the lemma follows from \erf{otht>1}and \erf{otht01} together with the energy inequality in Lemma \ref{endm} (c).
\qed

Now we are in a position to prove $L^\infty$ bounds for $v_r$ and $v_3$.

\begin{lemma}
\lab{levr3inft}
Let $v$ be a solution to Problem \ref{wtdm} coming from Proposition \ref{prsolS}. Then, there is an absolute constant $C$ such that
\be
\lab{vrinft}
\Vert v_r(\cdot, t) \Vert_{L^\infty(D_m)} \le C C^{13/6}_s \Vert \o_\th(\cdot, t) \Vert^{13/3}_{L^\infty(D_m)}
\left(\int_{D_m}  (|v| + 1)^{2}(\cdot, 0) dx \right)^{1/2},
\ee
\be
\lab{v3inft}
\Vert v_3(\cdot, t) \Vert_{L^\infty(D_m)} \le C \mathbf{s}^{3/2}_0 \Vert \o_\th(\cdot, t) \Vert^{3/2}_{L^\infty(D_m)}
\left(\int_{D_m}  (|v| + 1)^{2}(\cdot, 0) dx \right)^{1/2}.
\ee Here $C_s $ and $\mathbf{s}_0$ are the Sobolev constants in Lemma \ref{lewsob} and Lemma \ref{sobdm} respectively.
\end{lemma}
\proof According to the Biot-Savart law and the Navier boundary condition, the following equations hold for each $t>0$.
\be
\lab{equrt0}
\begin{cases}
\left( \Delta-\frac{1}{r^2} \right) v_r =  \p_{x_3} \o_\th, \quad \text{in} \quad D_m,\\
\p_{x_3} v_r = 0, \quad \text{on} \quad \p^H D_m; \quad v_r = 0, \quad \text{on} \quad \p^V D_m,
\end{cases}
\ee
\be
\lab{equ3t0}
\begin{cases}
\Delta v_3 =  - \left(\p_r \o_\th + \frac{\o_\th}{r}\right), \quad \text{in} \quad D_m,\\
 v_3 = 0, \quad \text{on} \quad \p^H D_m; \quad \p_r  v_3 = 0, \quad \text{on} \quad \p^V D_m.
\end{cases}
\ee where the $t$ variable is omitted for simplicity.

First let us bound $v_r$. By the a priori estimate in Proposition \ref{prsolS}, for any $p \ge 0$, the
function $v^{2p+1}_r$ can be used as a test function of \erf{equrt0}, which infers, after using the boundary condition,
\[
\int_{D_m} |\na v^{p+1}_r|^2 dx + \frac{(p+1)^2}{2p+1} \int_{D_m} \frac{1}{r^2} v^{2(p+1)}_r dx
= (p+1) \int_{D_m} \o_\th v^p_r \p_{x_3} v^{p+1}_r dx.
\]By Cauchy Schwarz, we deduce
\[
\int_{D_m} [ |\na v^{p+1}_r|^2 + \frac{1}{r^2} v^{2(p+1)}_r] dx \le 2 (p+1)^2 \int_{D_m} \o^2_\th v^{2p}_r dx.
\]From the Sobolev inequality in Lemma \ref{lewsob}, this shows,
\[
\left(\int_{D_m} |v_r|^{(p+1) 13/5} dx \right)^{10/13} \le 2 (p+1)^2 C_s \Vert \o_\th \Vert^2_{L^\infty(D_m)}
\int_{D_m}  |v_r|^{2p} dx.
\]
Hence, for some absolute constant $C$, we have
\[
\left(\int_{D_m} (|v_r| \vee 1)^{(p+1) 13/5} dx \right)^{10/13} \le C (p+1)^2 C_s \Vert \o_\th \Vert^2_{L^\infty(D_m)}
\int_{D_m}  (|v_r| \vee 1)^{2(p+1)} dx.
\]Taking $p+1= (13/10)^j$, $j=0, 1, 2, ...$ in Moser's iteration, we conclude, after inserting the $t$ variable and using the energy inequality, that
\[
\al
\Vert v_r(\cdot, t) \Vert_{L^\infty(D_m)} &\le C C^{13/6}_s \Vert \o_\th(\cdot, t) \Vert^{13/3}_{L^\infty(D_m)}
\left(\int_{D_m}  (|v_r| + 1)^{2}(\cdot, t) dx \right)^{1/2}\\
&\le C C^{13/6}_s \Vert \o_\th(\cdot, t) \Vert^{13/3}_{L^\infty(D_m)}
\left(\int_{D_m}  (|v| + 1)^{2}(\cdot, 0) dx \right)^{1/2},
\eal
\] which is \erf{vrinft}.

To bound $v_3$, we use $v^{2p+1}_3$ as a test function on equation \erf{equ3t0} to obtain
\[
\int_{D_m} |\na v^{p+1}_3|^2 dx
= (p+1) \int_{D_m} \o_\th v^p_3 \p_r v^{p+1}_3 dx.
\]
In the above we have used the boundary condition for $v_r$ so that
boundary integrals vanish.
By Cauchy Schwarz, we deduce
\[
\int_{D_m}  |\na v^{p+1}_3|^2 dx \le 2 (p+1)^2 \int_{D_m} \o^2_\th v^{2p}_3 dx.
\]Since $v_3=0$ on the horizontal boundary $\p^H D_m$, the Sobolev inequality in Lemma \ref{sobdm} can be used for $v^{p+1}_3$ in the above identity. Therefore
\[
\left(\int_{D_m} |v_3|^{6(p+1)} dx \right)^{1/3} \le 2 \mathbf{s}^2_0 (p+1)^2  \Vert \o_\th \Vert^2_{L^\infty(D_m)}
\int_{D_m}  |v_3|^{2p} dx.
\]Hence, for some absolute constant $C$, we have
\[
\left(\int_{D_m} (|v_3| \vee 1)^{6(p+1)} dx \right)^{1/3} \le C \mathbf{s}^2_0 (p+1)^2  \Vert \o_\th \Vert^2_{L^\infty(D_m)}
\int_{D_m}  (|v_3| \vee 1)^{2(p+1)} dx.
\]Taking $p+1= 3^j$, $j=0, 1, 2, ...$ in Moser's iteration, we conclude, after inserting the $t$ variable, that
\[
\Vert v_3(\cdot, t) \Vert_{L^\infty(D_m)} \le C \mathbf{s}_0^{3/2} \Vert \o_\th(\cdot, t) \Vert^{3/2}_{L^\infty(D_m)}
\left(\int_{D_m}  (|v_3| + 1)^{2}(\cdot, t) dx \right)^{1/2},
\] which implies \erf{v3inft} by the energy inequality. \qed

\subsection{Completion of the proof}

Pick any $v_0 \in C^2_{nb}(\overline{D_*})$. By definition, there exists a sequence of $v_0^{\{m\}} \in C^2(\overline{D_m})$ satisfying the Navier boundary condition such that $v_0^{\{m\}}$ converge to $v_0$ in $C^2$ norm. According to Proposition
\ref{prsolS}, Lemma \ref{levthinft} and Lemma \ref{levr3inft}, equation \erf{eqasns} with initial value
$v_0^{\{m\}}$ and Navier slip boundary condition on $D_m$ has a finite energy solution $v^{\{m\}}$ such that $\Vert v^{\{m\}}(\cdot, t) \Vert_{L^\infty(D_m)}$ are uniformly bounded for each $t \ge 0$ and all $m$. Also one can see from Lemma \ref{lemuneng} that the energy norm of $v^{\{m\}}$ is uniformly bounded for any finite time. This assertion is clear for $\na v^{\{m\}}_r$, $\p_{x_3} v^{\{m\}}_\theta$ and $\na v^{\{m\}}_3$ from statement (c) in that lemma. To see that $\p_r v^{\{m\}}_\theta$ is in $L^2_tL^2_x$ uniformly, we argue as follows.  From statement (b) in Lemma \ref{lemuneng}, we know that
\[
 \int^T_0 \int_{D_m} \left(|\p_r v^{\{m\}}_\th|^2 +  \frac{1}{r^2} |v^{\{m\}}_\th |^2 \right) dxdt
 - 2  \int^T_0 \int_{D_m} \p_r v^{\{m\}}_\th   \frac{v^{\{m\}}_\th}{r} dxdt
 \le \frac{1}{2} \int_{D_m} |v^{\{m\}}(x, 0)|^2 dx.
\]After integration by parts in the polar coordinates, we infer
\[
\al
 \int^T_0 \int_{D_m} &\left(|\p_r v^{\{m\}}_\th|^2 +  \frac{1}{r^2} |v^{\{m\}}_\th |^2 \right) dxdt
 + \sum^{m+1}_{j=2} \int^T_0 \int_{L_j} |v^{\{m\}}_\th|^2  dx_3dt\\
 &\le \frac{1}{2} \int_{D_m} |v^{\{m\}}(x, 0)|^2 dx + \int^T_0 \int_{L_1} |v^{\{m\}}_\th|^2  dx_3dt.
\eal
\]Again, $L_1$ is the right most vertical boundary of $D_m$ and $L_2, ..., L_{m+1}$ are the connected segments of the left boundary of $D_m$. Since, for each fixed $T$, the right hand side of the above inequality is uniformly bounded due to Lemma \ref{levth1}, we see that $\p_r  v^{\{m\}}_\th$ is uniformly bounded in $L^2_tL^2_x$ norm for each fixed $T>0$. This proves the assertion.
 Therefore we can extract a subsequence, still denoted by $\{ v^{\{m\}} \}$ which converges weakly in the energy norm to a solution $v$ in $\mathbf{E} \cap L^\infty$. Due to standard regularity theory, the convergence can also be made point wise except at the corners of the boundary of $D_*$. This concludes the proof of the theorem.
\qed

\section{Uniqueness results for elliptic equations in polygons with low integrability}

In this section we state and prove the uniqueness results for elliptic equations on domains with non-convex corners, which were used in Section 2. The study of elliptic equations in rough domains has been  very active over the years.  Comparing with standard theory for elliptic equations on Lipschitz domains, the integrability assumption is lower than usual, in view of the non-convex corners with angle $3 \pi/2$. We feel that these results may be known, but are unable to find them in the literature.

\begin{proposition}
\lab{prunid} Let an axially symmetric function $u \in W^{1, s}_0(D_m)$, $s=3^-/2$, be a solution to the Dirichlet problem
\be
\lab{dduvu}
\begin{cases}
\Delta u(x) - V(x) u(x) = 0, \quad x \in D_m,\\
u(x)=0, \quad x \in \p D_m;
\end{cases}
\ee where $V$ is a given axially symmetric function such that $V \ge 1/( 4 r^2)$ and $V \in L^\infty(D_m)$. Then $u=0$.
\end{proposition}

Notice that we are not making any assumption on the non-tangential maximal function of $|\na u|$. Therefore the established theory of elliptic equations in Lipschitz domains  does not seem to apply directly. The idea of the proof is to regard the equation as a two dimensional one in the $r, x_3$ variable first. Next  the problem is converted via conformal mapping to the case with smooth boundary but weaker integrability conditions on the solution and its gradient. Then a duality argument will infer uniqueness.
One can also use the regularity result for very weak solutions in \cite{Sha} after the conversion, even replacing the condition $V \ge 1/(4 r^2)$ by $V \ge 0$. See the proof of Proposition \ref{prunidn} below.

Before proving the proposition, we need the following lemma on uniqueness of $L^{2+}$ solutions on $C^3$ domains in $\mathbb{R}^2$, which may be of independent interest.

\begin{lemma}
\lab{leuni2d}
Let $U \subset \mathbb{R}^2$ be a bounded $C^3$ domain and $Q =\{x_1,, ..., x_k \}$ be finitely many points on the
boundary $\p U$. Given any fixed $p>2$, let
$u \in L^p(U)  \cap C^3(\bar{U} \backslash Q)$
be a solution to
\be
\lab{uv0}
\begin{cases}
\Delta u(x) - V(x) u(x) = 0, \quad x \in U,\\
u(x)=0, \quad x \in \p U \backslash Q,
\end{cases}
\ee where $0 \le V(x) \le \sum^k_{i=1} \frac{C}{|x-x_i|}$ for all $x \in U$ and a constant $C$.  Then $u=0$.
\end{lemma}

{\remark Comparing with standard results, the solution $u$ is not assumed to be bounded near $x_i$, $i=1, 2, ..., k$, nor is there any assumption on $|\na u|$. So the result is more akin to removable singularity theorems, but  with a boundary condition.}

\proof (of Lemma \ref{leuni2d}). Let $u$ be a solution to \erf{uv0}. We first find a solution $\phi \in W^{1,2}_0(U)$ to the standard nonhomogeneous Dirichlet problem
\be
\begin{cases}
\Delta \phi(x) - V(x) \phi(x) = u(x), \quad x \in U,\\
\phi(x)=0, \quad x \in \p U.
\end{cases}
\ee Since $V \ge 0$, $V \in L^{2^-}(U)$ and $u \in L^p(U)$ with $p>2>1$, which is $1/2$ of the dimension,
the potential function $V$ and nonhomogeneous term $u$ are in the regularity class. Therefore
the problem admits a unique solution in $W^{1, 2}_0(U) \cap C^\alpha(U)$ by standard theory. Here $\alpha \in (0, 1)$.

 Now we show that
\be
\lab{dphijie}
|\na \phi | \in L^\infty(U).
\ee It is clear that we just need to prove this in small balls around the potential singular points $x_1, ..., x_k, y_1, ..., y_l$.

Let $\Gamma=\Gamma(x, y)$ be the Dirichlet Green's function on $U$. Then
\be
\lab{phi=gu}
\phi(x)= - \int_U \Gamma(x, y) u(y) dy - \int_U \Gamma(x, y) V(y) \phi(y) dy.
\ee By standard theory, there exists a constant $\bar{C_0}$ depending only on $U$ such that
\[
|\Gamma(x, y) | \le \bar{C_0} [|\ln |x-y|| + 1], \quad |\na_x \Gamma(x, y) | \le \bar{C_0}/|x-y|, \quad
x, y \in U.
\]From these bounds, $\phi \in W^{1, 2}_0(U) \cap C^\alpha(U) \subset L^q(U)$, $\forall q>1$, the bound on $V$, and $u \in L^p(U)$, $p>2$, it is easy to see by \erf{phi=gu} that
\[
\al
|\na \phi(x) | &\le \bar{C_0} \int_U \frac{|u(y)|}{|x-y|} dy + C \bar{C_0}   \int_U \frac{1}{|x-y|} \sum^k_{i=1} \frac{|\phi(y)-\phi(x_i)|}{|x-x_i|}
dy\\
&\le \bar{C_0} \int_U \frac{|u(y)|}{|x-y|} dy + C \bar{C_0}   \int_U \frac{1}{|x-y|} \sum^k_{i=1} \frac{1}{|y-x_i|^{1-\alpha}} dy \, \Vert \phi \Vert_{C^\alpha(U)}
\\
&\le  C \bar{C_0} \left( \Vert u \Vert_{L^p(U)} +  \Vert \phi \Vert_{W^{1, 2}(U)} \right),
\eal
\] proving \erf{dphijie}. In the last step, we have used the standard H\"older estimate for $\phi$ using
its $W^{1, 2}$ norm and $L^{1+}$ norm of $V$.

Pick a small number $\e>0$, denote $U_\e= U \backslash \cup^k_{i=1} B(x_i, \e) .$ Since all possible singularities of $u$ are outside of $U_\e$, we can compute:
\[
\al
\int_{U_\e} u^2 dx &= \int_{U_\e} u(x) \left( \Delta  \phi(x) - V(x) \phi(x)\right) dx\\
&=\int_{\p U_\e} u  \p_n \phi dS - \int_{\p U_\e} \phi  \p_n u dS +
\int_{U_\e} \left( \Delta u(x) - V(x) u(x) \right) \phi(x) dx,
\eal
\]where $n$ is the outward normal of $\p U_\e$. Due to the boundary condition, the above identity becomes
\be
\lab{u2de}
\al
\int_{U_\e} u^2 dx &= -\sum^k_{i=1} \int_{\p B(x_i, \e) \cap U} u  \p_n \phi  dS +
\sum^k_{i=1} \int_{\p B(x_i, \e) \cap U} \phi  \p_n u  dS.
\eal
\ee where $n$ is the outward normal of $\p B(x_i, \e)$ or $\p B(y_j, \e)$. Next we argue that the right hand side of the preceding identity goes to $0$ as $\e \to 0$.

Pick any $x \in \p B(x_i, \e)$.   Then either $B(x, \e/2)$ intersects with $U^c$, then we make $0$ extension for $u$ across $B(x, \e/2) \cap \p U$, or  $B(x, \e/2) \subset U$. We can choose $\e$ sufficiently small so that $B(x, \e/2) \cap \{x_1, x_2, ..., x_k \}$ is empty. So
in either case we can regard $u^2$ as a subsolution, namely
\[
\Delta u^2 \ge 2 V u^2 \ge 0, \quad \text{in} \quad B(x, \e/2).
\]Hence the standard mean value inequality states, for some absolute constant $\bar{C}_1$,   that
\[
|u(x)| \le \bar{C}_1 \left( |B(x, \e/4)|^{-1} \int_{B(x, \e/4)} u^2(y) dy \right)^{1/2}.
\] By the assumption that $u \in L^p(U)$, $p>2$, this implies that
\be
\lab{uo-e}
u(x) = o(1)/\e \quad \text{as} \quad  \e \to 0.
\ee

Further more, when $\e$ is sufficiently small, the region $[B(x_i, 1.5 \e)-B(x_i, 0.5\e)] \cap U$ is free of potential singular points in $Q$. Since $\p U$ is $C^3$ and $V$ has mild singularity near $x_i$, we can obtain a gradient estimate for $u$ on $[B(x_i, 1.2 \e)-B(x_i, 0.8\e)] \cap U$, by straighten out the boundary.
Since $V$ may be unbounded near $y_j$, we give a proof here. Pick  $x \in \p B(x_i, \e)$. For any small $\e$, we choose a smooth cut off function $\eta=\eta(x)$ supported in $B(x, \e/4)$ such that $\eta=1$ on $B(x, \e/8)$ and that
$|\na \eta| \le C/\e$ and $|\Delta \eta | \le C/\e^2$. If $B(x, \e/4)$ straddles the boundary of $U$, we can, as usual, straighten out the boundary and make odd reflection for $u$ and even reflection for $V$ afterwards. So we can just assume $u$ is a classical solution in the full ball $B(x, \e/4)$.
By \erf{uv0}, $u\eta$ satisfies
\[
\Delta (u \eta) - V u \eta = 2 \na \eta \na u + u \Delta \eta.
\]Let $\Gamma_0$ be the Green's function of the Laplacian in $\mathbb{R}^2$, we have, $\forall z \in \mathbb{R}^2$,
\[
\al
u \eta (z) &= - \int \Gamma_0(z, y) V \, u \, \eta (y) dy - 2 \int  \Gamma_0(z, y) \na \eta \, \na u(y) dy -
\int \Gamma_0(z, y) \Delta \eta \, u(y) dy \\
 &= - \int \Gamma_0(z, y) V \, u \, \eta (y) dy + 2 \int \na_y \Gamma_0(z, y) \na \eta \, u(y) dy +
\int \Gamma_0(z, y) \Delta \eta \, u(y) dy.
\eal
\]By our assumptions,  $0 \le V \le C/\e$ in $B(x, \e/4)$ and $\eta=1$ in $B(x, \e/8)$. Thus we deduce, after using  the usual bounds for $|\na \Gamma_0|$ and $|\na^2 \Gamma_0|$, that
\[
\al
|\na u(x) | &\le \frac{C}{\e} \int_{B(x, \e/4)} \frac{ |u(y)|}{|x-y|} dy +
\frac{C}{\e} \int_{B(x, \e/4)\backslash B(x, \e/8)} \frac{ |u(y)|}{|x-y|^2}  dy
+\frac{C}{\e^2} \int_{B(x, \e/4)\backslash B(x, \e/8)} \frac{ |u(y)|}{|x-y|} dy\\
&\le \frac{C}{\e} \int_{B(x, \e/4)} \frac{ |u(y)|}{|x-y|} dy +
\frac{C}{\e^3} \int_{B(x, \e/4)}|u(y)|  dy
\eal
\]This shows, since $u \in L^p$, $p>2$, that for each $x \in \p B(x_i, \e)$, we have
\be
\lab{duo-e2}
|\na u(x)| \le \bar{C}_1 \e^{-1} \sup_{B(x, \e/4)} |u(y)| = o(1)/\e^2.
\ee For the same $x$, by the gradient bound \erf{dphijie} and the boundary condition that $\phi=0$ on $\p U$, we know that
\be
\lab{phiOe}
|\phi(x)| = O(1) \e.
\ee
Substituting  \erf{uo-e}, \erf{duo-e2} and  \erf{phiOe} to the right side of \erf{u2de}, we conclude that $u=0$ after letting $\e \to 0$. This proves the lemma.
\qed

Now we are in a position to give a
\medskip

\noindent {\it Proof of Proposition \ref{prunid}.}

The idea of the proof is to regard the equation in \erf{dduvu} as a two dimensional elliptic equation in the $rx_3$ plane, and the domain $D_m$ as a polygon. Then we use the Schwarz-Christoffel mapping to conformally transform $D_m$ into the upper half plane and eventually to the unit disk. Problem \ref{dduvu} is thus converted to an elliptic problem on the unit disk. The condition that the solution $u \in W^{1, s}_0(D_m)$ will imply that the transformed solution, say $f$,  is in $L^{2+}$ on the unit disk.
Here and later $2^+$ is a number close to and greater than $2$. Then we can conclude from Lemma \ref{leuni2d} that $f$  is $0$. Thus the original solution is also $0$.

Since there are only finitely many corners on $D_m$ and the Schwarz-Christoffel mapping is smooth except at the corners, we can just prove that the transformed solution $f$ is in $L^{2+}$ near the images of these corners, which are points on the boundary of the upper half plane. The worst corners are those with $3\pi/2$ angles. Since the corners are isolated, without loss of generality, we can just carry out the proof of $L^{2+}$ integrability of $f$ near the image of one of the corners. After a change of variable,
the Schwarz-Christoffel mapping is essentially a power function of the complex variable and we can carry out the calculation explicitly as follows.

In $rx_3$ coordinate, the equation for $u$ in \erf{dduvu} is
\[
\p^2_r u + \frac{1}{r} \p_r u + \p^2_{x_3} u - V u =0.
\]In the identity,
\[
\p^2_r ( r^\beta u) = r^\beta \p^2_r u + 2 \beta r^{\beta -1} \p_r u + \beta (\beta-1) r^{\beta-2} u,
\]we take the constant $\beta = 1/2$. Then the above equation for $u$ is transformed into
\be
\lab{eqrru}
\p^2_r ( \sqrt{r} \, u) + \p^2_{x_3} ( \sqrt{r} \, u) -  (V - \frac{1}{4 r^2}) ( \sqrt{r} \, u) =0.
\ee Recall that in $D_m$, the variable $r$ is bounded away from $0$. So there is no singularity in the coefficients of this equation.

{\it Case 1}.  We  work near one of the nonconvex corners, located say at $(r^0, x^0_3)$.

 Let us  make a translation
$r-r^0$, $x_3-x^0_3$ and rotation of angle $\pi$ clockwise, namely
\be
\lab{rx3til}
\tilde{r} = - (r-r^0), \quad \tilde{x_3} = - (x_3-x^0_3).
\ee Then, in the $\tilde{r}\tilde{x_3}$ plane,  we can assume that the corner is located at the origin and the boundary of $D_m$ near the corner coincides with the positive horizontal  ($ \tilde{r} $) axis  and negative vertical ($ \tilde{x_3} $) axis respectively. Consider the functions
\be
\lab{fvrru}
f=f(\tilde{r}, \tilde{x_3})=  \sqrt{r} \, u(r, x_3) = \sqrt{r^0-\tilde{r}} \,  u(r^0-\tilde{r}, x^0_3-\tilde{x_3}), \quad \tilde{V}=\tilde{V}(\tilde{r}, \tilde{x_3})= V(r, x_3)-\frac{1}{4 r^2}.
\ee Then in the $\tilde{r} \tilde{x_3}$ variables, $f$ satisfies
\be
\lab{eqftilrx3}
\p^2_{\tilde{r}} f + \p^2_{\tilde{x_3}} f -  \tilde{V} f=0, \qquad f=0, \quad \text{on} \quad \p \tilde{D_m}.
\ee Here $\tilde{D_m}$ is the region $D_m$ in the $\tilde{r}\tilde{x_3}$ plane, which is still a rectangular polygon. Let us introduce the complex variable
\be
\lab{z=tilrx3}
z= \tilde{r} + i \tilde{x_3}.
\ee In terms of $z$, the equation for $f$ becomes, since
$\frac{\p}{\p z} = \frac{1}{2} \left(\frac{\p}{\p \tilde{r} } - i \frac{\p}{\p \tilde{x_3} } \right)$,
\be
\lab{eqfz}
\frac{\p^2 f}{\p z \p \bar{z}} -  \frac{\tilde{V}}{4} f=0.
\ee Now let $w = \mathbb{M}(z)$ be a Schwarz-Christoffel mapping that maps $\tilde{D_m}$ to a unit circle
 $\mathbb{D}_1$ in the $w=w_1 + i w_2$ plane, whose center is at $(w_1, w_2)=(0,  1)$ and that the corner point $z=0$ is mapped to the point $w=0$. Near the corner $z=0$, say in the region $B(0, \e) \cap D_m$, the mapping is given by, for a constant $\sigma$,
\be
\lab{w=z23}
\al
&w= \sigma \, z^{2/3}(1 + \text{lower order terms}),\\
&\frac{\p w}{\p z} = \frac{2 \sigma}{3} \, z^{-1/3}(1 + \text{lower order terms}),\\
&\frac{\p z}{\p w} = \frac{3 \sigma^{-3/2}}{2} \, w^{1/2}(1 + \text{lower order terms})=
\frac{3 \sigma^{-3/2}}{2} \, w^{1/2}(1 + a_1 w  + a_2 w^2+...).
\eal
\ee Note $|\frac{\p z}{\p w}|$ is bounded when $z$ is away from the corners since $z$ is an analytic function of $w$ and vice versa.

 Let
\be
\lab{domUbb}
U = \mathbb{M}(B(0, \e) \cap \tilde{D_m})
\ee where the ball $B(0, \e)$ in the $z$ plane does not contain any other corners with $\e$ sufficiently small. So $U$ is the image under $\mathbb{M}$ of a small corner in $\tilde{D_m}$. Consider the function
 \[
 g=g(w)= f(z),  \quad w \in \mathbb{D}_1,   \quad w=\mathbb{M}(z).
 \]Then, for any $p \in (2, 4)$, we have, from \erf{w=z23}
\[
\al
\int_U |g|^p dwd\bar{w}&= \frac{4}{9} \sigma^2 \int_{B(0, \e) \cap \tilde{D_m}} |f|^p \, |z|^{-2/3} (1+ o(1)) dzd\bar{z}\\
&=\frac{4}{9} \sigma^2  \left(\int_{B(0, \e) \cap \tilde{D_m}} |f|^{6^-} dzd\bar{z} \right)^{2/3^-}
 \, \left(\int_{B(0, \e) \cap \tilde{D_m}} |z|^{-2^-} (1+ o(1)) dzd\bar{z} \right)^{1/3^-}\\
& \le C \Vert f \Vert^p_{W^{1, s}(\tilde{D_m})}=C \Vert u \Vert^p_{W^{1, s}(D_m)}.
\eal
\]Here we just used the 2 dimensional Sobolev imbedding and the assumption that $s=3^-/2$.
Since there are only finitely many corners, this tells us that
\[
\int_{\mathbb{N}_\e \cap \mathbb{D}_1} |g|^p dwd\bar{w} \le  C \Vert u \Vert^p_{W^{1, s}(D_m)}<\infty,
\]where  $ \mathbb{N}_\e$ is the $\e$ neighborhood of images of nonconvex corners in the $w$ plane.
From the equation for $f$ \erf{eqfz} and the conformal nature of $w = \mathbb{M}(z)$, we know that
$g$ satisfies the equation
\[
\begin{cases}
\frac{\p^2 g}{\p w \p \bar{w}} -  \left|\frac{\p z}{\p w} \right|^2 \frac{\tilde{V}}{4} g=0, \quad w \in \overline{\mathbb{D}_1} \backslash \text{images of corners},\\
 \quad g(w)=0, \quad w \in \p \mathbb{D}_1 \backslash \text{images of corners}.
\end{cases}
\]By \erf{w=z23} the potential function in front of $g$ is bounded and nonnegative near the images of the nonconvex corners.

{\it Case 2}.  We  work near one of the convex corners, located say at $(r^0, x^0_3)$.

After similar translation, we have the following formula for the Schwarz-Christoffel mapping near these corners.
\be
\lab{w=z2390}
\al
&w= \sigma \, z^2(1 + \text{lower order terms}),\\
&\frac{\p w}{\p z} = 2 \, z(1 + \text{lower order terms}),\\
&\frac{\p z}{\p w} = \frac{1}{2 \sigma} \, w^{-1/2}(1 + \text{lower order terms})=
\frac{1}{2 \sigma} \, w^{-1/2}(1 + a_1 w  + a_2 w^2+...).
\eal
\ee  Again, near the images of the convex corners in the $w$ plane, it is easy to check that that $g \in L^p$ with some $p>2$ and that
\[
\begin{cases}
\frac{\p^2 g}{\p w \p \bar{w}} -  \left|\frac{\p z}{\p w} \right|^2 \frac{\tilde{V}}{4} g=0, \quad w \in \overline{\mathbb{D}_1} \backslash \text{images of corners},\\
 \quad g(w)=0, \quad w \in \p \mathbb{D}_1 \backslash \text{images of corners}.
\end{cases}
\]In fact, by standard theory, $u$ is smooth near the convex (i.e. $\pi/2$) corners of $D_m$, which infers that $g$ is H\"older continuous near the image of the convex corners in the $w$ plane. Recall that $w=0$ is actually the image of one of the corners in the $w$ plane.
By \erf{w=z2390} for images of convex corners and Case 1 for those of nonconvex corners, we know that the potential function $\left|\frac{\p z}{\p w} \right|^2 \frac{\tilde{V}}{4}$ satisfies
\be
\lab{vxxi1}
0 \le \left|\frac{\p z}{\p w} \right|^2 \frac{\tilde{V}}{4}(y) \le C \sum^k_{x=1} \frac{1}{|y-x_i|}, \qquad \forall y \in \mathbb{D}_1.
\ee Here $x_i$ are the images of corners in the $w$ plane.  Recall that near the images of  nonconvex corners, by Case 1, we know that the potential function is bounded. Since both $u$ and the Schwarz-Christoffel mapping are smooth except at the corners, we see that
\[
g \in L^p(\mathbb{D}_1) \cap C^3(\mathbb{D}_1 \backslash \text{images of  corners }).
\]
Now we can  use Lemma \ref{leuni2d} on $g$ to conclude that $g=0$ and hence $f=0$ and $u$=0, completing the proof of the proposition. In that lemma, we take $x_i$ as the images of the corners and the potential function as $\left|\frac{\p z}{\p w} \right|^2 \frac{\tilde{V}}{4}$ which satisfies \erf{vxxi1}.
\qed

Similarly we have the following uniqueness result for a mixed Dirichlet-Neumann problem, which was also used in Section 2.

\begin{proposition}
\lab{prunidn} Let an axially symmetric function $u \in W^{1, s}(D_m)$, $s=3^-/2$, be a solution to the mixed Dirichlet-Neumann problem
\be
\lab{ddubvu}
\begin{cases}
\Delta u(x)  - V(x) u(x) = 0, \quad x \in D_m,\\
u(x)=0, \quad x \in \p^H D_m, \quad \p_n u(x) =0, \quad x \in \p^V D_m.
\end{cases}
\ee Here  $V \ge 0$ is a $L^\infty(\bar{D_m})$ axially symmetric  function. Then $u=0$.
\end{proposition}

\proof The main task is to prove that $u \in C(\bar{D_m})$. Taking this statement for granted, one can quickly prove the uniqueness as follows.  If $u \neq 0$, we can assume without loss of generality that $\sup u>0$. Since the potential $V \ge 0$, by the maximum principle, $\sup u$ must occur at $\p D_m$. By the boundary condition $u=0$ on $\p^H D_m$, $\sup u$ must occur at a point, say  $x_0 \in \p^V D_m$ . This point $x_0$ can not be a corner point due to the continuity of $u$ in $\bar{D_m}$. So $\p_n u (x_0)=0$.  But this contradicts with the Hopf maximum principal.

Now let us prove $u \in C(\bar{D_m})$. It is clear we only need to prove $u$ is continuous in a
neighborhood of non-convex corners with angle $3 \pi/2$. In other places, standard reflection method can show that $u \in C^1$.  The proof is similar to that of Proposition \ref{dduvu}. We use a Schwarz-Christoffel mapping to convert the problem to a mixed  Dirichlet-Neumann problem on a smooth domain: the unit disk, which was studied long time ago. The transformed solution will inherit enough integrability property to allow us to use an earlier result to finish the proof.

As in the proof of the previous proposition, we may just work near one of the corners located, say, at $(r^0, x^0_3)$. Using exactly the same notations as before, we see the function $g=g(w)=f(z) = r^{1/2} u$
satisfies
\be
\lab{gwwbar}
\begin{cases}
\frac{\p^2 g}{\p w \p \bar{w}} -  \left|\frac{\p z}{\p w} \right|^2 \frac{\tilde{V}}{4} g=0, \quad w \in \mathbb{D}_1,\\
g(w)=0, \quad \text{or} \quad \p_n g(w)=0, \quad w \in \p \mathbb{D}_1 \backslash \text{images of corners}.
\end{cases}
\ee Due to the boundary condition $u=0$ on $\p^H D_m$ and the assumption that $u \in W^{1, s}$, $s=3^-/2$, we still can apply the 2 dimensional Sobolev inequality to obtain $u \in L^{6^-}(D_m)$. Hence the same proof as before shows that $g \in L^{4^-}(\mathbb{D}_1)$. Next we prove that $|\na g| \in L^{4^-/3}(\mathbb{N}_\epsilon)$, where $ \mathbb{N}_\e$ is again the $\e$ neighborhood of images of nonconvex corners in the $w$ plane..

Let $U$ be the domain in \erf{domUbb}. Then, for $q=4^-/3$, we compute, from \erf{w=z23},
\[
\al
&\int_U \left|\frac{\p g}{\p w} \right|^q dwd\bar{w}= \frac{4}{9} \sigma^{2-q} \int_{B(0, \e) \cap \tilde{D_m}} \left|\frac{\p f}{\p z} \right|^q \, \left(\frac{3}{2}\right)^q |z|^{q/3} \, |z|^{-2/3} (1+ o(1)) dzd\bar{z}\\
&= \left( 2/3 \right)^{2-q} \sigma^{2-q}  \left(\int_{B(0, \e) \cap \tilde{D_m}}
\left| \frac{\p f}{\p z} \right| ^{3^-/2} dzd\bar{z} \right)^{2q/3^-}
 \, \left(\int_{B(0, \e) \cap \tilde{D_m}} |z|^{\frac{3^-(q-2)}{3 (3^{-}-2q)}} (1+ o(1)) dzd\bar{z} \right)^{(3^{-}-2q)/3^-}\\
& \le C \Vert f \Vert^q_{W^{1, s}(\tilde{D_m})} \le C \Vert u \Vert^q_{W^{1, s}(D_m)},
\eal
\] where we have used $\frac{3^-(q-2)}{3 (3^{-}-2q)}>-2$ when $3^-$ is sufficiently close to $3$. Likewise the same bound holds for $|\frac{\p g}{\p \bar{w}}|$. Hence $g |_{\mathbb{N}_\e} \in W^{1, 4^-/3}(\mathbb{N}_\e) \subset L^{4^-}(\mathbb{N}_\e)$, where the inclusion is due to the 2 dimensional Sobolev inequality.
Now we can apply the main result in Shamir \cite{Sha} on \erf{gwwbar} to conclude that $g \in C^{1^-/2}(\bar{\mathbb{D}_1})$. The detail is as follows.
Notice that in small neighborhoods around the images of the nonconvex corners, we have $|\frac{\p z}{\p w}|^2 = C |w| (1 + \text{lower order terms})$ by \erf{w=z23}. So the coefficient for $g$ in \erf{gwwbar} is Lipschitz in these neighborhoods. By bootstrapping from $g|_{\mathbb{N}_\e} \in W^{1, 4^-/3}( \mathbb{N}_\e) \subset L^{4^-}(\mathbb{N}_\e)$, we can treat the term $ \left|\frac{\p z}{\p w} \right|^2 \frac{\tilde{V}}{4} g |_{\mathbb{N}_\e}  \in L^{4^-}(\mathbb{N}_\e)$ as an inhomogeneous term and apply local versions of Lemma 5.1 or Theorem 5.3 there with large $p$ to deduce $g|_{\mathbb{N}_\e} \in W^{1, 4^-}(\mathbb{N}_\e)$ and so on, giving us $g|_{\mathbb{N}_\e} \in C^{1^-/2}(\mathbb{N}_\e)$.  As mentioned, outside of $\mathbb{N}_\e$, we already know $u$ is in $C^1$. Therefore $u \in C(\bar{D_m})$ which completes the proof of the proposition.
\qed

\section*{Acknowledgments} We wish to thank  Professors M. Benes,  Hongjie Dong,  Zhen Lei, Zijin Li,
  Xinghong Pan, Xin Yang and Na Zhao and Mr. Chulan Zeng for helpful discussions. The support of Simons Foundation grant 710364 is gratefully acknowledged.
  \medskip
  
\noindent{\it{Competing interests statement.}} The author has no competing interests to declare.

\bibliographystyle{plain}



\def\cprime{$'$}

\end{document}